\documentclass[a4paper,10pt,fleqn]{article}
\usepackage{amsmath,amssymb,amsthm,graphicx,subfigure,float,caption,epstopdf,tabularx,color, bm,amsfonts,epic}
\usepackage[top=1.25in, bottom=1.0in, left=1.0in, right=1.0in]{geometry}
\usepackage{appendix}
\usepackage{amssymb}
\usepackage{multirow}
\usepackage{mathrsfs}
\usepackage[table]{xcolor}
\usepackage[numbers,sort&compress]{natbib}
\usepackage{comment,enumerate,multicol,xspace}
\usepackage{fancyhdr}
\usepackage{soul}
\usepackage{enumerate}
\usepackage{makecell}
\definecolor{Cobalt}{rgb}{0.25,0.41,0.88}
\usepackage[linkcolor=Cobalt,anchorcolor=Cobalt,citecolor=Cobalt,colorlinks]{hyperref}
\usepackage{bbm}
\usepackage{algorithm}
\usepackage{algorithmic}
\usepackage{booktabs}
\usepackage{lineno}
\usepackage{dsfont}
\usepackage{cases}
\usepackage{pifont}
\usepackage{empheq}

\allowdisplaybreaks
\newtheorem{theorem}{Theorem}[section]

\newtheorem{lemma}[theorem]{Lemma}

\numberwithin{equation}{section}

\usepackage{indentfirst}
\graphicspath{{Fig}}
\begin{document}
\title{Galerkin finite element method for time fractional diffusion-wave equations with constant time delay}
\author{Xuan Li$^{1}$,\quad Kaile Qiu$^{1}$,\quad Chaobao Huang$^2$,\quad Weiping Bu$^{1,3}$\footnote{Correspondence author. Email: weipingbu@xtu.edu.cn}
\\{\small $^1$ School of Mathematics and Computational Science, Xiangtan University, }
\\{\small Hunan 411105, China}
\\{\small$^2$ School of Statistics and Mathematics, Shandong University of Finance and Economics,}
\\{\small Jinan 250014, China}
\\{\small $^3$ Hunan Key Laboratory for Computation and Simulation in Science and Engineering,}
\\{\small Hunan 411105, China}
}
\date{}
\maketitle
\begin{abstract}
	This paper studies numerical methods for the time fractional diffusion-wave equations with the constant delay. Using a symmetric fractional-order reduction technique and introducing the auxiliary variable $w$, we reformulate the original problem as an equivalent coupled system and for the first time carry out a rigorous regularity analysis of $w$. Based on the symmetric graded time mesh, the coupled problem is discretised by the $L1$ scheme in time and a standard finite element method in space, and the stability result of the developed fully discrete finite element scheme is established by means of a discrete fractional Gr\"onwall inequality.
    Furthermore, we derive a new comparison principle and reconstruct a new barrier function on each interval $((k-1)\tau,k\tau]$, $1\le k\le K$  according to the properties of the symmetric graded mesh. Based on the new barrier function, we obtain new relevant properties, which lead to a novel and crucial stability property for the discrete fractional operator. Combining the obtained local truncation error in time gives a local error estimate for the fully discrete finite element scheme, which shows that our convergence result provides a unified framework that seamlessly encompasses both pointwise-in-time error estimates on uniform meshes and global-in-time error estimates on symmetrical graded meshes.
    Numerical experiments confirm the theoretical results.\\[2ex]
\end{abstract}
\section{Introduction}
As is known to all, there exist many physical and engineering systems
exhibiting long-range temporal memory
and spatial interactions. For these systems,
classical differential equations often fail to
accurately capture such long-range dependencies,
while fractional calculus
owing to their non-locality
provides a natural description of these  problems \cite{Hilfer2000Applications,Sun2018A,Kilbas2006Theory}.
Furthermore,
when the evolution of these systems also inherently exhibits time delays,
fractional delay differential equations provide a more precise description of their dynamics, and have been widely applied in some fields such as biology, automatic control systems with feedback, economics, and medicine \cite{Rihan2021Delay,Schiesser2020Time,Smith20101An,Erneux2009Applied,Driver2012Ordinary}
In this paper,
we consider the following time fractional diffusion-wave equation
with constant time delay
\begin{equation}\label{eq1_1}
{}^{C}_{0}D_{t}^{\alpha} u(x,t) = p \Delta u(x,t) + a u(x,t) + b u(x,t-\tau) + f(x,t), \quad 0 < x < L, \; 0 < t \le K\tau,
\end{equation}
subject to the initial and boundary conditions
\begin{equation}\label{eq1_2}
u(x,t) = \phi(x,t),\quad u_t(x,t) = \psi(x,t), \quad 0 \leq x \leq L, \quad-\tau \leq t \leq 0,
\end{equation}
%\begin{equation}\label{eq1_3}
%u_t(x,t) = \psi(x,t), \quad 0 \leq x \leq L, \quad -\tau \le t \leq  0
%\end{equation}
\begin{equation}\label{eq1_4}
u(0,t) = 0, \quad u(L,t) = 0, \quad -\tau \leq t \leq K\tau,
\end{equation}
where \(\alpha\in(1,2)\), \(p\) is a positive constant,
\(a \le 0\), \(\tau > 0\) is the time delay parameter,  \(K\) is a given positive integer,
 \(f(x,t)\), \(\phi(x,t)\) and \(\psi(x,t)\) are some continuous functions.
Here, the time Caputo fractional derivative is defined as
\begin{equation}\label{eq1_5}
{}^{C}_{0}D_{t}^{\gamma} u(x,t) = \int_{0}^{t} \omega_{n-\gamma}(t-s) \partial_s^{n} u(x,s)\, ds, \quad n-1 < \gamma < n, \quad n\in \mathbb{N},
\end{equation}
where \(\omega_{\rho}(t) := t^{\rho-1} / \Gamma(\rho)\) for \(\rho > 0\), $\partial_s^n(\cdot)$ denotes the $n-$th order derivative with respect to the variable $s$, and $\partial_s(\cdot):=\partial_s^1(\cdot)$.

At present, there have existed some theoretical studies on time fractional delay diffusion-wave equations such as those in \cite{Yang2023Stability,Li2024Existence}. However, the difficulty of the practical applicability of analytical solutions has motivated a strong focus on numerical techniques. In \cite{Zhang2022Compact}, a compact finite difference method was developed for a one-dimensional variable-coefficient time fractional delay diffusion-wave equation. For the two-dimensional case, an alternating direction implicit method was proposed in \cite{Zhang2023Numerical} to enhance the computational efficiency. By combining the Crank-Nicolson method and the trigonometric collocation method, an efficient numerical method was proposed in \cite{Chawla2025A} for a fractional delay diffusion-wave equation. In Refs. \cite{Hendy2017semi,Abbaszadeh2021Galerkin}, when the delay diffusion-wave equations with distributed-order time-fractional derivative are considered, the linearized compact finite difference method and the Galerkin meshless method were proposed. Furthermore, Ref. \cite{Alikhanov2024second} discussed a delay diffusion-wave equation with generalized time-fractional operator, and developed a second-order finite difference method.
In addition, error estimates for time-fractional subdiffusion problems also have drawn considerable attention. Li et al. \cite{Li2021sharp} established sharp pointwise-in-time error estimates for nonlinear subdiffusion equations. Ref. \cite{Peng2024Unconditionally} obtained the unconditional convergence and superconvergence of a Newton linearized finite element scheme for nonlinear time-fractional parabolic equations with distributed delay. Jiang and Chen \cite{Jiang2024Local} analyzed the global and local convergence of two-dimensional time-fractional reaction-subdiffusion equations with weakly singular solutions.

In recent research on some time fractional delay differential equations, Refs. \cite{Tan2022L1,Cen2023The} discovered that the exact solutions exhibited multiple singularities in time. Subsequently, Refs. \cite{Ou2025Modified,Tan2026Order} extended the discussion to time fractional diffusion wave equation with with constant time delay, and found that the exact solution satisfied
\begin{equation}\label{eq1_6}
	|\partial_{t}^{l}u_{k\tau}(x,t)|\leq C(1+(t-(k-1)\tau)^{k\alpha-l}), \quad l=0,1,2,
\end{equation}
where $u_{k\tau}(x,t):=u(x,t),t\in((k-1)\tau,k\tau].$ This results implies that the exact solution has potential multiple temporal singularities for $l>2.$ From the above discussion, it is obvious that numerical treatments of time fractional delay diffusion wave equation based on  (\ref{eq1_6}) remain very limited. The studies mentioned in Refs. \cite{Ou2025Modified,Tan2026Order} represent pioneering efforts in this respect. Under the uniform time mesh, the former adopts the convolution quadrature method to discretize the fractional derivative, whereas the latter treats the fractional derivative via an order reduction strategy combined with the $L2-1_\sigma$ formula.

In our view, while uniform time mesh offers the benefit of avoiding round-off errors caused by overly concentrated local mesh points, non-uniform mesh provides a compelling advantage by adapting to the low regularity of the solution. However, to the best of our knowledge, existing convergence results for numerical methods for time fractional delay differential equations are confined to either local-in-time convergence on uniform meshes or global-in-time convergence on non-uniform meshes. No work has yet considered the local error estimate on non-uniform meshes.
In \cite{Tan2022L1}, a symmetric graded time mesh is proposed and successfully applied to handle the multiple temporal singularities of the solution.
Motivated by the above discussion, we develop a fully discrete finite element scheme based on the symmetric graded time mesh to solve problem (\ref{eq1_1})--(\ref{eq1_4}). The main contributions of this paper include:

$\bullet$ By employing the order reduction technique mentioned in \cite{Lyu2022A} and introducing a intermediate variable $w = D_t^{\alpha/2}(u - t\psi(x,0)) $, we reformulate the original problem (\ref{eq1_1})--(\ref{eq1_4}) into an equivalent coupled system. Unlike in Ref. \cite{Tan2026Order}, we first derive the regularity of $w$, which is essential for the subsequent local truncation error analysis in time.
%for which the regularity of the intermediate variable $w$ is analyzed in detail.

$\bullet$ By discretizing the resulting coupled system using the L1 formula on the symmetric graded time mesh in the temporal direction and the finite element method in the spatial direction, we propose a fully discrete finite element scheme for an equivalent coupled system. Moreover, using the discrete fractional Gr\"onwall inequality, we confirm that the proposed scheme is stable.

$\bullet$ To the best of our knowledge, we are the first to develop the work of \cite{Kopteva2025Pointwise} to the delay case. Different from \cite{Kopteva2025Pointwise}, we derive a new comparison principle. In addition, we reconstruct a new barrier function in each interval $((k-1)\tau,k\tau]$, $1\le k\le K$ according to the properties of the symmetric graded mesh. Based on the new barrier function, we obtain new relevant properties, which lead to a novel and crucial stability property for the discrete fractional operator. By using the above obtained result, we establish a local temporal error estimate for the proposed numerical scheme on symmetric graded time mesh.
%To the best of our knowledge, existing convergence results for numerical methods for time fractional delay differential equations are confined to either local-in-time convergence on uniform meshes or global-in-time convergence on non-uniform meshes.
Observing that the symmetric graded time mesh can reduce to a uniform one when the mesh parameter is appropriately chosen, we discover that our convergence result provides a unified framework that seamlessly encompasses both pointwise-in-time error estimates on uniform meshes and global-in-time error estimates on symmetrical graded meshes.

The remainder of this paper is organized as follows.
In Section 2, some notations are introduced, an equivalent coupled system is constructed by introducing an intermediate variable $w,$ and the regularity of $w$ is analyzed.
In Section 3, a fully discrete finite element method for the coupled problem is proposed, and by using the discrete fractional Gr\"onwall inequality, a new stability result for the proposed numerical scheme is presented.
In Section 4, we establish a new stability property of the discrete fractional derivative operator in the delayed case. The local truncation error of the L1 scheme on the symmetric graded time mesh is derived, and the local temporal error estimate for the proposed numerical scheme is obtained.
In Section 5, the obtained theoretical results are verified by some numerical experiments. Finally, a conclusion is give in Section 6.

\textbf{Notation.} In this paper, $C$ denotes a positive constant, whose value may differ in different cases. Denote \(y \lesssim z\) to mean \(y \leq Cz\), and \(y \simeq z\) to mean \(y \lesssim z\) and \(z \lesssim y\).

\section{Preliminaries}
In this section, we first introduce some notations.
Define the inner product $(\cdot, \cdot)$ and  norm $\|\cdot\|$ of  $L^2([0,L])$ by
$(u,v) = \int_{0}^{L} uvdx$, $\|v\| =  \left(v,v\right)^{1/2}, u,v \in L^{2}([0,L]).$
Let $G(u)= -(p\Delta u(x,t) + a u(x,t))$, and $X_i$ be  a eigenfunction corresponding to  eigenvalue $\lambda_i$  satisfying
%\begin{equation}\label{eq2_3}
$G(X_{i}) = \lambda_{i} X_{i}$,
$X_{i}(0) =  X_{i}(L) = 0$ and $\|X_i\|=1$,
where $i=1,2,\cdots.$ For $\gamma>0$, denote
$\|g\|_{G^{\gamma}}=\left(\sum_{i=1}^{\infty}\lambda_{i}^{2\gamma} \left|(g, X_{i}(x))\right|^{2}\right)^{1/2}
$
and
$
D(G^{\gamma})=\left\{g\in L^{2}(0,L):\sum_{i=1}^{\infty}\lambda_{i}^{2\gamma}|(g,X_{i}(x))|^{2}<\infty\right\}.
$

In order to construct finite element scheme for \eqref{eq1_1}--\eqref{eq1_4}, we first transform the considered problem into an equivalent coupled system by using the symmetric fractional-order reduction technique proposed in \cite{Lyu2022A}.
Let $\tilde{u}(x,t) = u(x, t) - t\psi(x,0)$. Then \eqref{eq1_1}--\eqref{eq1_4} can be rewritten as
	\begin{align}
		&{}^{C}_{0}D_{t}^{\alpha} \tilde{u}(x,t)= p\Delta \tilde{u}(x,t) + a\tilde{u}(x, t) + b \tilde{u}(x, t-\tau) + g(x, t),  \quad 0<x<L,\quad 0<t<K \tau ,\label{eq2_36a1}\\
		&\tilde{u}(x, t) = \phi(x, t) - t\psi(x,0),\quad \tilde{u}_t(x, t) =\psi(x,t) -\psi(x,0), \quad 0 \leq x \leq L, \quad -\tau \leq t \leq 0 ,\label{eq2_36a2}\\
		%&\tilde{u}_t(x, t) =\psi(x,t) -\psi(x,0), \quad 0 \leq x \leq L ,\quad -\tau \leq t \leq 0,\label{eq2_36a3}\\
		&\tilde{u}(0,t) = 0, \quad \tilde{u}(L, t) = 0, \quad -\tau \leq t \leq K\tau,\label{eq2_36a4}
	\end{align}
where
\begin{equation}\label{eq2_37}
	g(x, t) = pt\Delta\psi(x,0)+at \psi(x,0) + b(t - \tau) \psi(x,0)  + f(x, t).
\end{equation}
Denote $\beta=\frac{\alpha}{2}$ and $w(x,t)= {}^{C}_{0}D_{t}^{\beta} \tilde{u}(x,t)$. As proved in \cite[Lemma 2.1]{Lyu2022A}, one has
\begin{equation}\label{eq2_36a}
{}^{C}_{0}D_{t}^{\alpha} \tilde{u}(x,t) = {}^{C}_{0}D_{t}^{\beta} w(x,t).
\end{equation}
Therefore, it follows from  \eqref{eq2_36a} that the problem \eqref{eq2_36a1}--\eqref{eq2_36a4} can be reformulate as the following coupled system
	\begin{align}
	&{}^{C}_{0}D_{t}^{\beta}w(x,t) = p\Delta \tilde{u}(x,t) + a\tilde{u}(x, t) + b \tilde{u}(x, t-\tau) + g(x, t), \quad 0<x<L,\quad 0<t<K \tau , \label{eq2_36b1}\\
	&{}^{C}_{0}D_{t}^{\beta}\tilde{u}(x,t) - w(x,t) = 0,  \quad 0<x<L,\quad 0<t<K \tau ,\label{eq2_36b2}\\
	&\tilde{u}(x, t) = \phi(x, t) - t\psi(x,0),\quad w(x, t) = {}^{C}_{0}D_{t}^{\beta}(\phi(x, t) - t\psi(x,0)),\quad 0 \leq x \leq L, \quad -\tau \leq t \leq 0,\label{eq2_36b3}   \\
	%&w(x, t) = {}^{C}_{0}D_{t}^{\beta}(\phi(x, t) - t\psi(x,0)),   \quad 0 \leq x \leq L, \quad -\tau \leq t \leq 0,   \label{eq2_36b4}\\
	&\tilde{u}_t(x, t) = \psi(x,t)-\psi(x,0), \quad 0 \leq x \leq L,\quad -\tau \leq t \leq 0,\label{eq2_36b5}\\
	&w(0, t)= w(L, t) =\tilde{u}(0, t) = \tilde{u}(L, t) = 0, \quad -\tau \leq t \leq K\tau\label{eq2_36b6}.
\end{align}

Similar to \cite[Section 2]{Tan2022L1} and \cite[Section 2.1]{Tan2026Order}, the method of separation of variables implies that, the solution $\tilde{u}_{k\tau}(x, t):=\tilde{u}(x, t), (x,t)\in [0,L]\times((k-1)\tau, k\tau]$ of \eqref{eq2_36a1}--\eqref{eq2_36a4} can be written as
\begin{equation}\label{eq2_40}
	\tilde{u}_{k\tau}(x, t) = \sum_{i=1}^{\infty} \left[\phi_i(0) E_{\alpha,1}(-\lambda_i t^{\alpha})
	+ \int_0^t s^{\alpha-1} E_{\alpha,\alpha}(-\lambda_i s^{\alpha})H_{k\tau,i}(t-s) ds\right]X_i,
\end{equation}
where $E_{\xi_1,\xi_2}(z):= \sum_{j=0}^{\infty} \frac{z^j}{\Gamma(\xi_1 j + \xi_2)}$ is the Mittag-Leffler function, $k=1,2,3,\cdots,K,$ $\phi_i(0)=(\phi(x,0),X_{i}(x)),$ and
$H_{k\tau,i}(t)=(H_{k\tau}(x,t),X_{i}(x))$ with $\tilde{u}_{0}(x,t):=\phi(x,t)-t\psi(x,0), t\in[-\tau,0]$ and
\begin{equation}\label{eq2_41}
	H_{k\tau}(x,t)=
	\begin{cases}
		b\tilde{u}_{0}(x,t-\tau)+g(x,t),&0\leq x \leq L,\quad 0\leq t\leq\tau,\\ b\tilde{u}_{\tau}(x,t-\tau)+g(x,t),&0\leq x\leq L,\quad\tau<t\leq 2\tau,\\
         \vdots \\
        b\tilde{u}_{(k-1)\tau}(x,t-\tau)+g(x,t),&0\leq x\leq L,\quad(k-1)\tau<t\leq k\tau.
	\end{cases}
\end{equation}
From \cite{Tan2026Order}, we can obtain the temporal regularity of the solution $\tilde{u}(x, t)$ of the above coupled system. Therefore, we now discuss the temporal regularity of the introduced variable $w$. First, we give three useful lemmas.
\begin{lemma}\label{lemma2.1}(\cite{Sakamoto2011Initial},\cite[p.83]{Gorenflo2020Mittag-Leffler})
Let $\xi_1>0$, $\xi_2>0.$ For $t\geq 0$ and $\xi_1<2,$ the following boundedness of Mittag-Leffler function holds
\begin{equation}\label{eq2_8}
	E_{\xi_1,\xi_2}(-t)\leq\frac{C}{1+t}.
\end{equation}
What is more, for given positive integer $m$ and $t>0$, one has %$\frac{d^m}{dt^m} E_{\xi_1,1}(-\lambda t^{\xi_1}) = -\lambda t^{\xi_1-m} E_{\xi_1,\xi_1-m+1}(-\lambda t^{\xi_1})$ and $\frac{d^m}{dt^m} \left[ t^{\xi_2-1} E_{\xi_1,\xi_2} \left( -\lambda t^{\xi_1} \right) \right] = t^{\xi_2-m-1} E_{\xi_1,\xi_2-m}(-\lambda t^{\xi_1}).$
\begin{align*}
	&\frac{d^m}{dt^m} E_{\xi_1,1}(-\lambda t^{\xi_1}) = -\lambda t^{\xi_1-m} E_{\xi_1,\xi_1-m+1}(-\lambda t^{\xi_1}),\\
	&\frac{d^m}{dt^m} \left[ t^{\xi_2-1} E_{\xi_1,\xi_2} \left( -\lambda t^{\xi_1} \right) \right] = t^{\xi_2-m-1} E_{\xi_1,\xi_2-m}(-\lambda t^{\xi_1}).
\end{align*}  	
\end{lemma}

\begin{lemma}\label{lemma2.2}
Let $1<\alpha<2.$ Then, for $t\in((j-1)\tau,j\tau], 2\leq j\leq k$,
the following inequality holds
\begin{equation}\label{eq2_13b}
\begin{aligned}
&IV_{\alpha,\theta_1,\theta_2}:=\int_{0}^{\tau}(t-s)^{\alpha-\theta_1}ds+\sum_{j=2}^{k}\int_{(j-1)\tau}^{\min\{j\tau,t\}}(t-s)^{\alpha-\theta_1}\left(1+\sum_{m=2}^{j}(s-(m-1)\tau)^{(m-1)\alpha-\theta_2}\right)ds\\
\leq& C\left(t^{\alpha-\theta_1+1}+\sum_{m=2}^{k}(t-(m-1)\tau)^{m\alpha-(\theta_1+\theta_2)+1}\right),
\end{aligned}
\end{equation}
where $\theta_1=1,2$,
$\theta_2=0,1,2$,
$k=2,3,\dots,K.$
\end{lemma}
\begin{proof}
Since
\begin{align*}
    %&\int_{0}^{\tau}(t-s)^{\alpha-\theta_1}ds+\sum_{j=2}^{k}\int_{(j-1)\tau}^{\min\{j\tau,t\}}(t-s)^{\alpha-\theta_1}\left(1+\sum_{m=2}^{j}(s-(m-1)\tau)^{(m-1)\alpha-\theta_2}\right)ds\\
    IV_{\alpha,\theta_1,\theta_2}=&\int_{0}^{\tau} (t-s)^{\alpha-\theta_1} ds
    + \int_{\tau}^{2\tau} (t-s)^{\alpha-\theta_1} \left(1+(s-\tau)^{\alpha-\theta_2}\right) ds\\
    %+ \int_{2\tau}^{3\tau} (t-s)^{\alpha-\theta_1} \left(1+(s-\tau)^{\alpha-\theta_2}+(s-2\tau)^{2\alpha-\theta_2}\right) ds \\
    &+ \cdots
    + \int_{(k-1)\tau}^{\min\{k\tau,\,t\}}
    (t-s)^{\alpha-\theta_1}
    \left(1+(s-\tau)^{\alpha-\theta_2}+\cdots+\left(s-(k-1)\tau\right)^{(k-1)\alpha-\theta_2}\right) ds \\
    =& \int_{0}^{t} (t-s)^{\alpha-\theta_1} ds
    + \int_{\tau}^{t} (t-s)^{\alpha-\theta_1} (s-\tau)^{\alpha-\theta_2} ds + \cdots
    + \int_{(k-1)\tau}^{t}
    (t-s)^{\alpha-\theta_1} \left(s-(k-1)\tau\right)^{(k-1)\alpha-\theta_2} ds,
\end{align*}
by applying the Euler-beta function, it is easy to obtain
\begin{align*}
    %&\int_{0}^{\tau}(t-s)^{\alpha-\theta_1}ds+\sum_{j=2}^{k}\int_{(j-1)\tau}^{\min\{j\tau,t\}}(t-s)^{\alpha-\theta_1}\left(1+\sum_{m=2}^{j}(s-(m-1)\tau)^{(m-1)\alpha-\theta_2}\right)ds\\
    %\leq& C\left( t^{\alpha-\theta_1+1}
    %+ (t-\tau)^{2\alpha-(\theta_1+\theta_2)+1}
    %+\cdots
    %+\left(t-(k-1)\tau\right)^{k\alpha-(\theta_1+\theta_2)+1}\right)\\
    IV_{\alpha,\theta_1,\theta_2}\leq& C\left(t^{\alpha-\theta_1+1}+\sum_{m=2}^{k}(t-(m-1)\tau)^{m\alpha-(\theta_1+\theta_2)+1}\right).
\end{align*}
\end{proof}
\begin{lemma}\label{lemma2.3}
Let $1<\alpha<2$, $0<\epsilon\ll 1$. Suppose that
\begin{enumerate}
    \item[(I)]
    $\phi(x,t-\tau)\in D(G^{\frac{5+\epsilon}{4}})$,
    $\partial_t\phi(x,t-\tau)\in D(G^{\frac{1+\epsilon}{4}})$,
    $\partial_t^2\phi(x,t-\tau)\in D(G^{\frac{1+\epsilon}{4}})$,
     $\psi(x,t-\tau)\in D(G^{\frac{5+\epsilon}{4}})$ for each $t\in[0,\tau]$,
     and $\psi(0,0)=\psi(L,0)=0$;
    \item[(II)]
    $f(x,t)\in D(G^{\frac{1+\epsilon}{4}})$,
    $\partial_t f(x,t)\in D(G^{\frac{1+\epsilon}{4}})$
    and $\partial_t^2 f(x,t)\in D(G^{\frac{1+\epsilon}{4}})$, for each $t\in[0,K\tau]$.
\end{enumerate}
Then one has
$\sum_{i=1}^{\infty} |H_{\tau,i}(0)|
    \le C,$
$\sum_{i=1}^{\infty} |\partial_t H_{\tau,i}(0)|
    \le C$
and for $t\in((j-1)\tau,j\tau], 1\leq j\leq k$,
\begin{align}
%    &\sum_{i=1}^{\infty} |H_{k\tau,i}(t)|
%    \le C\left(1+\sum_{i=2}^{j}(t-(i-1)\tau)^{(i-1)\alpha}\right),\label{lemma2.3eq1}\\
%   &\sum_{i=1}^\infty |\partial_t H_{k\tau,i}(t)|
%   \le C\left(1+\sum_{i=2}^{j}(t-(i-1)\tau)^{(i-1)\alpha-1}\right),\label{lemma2.3eq2}\\
    \sum_{i=1}^\infty |\partial_t^l H_{k\tau,i}(t)|
    \le C\left(1+\sum_{m=2}^{j}(t-(m-1)\tau)^{(m-1)\alpha-l}\right),\label{lemma2.3eq3}
\end{align}
where $l=0,1,2$ and $k=1,2,\dots,K.$
\end{lemma}
\begin{proof}
When $t\in[0,\tau]$, it follows from \eqref{eq2_37} and \eqref{eq2_41} that
\begin{equation}\label{eq2_14a1}
\begin{aligned}
\left(H_\tau(x,t),X_i(x)\right)
%&=\left(b\phi(x,t-\tau)+pt\Delta\psi(x,0)+at\psi(x,0)+f(x,t),X_i(x)\right)\\
&=\left(b\phi(x,t-\tau),X_i(x)\right)
+\left(pt\Delta\psi(x,0),X_i(x)\right)+\left(at\psi(x,0),X_i(x)\right)
+\left(f(x,t),X_i(x)\right).
\end{aligned}
\end{equation}
According to the assumptions (I), we know that
$\phi(x,t-\tau)\in D(G^{\frac{5+\varepsilon}{4}})\subset D(G^{\frac{1+\varepsilon}{4}})$. Therefore,
applying the Cauchy-Schwarz inequality yields
\begin{equation}\label{eq2_14a4}
\sum_{i=1}^{\infty}\left|\left(\phi(x,t-\tau),X_i(x)\right)\right|
\leq \left(\sum_{i=1}^{\infty}\frac{1}{\lambda_i^{\frac{1+\varepsilon}{2}}}\right)^{\frac12}
\left(\sum_{i=1}^{\infty}\lambda_i^{\frac{1+\varepsilon}{2}}\left|\left(\phi(x,t-\tau),X_i(x)\right)\right|^2\right)^{\frac12}\leq C
\end{equation}
and
\begin{equation}\label{eq2_14a5}
\sum_{i=1}^{\infty}\left|\lambda_i\left(\phi(x,t-\tau),X_i(x)\right)\right|
\leq \left(\sum_{i=1}^{\infty}\frac{1}{\lambda_i^{\frac{1+\varepsilon}{2}}}\right)^{\frac12}
\left(\sum_{i=1}^{\infty}\lambda_i^{\frac{5+\varepsilon}{2}}\left|\left(\phi(x,t-\tau),X_i(x)\right)\right|^2\right)^{\frac12}\leq C,
\end{equation}
where the fact $\lambda_i \approx i^2$ \cite[p. 415]{R1953Methods} is used.
Due to
$\psi(x,0)\in D(G^{\frac{5+\varepsilon}{4}})\subset D(G^{\frac{1+\varepsilon}{4}})$
and $f(x,t)\in D(G^{\frac{1+\varepsilon}{4}})$, an analogous argument leads to
\begin{equation}\label{eq2_14a10}
\sum_{i=1}^{\infty}\left|\lambda_i\left(\psi(x,0),X_i(x)\right)\right|\leq C,
\end{equation}
and
\begin{equation}\label{eq2_14a9}
\sum_{i=1}^{\infty}\left|\left(\psi(x,0),X_i(x)\right)\right|\leq C,\quad
\sum_{i=1}^{\infty}\left|\left(f(x,t),X_i(x)\right)\right|\leq C.
\end{equation}
Moreover, since $\psi(0,0)=\psi(L,0)=0$ and $G(X_i(x))=\lambda_iX_i(x)$
imply that
$
(\Delta\psi(x,0),X_i(x))
=(\psi(x,0),\Delta X_i(x))
=-\frac{a+\lambda_i}{p}(\psi(x,0),X_i(x)),
$
by combining \eqref{eq2_14a10} and \eqref{eq2_14a9}, we have
\begin{equation}\label{eq2_14a11}
\sum_{i=1}^{\infty}\left|\left(\Delta\psi(x,0),X_i(x)\right)\right|
\leq\frac{-a}{p}\sum_{i=1}^{\infty}\left|\left(\psi(x,0),X_i(x)\right)\right|
+\frac{1}{p}\sum_{i=1}^{\infty}\left|\lambda_i\left(\psi(x,0),X_i(x)\right)\right|\leq C.
\end{equation}
Note that $H_{\tau,i}(t)=(H_\tau(x,t),X_i(x)).$ Thus, for $t\in[0,\tau],$ the combination of \eqref{eq2_14a4}, \eqref{eq2_14a9} and \eqref{eq2_14a11} immediately obtain
%\begin{equation}\label{eq2_14a3}
$\sum_{i=1}^{\infty}|H_{\tau,i}(t)|
\leq C.$
%\end{equation}
Now we differentiate $H_\tau(x,t),t\in[0,\tau]$ in \eqref{eq2_41} with respect to $t$, which gives
$
\partial_t H_\tau(x,t)
=b\partial_t \phi(x,t-\tau)+p\Delta\psi(x,0)+a\psi(x,0)+\partial_t f(x,t)
$
and
$
 \partial_t^2 H_\tau(x,t)=b\partial_t^2 \phi(x,t-\tau)+\partial_t^2 f(x,t).
$
As $\partial_t \phi(x,t-\tau)$, $\partial_t f(x,t)\in D(G^{\frac{1+\varepsilon}{4}})$ in (I) and (II)
show that
\begin{equation}\label{eq2_17a2}
\sum_{i=1}^{\infty}\left|\left(\partial_t \phi(x,t-\tau),X_i(x)\right)\right|
\leq C,\quad
\sum_{i=1}^{\infty}\left|\left(\partial_t f(x,t),X_i(x)\right)\right|\leq C,
\end{equation}
it follows from  \eqref{eq2_14a9}, \eqref{eq2_14a11} and \eqref{eq2_17a2} that
\begin{equation}\label{eq2_17a4}
\begin{aligned}
\sum_{i=1}^{\infty}\left|\partial_t H_{\tau,i}(t)\right|
=&\sum_{i=1}^{\infty}\left|(b\partial_t \phi(x,t-\tau)+p\Delta\psi(x,0)+a\psi(x,0)+\partial_t f(x,t),X_i(x))\right|\\
\leq &b\sum_{i=1}^{\infty}\left|\left(\partial_t \phi(x,t-\tau),X_i(x)\right)\right|
+p\sum_{i=1}^{\infty}\left|\left(\Delta\psi(x,0),X_i(x)\right)\right|\\
&+(-a)\sum_{i=1}^{\infty}\left|\left(\psi(x,0),X_i(x)\right)\right|
+\sum_{i=1}^{\infty}\left|\left(\partial_t f(x,t),X_i(x)\right)\right|\leq C.
\end{aligned}
\end{equation}
Similarly, in view of the fact that the assumption $\partial_t^2 \phi(x,t-\tau), \partial_t^2 f(x,t)\in D(G^{\frac{1+\varepsilon}{4}})$ imply
\begin{equation}\label{2_17aa1}
\sum_{i=1}^{\infty}\left|\left(\partial_t^2 \phi(x,t-\tau),X_i(x)\right)\right|
\leq C,\quad
\sum_{i=1}^{\infty}\left|\left(\partial_t^2 f(x,t),X_i(x)\right)\right|\leq C,
\end{equation}
therefore we have
\begin{equation}\label{eq2_25a1}
\sum_{i=1}^{\infty}\left|\partial_t^2 H_{\tau,i}(t)\right|
=\sum_{i=1}^{\infty}\left|\left(b\partial_t^2 \phi(x,t-\tau)+\partial_t^2 f(x,t),X_i(x)\right)\right|\leq C.
\end{equation}
The above discussion demonstrates $\sum_{i=1}^{\infty} |H_{\tau,i}(0)|
    \le C$,
$\sum_{i=1}^{\infty} |\partial_tH_{\tau,i}(0)|
    \le C$
and \eqref{lemma2.3eq3} holds for $[0,\tau]$.

Now we aim to bound $\sum_{i=1}^{\infty} |\partial_t^lH_{2\tau,i}(t)|, l=0,1,2$ for $t\in(\tau,2\tau]$. Due to
\begin{equation}\label{eq2_25aa1}
     \sum_{i=1}^{\infty} |\partial_t^lH_{2\tau,i}(t)|
    =\sum_{i=1}^{\infty} |(b \partial_t^l\tilde{u}_\tau(x,t-\tau)+\partial_t^lg(x,t),X_i(x))|,
\end{equation}
we bound $\sum_{i=1}^{\infty} |( \partial_t^l\tilde{u}_\tau(x,t-\tau),X_i(x))|$ and
$\sum_{i=1}^{\infty} |(\partial_t^lg(x,t),X_i(x))|,$ respectively. For simplicity, let $\tilde{t}=t-\tau$, $A_{i}(\tilde{t})=E_{\alpha,1}(-\lambda_i \tilde{t}^{\alpha})$, $B_i(\tilde{t})= \tilde{t}^{\alpha-1}E_{\alpha,\alpha}(-\lambda_i \tilde{t}^\alpha)$ and $C_i(\tilde{t})=\tilde{t}^{\alpha-2}E_{\alpha,\alpha-1}(-\lambda_i \tilde{t}^\alpha)$ in the following discussion.
First, it follows from \eqref{eq2_40} and Lemma \ref{lemma2.1} that
\begin{equation}\label{eq2_25aa2}
	\tilde{u}_{\tau}(x, \tilde{t}) = \sum_{i=1}^{\infty} \left[\phi_i(0) A_i(\tilde{t})
	+ \int_0^{\tilde{t}} B_i(\tilde{t}-s)H_{\tau,i}(s) ds\right]X_i(x),
\end{equation}
%\begin{equation}\label{eq2_25aa2}
	%\tilde{u}_{\tau}(x, \tilde{t}) = \sum_{i=1}^{\infty} \left[\phi_i(0) E_{\alpha,1}(-\lambda_i \tilde{t}^{\alpha})
	%+ \int_0^{\tilde{t}} (\tilde{t}-s)^{\alpha-1} E_{\alpha,\alpha}(-\lambda_i (\tilde{t}-s)^{\alpha})H_{\tau,i}(s) ds\right]X_i(x),
%\end{equation}
\begin{equation}\label{eq2_17}
		\begin{aligned}
			\partial_{t}\tilde{u}_\tau(x,\tilde{t})
			=&\partial_{t}\sum_{i=1}^{\infty} \left[\phi_i(0)A_i(\tilde{t})
	       + \int_0^{\tilde{t}} B_i(s)H_{\tau,i}(\tilde{t}-s) ds\right]X_i(x)\\
            =&\sum_{i=1}^{\infty}\left[-\phi_i(0)\lambda_iB_i(\tilde{t})
			+H_{\tau,i}(0)B_i(\tilde{t})
			+\int_0^{\tilde{t}} B_i(s)\partial_{t}  H_{\tau,i}(\tilde{t}-s)ds\right]X_i(x)\\
            =&\sum_{i=1}^{\infty}\left[-\phi_{i}(0)\lambda_iB_i(\tilde{t})
			+H_{\tau,i}(0)B_i(\tilde{t})
			+\int_0^{\tilde{t}} B_i(\tilde{t}-s) \partial_{s} H_{\tau,i}(s)ds\right]X_i(x)
		\end{aligned}
	\end{equation}
%\begin{equation}\label{eq2_17}
		%\begin{aligned}
			%\partial_{t}\tilde{u}_\tau(x,\tilde{t})
			%=&\partial_{t}\sum_{i=1}^{\infty} \left[\phi_i(0) E_{\alpha,1}(-\lambda_i \tilde{t}^{\alpha})
	       %+ \int_0^{\tilde{t}} s^{\alpha-1} E_{\alpha,\alpha}(-\lambda_i s^{\alpha})H_{\tau,i}(\tilde{t}-s) ds\right]X_i(x)\\
           % =&\sum_{i=1}^{\infty}\left[-\phi_{i}(0)\lambda_i \tilde{t}^{\alpha-1}E_{\alpha,\alpha}(-\lambda_i \tilde{t}^\alpha)
			%+\tilde{t}^{\alpha-1}E_{\alpha,\alpha}(-\lambda_i \tilde{t}^\alpha)H_{\tau,i}(0)
			%+\int_0^{\tilde{t}} s^{\alpha-1}E_{\alpha,\alpha}(-\lambda_i s^\alpha)\partial_{t} H_{\tau,i}(\tilde{t}-s)ds\right]X_i(x)\\
          %  =&\sum_{i=1}^{\infty}\left[-\phi_{i}(0)\lambda_i \tilde{t}^{\alpha-1}E_{\alpha,\alpha}(-\lambda_i \tilde{t}^\alpha)
			%+\tilde{t}^{\alpha-1}E_{\alpha,\alpha}(-\lambda_i \tilde{t}^\alpha)H_{\tau,i}(0)
			%+\int_0^{\tilde{t}} (\tilde{t}-s)^{\alpha-1}E_{\alpha,\alpha}(-\lambda_i (\tilde{t}-s)^\alpha)\partial_{s} H_{\tau,i}(s)ds\right]X_i(x)
		%\end{aligned}
	%\end{equation}
and
\begin{equation}\label{eq2_22}
		\begin{aligned}
			\partial_{t}^2 \tilde{u}_\tau(x,\tilde{t})
			=&\sum_{i=1}^\infty \left[
			-\phi_i(0)\lambda_i C_i(\tilde{t})
			+C_i(\tilde{t})H_{\tau,i}(0)
		      +B_i(\tilde{t})\partial_{t} H_{\tau,i}(0)\right.
			 \left.+\int_0^{\tilde{t}} B_i(s)\partial_{t}(\partial_s H_{\tau,i}(\tilde{t}-s))\,ds
			\right]X_i(x)\\
			=&\sum_{i=1}^{\infty}
			\left[
			-\phi_{i}(0)\lambda_i C_i(\tilde{t})
			+ C_i(\tilde{t})H_{\tau,i}(0)\right.
			+\left.\int_0^{\tilde{t}}C_i(\tilde{t}-s) \partial_s H_{\tau,i}(s)ds
			\right]X_i(x).
		\end{aligned}
	\end{equation}
%\begin{equation}\label{eq2_22}
		%\begin{aligned}
			%\partial_{t}^2 \tilde{u}_\tau(x,\tilde{t})
			%=&\sum_{i=1}^\infty \left[
			%-\phi_i(0)\lambda_i \tilde{t}^{\alpha-2}E_{\alpha,\alpha-1}(-\lambda_i \tilde{t}^\alpha)
			%+\tilde{t}^{\alpha-2}E_{\alpha,\alpha-1}(-\lambda_i \tilde{t}^\alpha)H_{\tau,i}(0)
			%+\tilde{t}^{\alpha-1}E_{\alpha,\alpha}(-\lambda_i \tilde{t}^\alpha)\partial_{t} H_{\tau,i}(0)\right.\\
			%& \left.+\int_0^{\tilde{t}} s^{\alpha-1}E_{\alpha,\alpha}(-\lambda_i s^\alpha)\partial_{t}(\partial_s H_{\tau,i}(\tilde{t}-s))\,ds
			%\right]X_i(x)\\
			%=&\sum_{i=1}^{\infty}
			%\left[
			%-\phi_{i}(0)\lambda_i \tilde{t}^{\alpha-2}E_{\alpha,\alpha-1}\left(-\lambda_i \tilde{t}^\alpha\right)
			%+ \tilde{t}^{\alpha-2}E_{\alpha,\alpha-1}\left(-\lambda_i \tilde{t}^\alpha\right)H_{\tau,i}(0)\right.\\
			%&+\left.\int_0^{\tilde{t}} (\tilde{t}-s)^{\alpha-2}E_{\alpha,\alpha-1}\left(-\lambda_i (\tilde{t}-s)^\alpha\right)\partial_s H_{\tau,i}(s)ds
			%\right]X_i(x).
		%\end{aligned}
	%\end{equation}
where the integration-by-parts formula is used. Combining $\sum_{i=1}^{\infty}|H_{\tau,i}(t)|
\leq C$, \eqref{eq2_8}, \eqref{eq2_14a4}, \eqref{eq2_14a5}
%\eqref{eq2_14a3}
and \eqref{eq2_17a4}
%$|X_i(x)|\le C$ \cite[p. 335]{Luchko2012Initial-boundary-value}
it is easy to obtain
	\begin{equation}\label{eq2_16}
		\begin{aligned}
         \sum_{i=1}^{\infty} |(\tilde{u}_\tau(x,\tilde{t}),X_i(x))|
			&\le C\left[\sum_{i=1}^{\infty}|\phi_{i}(0)|
			+\int_0^{\tilde{t}}(\tilde{t}-s)^{\alpha-1}\sum_{i=1}^{\infty}\left|H_{\tau,i}(s)\right|ds\right]
			\le C(1+\tilde{t}^{\alpha}),
		\end{aligned}
	\end{equation}
\begin{equation}\label{eq2_21}
		\begin{aligned}
        \sum_{i=1}^{\infty} |( \partial_{t}\tilde{u}_\tau(x,\tilde{t}),X_i(x))|
			&\le C\left[\tilde{t}^{\alpha-1}\sum_{i=1}^{\infty}\left|\lambda_i \phi_{i}(0)\right|
			+  \tilde{t}^{\alpha-1}\sum_{i=1}^{\infty}\left|H_{\tau,i}(0)\right|
			+ \int_0^{\tilde{t}} (\tilde{t}-s)^{\alpha-1}\sum_{i=1}^{\infty}\left|\partial_{s} H_{\tau,i}(s)\right|ds\right]
			\le C(1+\tilde{t}^{\alpha-1})
		\end{aligned}
\end{equation}
and
\begin{equation}\label{eq2_25}
		\begin{aligned}
        \sum_{i=1}^{\infty} |( \partial_{t}^2\tilde{u}_\tau(x,\tilde{t}),X_i(x))|
		&\le C\left[\tilde{t}^{\alpha-2}\sum_{i=1}^\infty |\lambda_i\phi_i(0)|
			+\tilde{t}^{\alpha-2}\sum_{i=1}^\infty |H_{\tau,i}(0)|
			+\int_0^{\tilde{t}} (\tilde{t}-s)^{\alpha-2}\sum_{i=1}^\infty |\partial_{s} H_{\tau,i}(s)|\,ds\right]
			\le C(1+\tilde{t}^{\alpha-2}).
		\end{aligned}
	\end{equation}
Next, we consider the boundedness of  $\sum_{i=1}^{\infty}|(\partial_t^lg(x,t),X_i(x))|$, $l=0,1,2$. Note that
$
g(x,t)= pt\Delta\psi(x,0)+at \psi(x,0) + b(t - \tau) \psi(x,0)  + f(x, t),
$
$
\partial_t g(x,t)=p\Delta\psi(x,0)+a\psi(x,0)+b\psi(x,0)+\partial_t f(x,t)
$
and
$
\partial_t^2 g(x,t)=\partial_t^2 f(x,t).
$
The combination of \eqref{eq2_14a9},\eqref{eq2_14a11}, \eqref{eq2_17a2} and \eqref{2_17aa1} can immediately obtain
\begin{equation}\label{7.13eq1}
 \sum_{i=1}^{\infty} |(\partial_t^lg(x,t),X_i(x))|\leq C, l=0,1,2.
\end{equation}
 Combining \eqref{eq2_25aa1},\eqref{eq2_16},\eqref{eq2_21} and \eqref{eq2_25}, we obtain
\begin{equation}\label{eq2_26a}
 \sum_{i=1}^{\infty} |\partial_t^lH_{2\tau,i}(t)|
		=\sum_{i=1}^{\infty} |(b \partial_t^l\tilde{u}_\tau(x,t-\tau)+\partial_t^lg(x,t),X_i(x))|
		\le C(1+(t-\tau)^{\alpha-l}),\quad l=0,1,2.
\end{equation}
Thus, the above derivation demonstrates that \eqref{lemma2.3eq3} is correct for $t\in(\tau,2\tau].$

Suppose that \eqref{lemma2.3eq3} is true for $k=n$, $1<n<K$, i.e., for $t\in((j-1)\tau,j\tau]$, $j=2,3,\dots,n,$
\begin{equation}\label{7.15eq2}
\begin{aligned}
\sum_{i=1}^{\infty} |\partial_t^lH_{n\tau,i}(t)|
    %&=\sum_{i=1}^{\infty} |(b \partial_t^l\tilde{u}_{(n-1)\tau}(x,t-%\tau)+\partial_t^lg(x,t),X_i(x))|\\
    \le C\left(1+\sum_{m=2}^{j}(t-(m-1)\tau)^{(m-1)\alpha-l}\right), l=0,1,2.
    \end{aligned}
\end{equation}
Then, when $k=n+1$, it is obvious that
\begin{equation}\label{eq2_27}
\sum_{i=1}^{\infty} |\partial_t^lH_{(n+1)\tau,i}(t)|
\le C\left(1+\sum_{m=2}^{j}(t-(m-1)\tau)^{(m-1)\alpha-l}\right),l=0,1,2,
\end{equation}
holds for $t\in((j-1)\tau,j\tau]$, $j=2,3,\dots,n.$ Now we estimate the boundedness of $\sum_{i=1}^{\infty} |\partial_t^lH_{(n+1)\tau,i}(t)|$ under $t\in (n\tau, (n+1)\tau].$ Since, when $t\in (n\tau, (n+1)\tau],$
\begin{equation}\label{7.15eq1}
\begin{aligned}
\sum_{i=1}^{\infty} |\partial_t^lH_{(n+1)\tau,i}(t)|
&=\sum_{i=1}^{\infty} |(b \partial_t^l\tilde{u}_{n\tau}(x,t-\tau)+\partial_t^lg(x,t),X_i(x))|\\
&\le \sum_{i=1}^{\infty} |(b \partial_t^l\tilde{u}_{n\tau}(x,t-\tau),X_i(x))|+\sum_{i=1}^{\infty} |(\partial_t^lg(x,t),X_i(x))|\\
&\le \sum_{i=1}^{\infty} |(b \partial_t^l\tilde{u}_{n\tau}(x,t-\tau),X_i(x))| + C,
\end{aligned}
\end{equation}
we only need to bound $\sum_{i=1}^{\infty} |( \partial_t^l\tilde{u}_{n\tau}(x,t-\tau),X_i(x))|,t\in (n\tau, (n+1)\tau],$
where \eqref{7.13eq1} is used.
From Lemma \ref{lemma2.2} and  \eqref{eq2_27}, an argument similar to that presented in \eqref{eq2_16}--\eqref{eq2_25} suggests that
\begin{equation}\label{eq2_31}
    \begin{aligned}
        \sum_{i=1}^{\infty}|(\tilde{u}_{n\tau}(x,\tilde{t}),X_i(x))|
        =&\sum_{i=1}^{\infty} \left|\phi_i(0) A_i(\tilde{t})+\int_0^\tau B_i(\tilde{t}-s)H_{\tau,i}(s)ds
       +\sum_{j=2}^{n}\int_{(j-1)\tau}^{\min\{j\tau,\tilde{t}\}}B_i(\tilde{t}-s)H_{n\tau,i}(s)\right|ds\\
        \le& C\left(1+IV_{\alpha,1,0}
        %\int_0^\tau (\tilde{t}-s)^{\alpha-1}\,ds+
       % +\sum_{j=2}^{n}\int_{(j-1)\tau}^{\min\{j\tau,\tilde{t}\}}(\tilde{t}-s)^{\alpha-1}
        %\left(1+\sum_{m=2}^{j}(\tilde{t}-(m-1)\tau)^{(m-1)\alpha}\right)ds
        \right)
        \le C\left(1+\tilde{t}^{\alpha} +\sum_{m=2}^{n}(\tilde{t}-(m-1)\tau)^{m\alpha}\right),\\
    \end{aligned}
\end{equation}
%\begin{equation}\label{eq2_31}
    %\begin{aligned}
       % \sum_{i=1}^{\infty}|(\tilde{u}_{n\tau}(x,\tilde{t}),X_i(x))|
        %=&\sum_{i=1}^{\infty} \left|\phi_i(0) E_{\alpha,1}(-\lambda_i \tilde{t}^{\alpha})+\int_0^\tau (\tilde{t}-s)^{\alpha-1}E_{\alpha,\alpha}(-\lambda_i (\tilde{t}-s)^{\alpha})H_{\tau,i}(s)\,ds\right.\\
       %&+\left.\sum_{j=2}^{n}\int_{(j-1)\tau}^{\min\{j\tau,\tilde{t}\}}(\tilde{t}-s)^{\alpha-1}E_{\alpha,\alpha}(-\lambda_i (\tilde{t}-s)^{\alpha})H_{n\tau,i}(s)\right|ds\\
       % \le& C\left(1+\int_0^\tau (\tilde{t}-s)^{\alpha-1}\,ds
       % +\sum_{j=2}^{n}\int_{(j-1)\tau}^{\min\{j\tau,\tilde{t}\}}(\tilde{t}-s)^{\alpha-1}
        %\left(1+\sum_{m=2}^{j}(\tilde{t}-(m-1)\tau)^{(m-1)\alpha}\right)ds\right)\\
        %\le& C\left(1+\tilde{t}^{\alpha} +\sum_{m=2}^{n}(\tilde{t}-(m-1)\tau)^{m\alpha}\right),\\
  %  \end{aligned}
%\end{equation}
\begin{equation}\label{eq2_32a1}
    \begin{aligned}
        \sum_{i=1}^{\infty}|(\partial_{t} \tilde{u}_{n\tau}(x,\tilde{t}),X_i(x))|
       % &=\sum_{i=1}^{\infty}\left|-\phi_{i}(0)\lambda_i \tilde{t}^{\alpha-1}E_{\alpha,\alpha}(-\lambda_i \tilde{t}^\alpha)
	    %+\tilde{t}^{\alpha-1}E_{\alpha,\alpha}(-\lambda_i \tilde{t}^\alpha)H_{n\tau,i}(0)\right.\\
	    %&+\left.\int_0^{\tilde{t}} (\tilde{t}-s)^{\alpha-1}E_{\alpha,\alpha}(-\lambda_i (\tilde{t}-s))^\alpha)\partial_{s} H_{n\tau,i}(s)ds\right|\\
       =&\sum_{i=1}^{\infty}\left|-\phi_{i}(0)\lambda_i B_i(\tilde{t})
     +B_i(\tilde{t})H_{n\tau,i}(0)\right.
	    +\int_0^{\tau} B_i(\tilde{t}-s)\partial_{s} H_{\tau,i}(s)ds\\
        &+\left.\sum_{j=2}^{n} \int_{(j-1)\tau}^{\min\{j\tau,\tilde{t}\}} B_i(\tilde{t}-s)
        \partial_{s} H_{n\tau,i}(s)ds\right|\\
         \leq& C\left( \tilde{t}^{\alpha-1}
        +IV_{\alpha,1,1}
        %\int_{0}^{\tau} (\tilde{t}-s)^{\alpha-1}ds
        %+\sum_{j=2}^{n} \int_{(j-1)\tau}^{\min\{j\tau,\tilde{t}\}} (\tilde{t}-s)^{\alpha-1} \left(1+\sum_{m=2}^{j} (s-(m-1)\tau)^{(m-1)\alpha-1} \right) ds
        \right)
        \leq C\left(\tilde{t}^{\alpha}+\tilde{t}^{\alpha-1}+\sum_{m=2}^{n}(\tilde{t}-(m-1)\tau)^{m\alpha-1}\right)
    \end{aligned}
\end{equation}
%\begin{equation}\label{eq2_32a1}
   % \begin{aligned}
     %   \sum_{i=1}^{\infty}|(\partial_{t} \tilde{u}_{n\tau}(x,\tilde{t}),X_i(x))|
       % &=\sum_{i=1}^{\infty}\left|-\phi_{i}(0)\lambda_i \tilde{t}^{\alpha-1}E_{\alpha,\alpha}(-\lambda_i \tilde{t}^\alpha)
	    %+\tilde{t}^{\alpha-1}E_{\alpha,\alpha}(-\lambda_i \tilde{t}^\alpha)H_{n\tau,i}(0)\right.\\
	    %&+\left.\int_0^{\tilde{t}} (\tilde{t}-s)^{\alpha-1}E_{\alpha,\alpha}(-\lambda_i (\tilde{t}-s))^\alpha)\partial_{s} H_{n\tau,i}(s)ds\right|\\
      % =&\sum_{i=1}^{\infty}\left|-\phi_{i}(0)\lambda_i \tilde{t}^{\alpha-1}E_{\alpha,\alpha}(-\lambda_i \tilde{t}^\alpha)
     %+\tilde{t}^{\alpha-1}E_{\alpha,\alpha}(-\lambda_i \tilde{t}^\alpha)H_{n\tau,i}(0)\right.\\
	    %&+\int_0^{\tau} (\tilde{t}-s)^{\alpha-1}E_{\alpha,\alpha}(-\lambda_i (\tilde{t}-s)^\alpha)\partial_{s} H_{\tau,i}(s)ds\\
       % &+\left.\sum_{j=2}^{n} \int_{(j-1)\tau}^{\min\{j\tau,\tilde{t}\}} (\tilde{t}-s)^{\alpha-1}E_{\alpha,\alpha}(-\lambda_i (\tilde{t}-s))^\alpha)\partial_{s} H_{n\tau,i}(s)ds\right|\\
         %\leq& C\Bigg( \tilde{t}^{\alpha-1}
       % + \int_{0}^{\tau} (\tilde{t}-s)^{\alpha-1}ds
        %+\sum_{j=2}^{n} \int_{(j-1)\tau}^{\min\{j\tau,\tilde{t}\}} (\tilde{t}-s)^{\alpha-1} \left(1+\sum_{m=2}^{j} (s-(m-1)\tau)^{(m-1)\alpha-1} \right) ds\Bigg) \\
       % \leq &C\left(\tilde{t}^{\alpha}+\tilde{t}^{\alpha-1}+\sum_{m=2}^{n}(\tilde{t}-(m-1)\tau)^{m\alpha-1}\right)
   % \end{aligned}
%\end{equation}
and
\begin{equation}\label{eq2_32a3}
    \begin{aligned}
        \sum_{i=1}^{\infty}|(\partial_{t}^2 \tilde{u}_{n\tau}(x,\tilde{t}),X_i(x))|
        %&=\sum_{i=1}^{\infty}
			%\left|
			%-\phi_{i}(0)\lambda_i \tilde{t}^{\alpha-2}E_{\alpha,\alpha-1}\left(-\lambda_i \tilde{t}^\alpha\right)
			%+ \tilde{t}^{\alpha-2}E_{\alpha,\alpha-1}\left(-\lambda_i \tilde{t}^\alpha\right)H_{n\tau,i}(0)\right. \\
			%&\left.+\int_0^{\tilde{t}} (\tilde{t}-s)^{\alpha-2}E_{\alpha,\alpha-1}\big(-\lambda_i (\tilde{t}-s))^\alpha\big)\partial_s H_{n\tau,i}(s)ds\right|\\
         =& \sum_{i=1}^{\infty}
			\left|
			-\phi_{i}(0)\lambda_i C_i(\tilde{t})
			+ C_i(\tilde{t})H_{n\tau,i}(0)\right.
			+\int_0^{\tau} C_i(\tilde{t}-s) \partial_s H_{\tau,i}(s)ds\\
            &+\left.\sum_{j=2}^{n} \int_{(j-1)\tau}^{\min\{j\tau,\tilde{t}\}} C_i(\tilde{t}-s)\partial_s H_{n\tau,i}(s)ds\right|\\
         \leq& C\left( \tilde{t}^{\alpha-2} +
        IV_{\alpha,2,1}
        %\int_{0}^{\tau} (\tilde{t}-s)^{\alpha-2}ds
        %+\sum_{j=2}^{n} \int_{(j-1)\tau}^{\min\{j\tau,\tilde{t}\}} (\tilde{t}-s)^{\alpha-2} \left(1 +\sum_{m=2}^{j} (s-(m-1)\tau)^{(m-1)\alpha-1} \right) ds
        \right)
        \leq C\left(\tilde{t}^{\alpha-1}+\tilde{t}^{\alpha-2}+\sum_{m=2}^{n}(\tilde{t}-(m-1)\tau)^{m\alpha-2}\right),
    \end{aligned}
\end{equation}
%\begin{equation}\label{eq2_32a3}
    %\begin{aligned}
       % \sum_{i=1}^{\infty}|(\partial_{t}^2 \tilde{u}_{n\tau}(x,\tilde{t}),X_i(x))|
        %&=\sum_{i=1}^{\infty}
			%\left|
			%-\phi_{i}(0)\lambda_i \tilde{t}^{\alpha-2}E_{\alpha,\alpha-1}\left(-\lambda_i \tilde{t}^\alpha\right)
			%+ \tilde{t}^{\alpha-2}E_{\alpha,\alpha-1}\left(-\lambda_i \tilde{t}^\alpha\right)H_{n\tau,i}(0)\right. \\
			%&\left.+\int_0^{\tilde{t}} (\tilde{t}-s)^{\alpha-2}E_{\alpha,\alpha-1}\big(-\lambda_i (\tilde{t}-s))^\alpha\big)\partial_s H_{n\tau,i}(s)ds\right|\\
         %&= \sum_{i=1}^{\infty}
			%\left|
			%-\phi_{i}(0)\lambda_i \tilde{t}^{\alpha-2}E_{\alpha,\alpha-1}\left(-\lambda_i \tilde{t}^\alpha\right)
			%+ \tilde{t}^{\alpha-2}E_{\alpha,\alpha-1}\left(-\lambda_i \tilde{t}^\alpha\right)H_{n\tau,i}(0)\right. \\
			%&+\int_0^{\tau} (\tilde{t}-s)^{\alpha-2}E_{\alpha,\alpha-1}\big(-\lambda_i (\tilde{t}-s))^\alpha\big)\partial_s H_{\tau,i}(s)ds\\
           %& \left.+\sum_{j=2}^{n} \int_{(j-1)\tau}^{\min\{j\tau,\tilde{t}\}} (\tilde{t}-s)^{\alpha-2}E_{\alpha,\alpha-1}\big(-\lambda_i (\tilde{t}-s)^\alpha\big)\partial_s H_{n\tau,i}(s)ds\right|\\
        %& \leq C\Bigg( \tilde{t}^{\alpha-2} + \int_{0}^{\tau} (\tilde{t}-s)^{\alpha-2}ds
        %+\sum_{j=2}^{n} \int_{(j-1)\tau}^{\min\{j\tau,\tilde{t}\}} (\tilde{t}-s)^{\alpha-2} \left(1 +\sum_{m=2}^{j} (s-(m-1)\tau)^{(m-1)\alpha-1} \right) ds \Bigg) \\
        %&\leq C\left(\tilde{t}^{\alpha-1}+\tilde{t}^{\alpha-2}+\sum_{m=2}^{n}(\tilde{t}-(m-1)\tau)^{m\alpha-2}\right),
    %\end{aligned}
%\end{equation}
where $t\in(n\tau,(n+1)\tau]$, $n=2,3,\dots,K$.
Thus, by mathematical induction, and in conjunction with the preceding discussion in \eqref{7.15eq2}–-\eqref{eq2_32a3}, it follows that \eqref{lemma2.3eq3} holds.
\end{proof}

%\begin{remark}
%From the proof of Lemma \ref{lemma2.2}, we can also obtain
%$
%|\partial_t^l \tilde{u}_{k\tau}(x,t)|\le C(1+(t-(k-1)\tau)^{k\alpha-l}),\quad l=0,1,2,
%$
%for $t\in((k-1)\tau, k\tau]$, $k=1,2,3,\dots,K$.
%\end{remark}

Now we state the regularity of $w_{k\tau}(x,t)$ in time.
\begin{theorem}\label{thmorem2.1}
	Suppose that assumptions (I) and (II) mentioned in Lemma \ref{lemma2.3} hold. Then for $(x,t)\in[0,L]\times((k-1)\tau,k\tau], 1\leq k\leq K$,
	  $w_{k\tau}(x,t)$ satisfies
	\begin{align}
       &|\partial_{t}^l w_{k\tau}(x,t)|\leq C(1+(t-(k-1)\tau)^{(2k-1)\beta-l}), \quad l=0,1,2. \label{eq2_42}
       %&|\partial_{x}^{n}w_{k\tau}(x,t)|\leq C, \quad n=0,1,2.\label{eq2_42b}
	\end{align}
\end{theorem}
\begin{proof}
%To prove \eqref{eq2_42}, firstly, we will derive expressions for $\partial_t^lw_{k\tau}(x,t)$, $l=0,1,2$.
Since $w(x,t)={}^C_0D_t^\beta\tilde{u}(x,t)$, the combination of \eqref{eq1_5} and \eqref{eq2_40} yields
\begin{equation}\label{eq2_44}
		\begin{aligned}
			w_{k\tau}(x,t)
			=&\frac{1}{\Gamma(1-\beta)}\int_0^t (t-s)^{-\beta}\left[\sum_{i=1}^{\infty}\left(-\phi_{i}(0)\lambda_iB_i(s)
			+B_i(s)H_{k\tau,i}(0)\right.\right.
			\left.\left.+\int_0^{s} B_i(s-\mu)\partial_{\mu} H_{k\tau,i}(\mu)d\mu \right)X_i(x)
            \right]ds	\\
            =&\sum_{i=1}^\infty
		\left[-\lambda_i\phi_i(0)Q_{1,i}(t)+H_{k\tau,i}(0)Q_{1,i}(t)+Q_{2,i}(t)\right]X_i(x),\\
            %&=\frac{1}{\Gamma(1-\beta)}\sum_{i=1}^\infty\bigg[
		%	-\int_0^t (t-s)^{-\beta}\phi_i(0)\lambda_i s^{\alpha-1}E_{\alpha,\alpha}(-\lambda_i s^\alpha)\,ds\\
		%	&\quad+\int_0^t (t-s)^{-\beta}s^{\alpha-1}E_{\alpha,\alpha}(-\lambda_i s^\alpha)H_{k\tau,i}(0)\,ds\\
			%&\quad+\int_0^t (t-s)^{-\beta}\int_0^s (s-\mu)^{\alpha-1}E_{\alpha,\alpha}(-\lambda_i (s-\mu)^\alpha)\partial_\mu H_{k\tau,i}(\mu)\,d\mu\,ds
			%\bigg]X_i(x)
		\end{aligned}
	\end{equation}
	%\begin{equation}\label{eq2_44}
		%\begin{aligned}
			%w_{k\tau}(x,t)
			%=&\frac{1}{\Gamma(1-\beta)}\int_0^t (t-s)^{-\beta}\left[\sum_{i=1}^{\infty}\left(-\phi_{i}(0)\lambda_i s^{\alpha-1}E_{\alpha,\alpha}(-\lambda_is^\alpha)
			%+s^{\alpha-1}E_{\alpha,\alpha}(-\lambda_i s^\alpha)H_{k\tau,i}(0)\right.\right.\\
			%&\left.\left.+\int_0^{s} (s-\mu)^{\alpha-1}E_{\alpha,\alpha}(-\lambda_i (s-\mu)^\alpha)\partial_{\mu} H_{k\tau,i}(\mu)d\mu \right)X_i(x)
           % \right]ds	\\
            %&\sum_{i=1}^\infty
		%\left[-\lambda_i\phi_i(0)Q_{1,i}(t)+H_{k\tau,i}(0)Q_{1,i}(t)+Q_{2,i}(t)\right]X_i(x),\\
            %&=\frac{1}{\Gamma(1-\beta)}\sum_{i=1}^\infty\bigg[
		%	-\int_0^t (t-s)^{-\beta}\phi_i(0)\lambda_i s^{\alpha-1}E_{\alpha,\alpha}(-\lambda_i s^\alpha)\,ds\\
		%	&\quad+\int_0^t (t-s)^{-\beta}s^{\alpha-1}E_{\alpha,\alpha}(-\lambda_i s^\alpha)H_{k\tau,i}(0)\,ds\\
			%&\quad+\int_0^t (t-s)^{-\beta}\int_0^s (s-\mu)^{\alpha-1}E_{\alpha,\alpha}(-\lambda_i (s-\mu)^\alpha)\partial_\mu H_{k\tau,i}(\mu)\,d\mu\,ds
			%\bigg]X_i(x)
		%\end{aligned}
	%\end{equation}
where
\begin{equation}\label{eq2_50}
		Q_{1,i}(t):=\frac{1}{\Gamma(1-\beta)}\int_0^t (t-s)^{-\beta}B_i(s)\,ds,
        \end{equation}
	%\begin{equation}\label{eq2_50}
		%Q_{1,i}(t):=\frac{1}{\Gamma(1-\beta)}\int_0^t (t-s)^{-\beta}s^{\alpha-1}E_{\alpha,\alpha}(-\lambda_i s^\alpha)\,ds,
        %\end{equation}
        \begin{equation}\label{eq2_51}
		Q_{2,i}(t):=\frac{1}{\Gamma(1-\beta)}\int_0^t (t-s)^{-\beta}\int_0^s B_i(\mu)\partial_s H_{k\tau,i}(s-\mu)\,d\mu\,ds.
        \end{equation}
        %\begin{equation}\label{eq2_51}
		%Q_{2,i}(t):=\frac{1}{\Gamma(1-\beta)}\int_0^t (t-s)^{-\beta}\int_0^s \mu^{\alpha-1}E_{\alpha,\alpha}(-\lambda_i \mu^\alpha)\partial_s H_{k\tau,i}(s-\mu)\,d\mu\,ds.
	%\end{equation}
Thus $\partial_{t}^l w_{k\tau}(x,t), l=0,1,2$ can be written as
	\begin{equation}\label{eq2_49}
    \begin{aligned}
		\partial_t^l w_{k\tau}(x,t)
       % &=\frac{1}{\Gamma(1-\beta)}\sum_{i=1}^\infty
       % \partial_t\bigg[
			%-\int_0^t (t-s)^{-\beta}\phi_i(0)\lambda_i s^{\alpha-1}E_{\alpha,\alpha}(-\lambda_i s^\alpha)\,ds\\
			%&\quad+\int_0^t (t-s)^{-\beta}s^{\alpha-1}E_{\alpha,\alpha}(-\lambda_i s^\alpha)H_{k\tau,i}(0)\,ds\\
			%&\quad+\int_0^t (t-s)^{-\beta}\int_0^s \mu^{\alpha-1}E_{\alpha,\alpha}(-\lambda_i \mu^\alpha)\partial_s H_{k\tau,i}(s-\mu)\,d\mu\,ds
			%\bigg]X_i(x)\\
		&=\sum_{i=1}^\infty
		\left[-\lambda_i\phi_i(0)\partial^l_tQ_{1,i}+H_{k\tau,i}(0)\partial_t^lQ_{1,i}+\partial^l_tQ_{2,i}\right]X_i(x).
    \end{aligned}
    \end{equation}
%and
%\begin{equation}\label{eq2_61}
%		\partial_t^2 w_{k\tau}(x,t)
%		=
%		\sum_{i=1}^\infty \left[-\lambda_i\phi_i(0)\partial_t^2Q_{1,i}
%		+H_{k\tau,i}(0)\partial_t^2Q_{1,i}
%		+\partial_t^2Q_{2,i}
%		\right]X_i(x),
%	\end{equation}
	
%Next, we compute $\partial_t^lQ_{1,i}$ and $\partial_t^lQ_{2,i},l=0,1,2.$
%which will be used to obtain the $\partial_tw_{k\tau}(x,t)$ and $\partial_t^2w_{k\tau}(x,t)$.
Using the fact \cite[p. 61]{Gorenflo2020Mittag-Leffler}
	\begin{equation*}
		\frac{1}{\Gamma(\rho)}\int_0^z (z-t)^{\mu-1}E_{\rho,\theta}(\lambda t^\rho)t^{\theta-1}\,dt
		=z^{\mu+\theta-1}E_{\rho,\mu+\theta}(\lambda z^\rho),\quad \mu>0,\ \rho>0,\ \theta>0,\ \lambda\in\mathbb C,
	\end{equation*}
we have
\begin{equation}\label{7.15eq3}
    Q_{1,i}(t)
		=\frac{\Gamma(\alpha)}{\Gamma(1-\beta)} t^\beta E_{\alpha,\beta+1}(-\lambda_i t^\alpha).
\end{equation}
Applying Lemma \ref{lemma2.1} to differentiating both sides of \eqref{7.15eq3} gives
	\begin{equation}\label{eq2_54}
    \begin{aligned}
		\partial^l_tQ_{1,i}(t)=\frac{\Gamma(\alpha)}{\Gamma(1-\beta)}t^{\beta-l}E_{\alpha,\beta+1-l}(-\lambda_i t^\alpha).
        \end{aligned}
	\end{equation}
 %    and
 %    \begin{equation}\label{eq2_62}
	% 	\partial_t^2Q_{1,i}
	% 	=\frac{\Gamma(\alpha)}{\Gamma(1-\beta)}t^{\beta-2}E_{\alpha,\beta-1}(-\lambda_i t^\alpha).
	% \end{equation}
Note that performing integration by parts to \eqref{eq2_51} leads to
\begin{equation}\label{eq2_55}
		\begin{aligned}
			Q_{2,i}(t)
			%&=\frac{1}{\Gamma(1-\beta)}\int_0^t (t-s)^{-\beta}\int_0^s \mu^{\alpha-1}E_{\alpha,\alpha}(-\lambda_i \mu^\alpha)\partial_s H_{k\tau,i}(s-\mu)\,d\mu\,ds\\
			%&=\frac{1}{\Gamma(1-\beta)}\int_0^t\int_0^s \mu^{\alpha-1}E_{\alpha,\alpha}(-\lambda_i \mu^\alpha)\partial_s H_{k\tau,i}(s-\mu)\,d\mu\,
			%d\left(-\frac{(t-s)^{1-\beta}}{1-\beta}\right)\\
			=&-\frac{(t-s)^{1-\beta}}{\Gamma(2-\beta)}
			\int_0^s B_i(\mu)\partial_s H_{k\tau,i}(s-\mu)\,d\mu
			+\int_0^t \frac{(t-s)^{1-\beta}}{\Gamma(2-\beta)}
			\partial_s\left(\int_0^s B_i(\mu)\partial_s H_{k\tau,i}(s-\mu)\,d\mu\right)ds\\
            =&\int_0^t \frac{(t-s)^{1-\beta}}{\Gamma(2-\beta)}
			\left(B_i(s)\partial_s H_{k\tau,i}(0)+\int_0^s B_i(\mu)\partial_s^2 H_{k\tau,i}(s-\mu)\,d\mu\right)ds.
		\end{aligned}
	\end{equation}
	%\begin{equation}\label{eq2_55}
		%\begin{aligned}
			%Q_{2,i}(t)
			%&=\frac{1}{\Gamma(1-\beta)}\int_0^t (t-s)^{-\beta}\int_0^s \mu^{\alpha-1}E_{\alpha,\alpha}(-\lambda_i \mu^\alpha)\partial_s H_{k\tau,i}(s-\mu)\,d\mu\,ds\\
			%&=\frac{1}{\Gamma(1-\beta)}\int_0^t\int_0^s \mu^{\alpha-1}E_{\alpha,\alpha}(-\lambda_i \mu^\alpha)\partial_s H_{k\tau,i}(s-\mu)\,d\mu\,
			%d\left(-\frac{(t-s)^{1-\beta}}{1-\beta}\right)\\
			%=&-\frac{(t-s)^{1-\beta}}{\Gamma(2-\beta)}
			%\int_0^s \mu^{\alpha-1}E_{\alpha,\alpha}(-\lambda_i \mu^\alpha)\partial_s H_{k\tau,i}(s-\mu)\,d\mu
			%+\int_0^t \frac{(t-s)^{1-\beta}}{\Gamma(2-\beta)}
			%\partial_s\left(\int_0^s \mu^{\alpha-1}E_{\alpha,\alpha}(-\lambda_i
            %\mu^\alpha)\partial_s H_{k\tau,i}(s-\mu)\,d\mu\right)ds\\
            %=&\int_0^t \frac{(t-s)^{1-\beta}}{\Gamma(2-\beta)}
			%\left(s^{\alpha-1}E_{\alpha,\alpha}(-\lambda_i s^\alpha)\partial_s H_{k\tau,i}(0)+\int_0^s \mu^{\alpha-1}E_{\alpha,\alpha}(-\lambda_i \mu^\alpha)\partial_s^2 H_{k\tau,i}(s-\mu)\,d\mu\right)ds.
		%\end{aligned}
	%\end{equation}
Furthermore, it follows from Lemma \ref{lemma2.1} and integration by parts that
% Applying  to differentiating both sides of \eqref{eq2_55} and applying the  formula yields
\begin{equation}\label{eq2_56}
\begin{aligned}
			\partial_tQ_{2,i}(t)
            =&\int_0^t \frac{(t-s)^{-\beta}}{\Gamma(1-\beta)}
			\left(B_i(s)\partial_s H_{k\tau,i}(0)
            +\int_0^sB_i(s-\mu)\partial_\mu^2 H_{k\tau,i}(\mu)d\mu \right)ds\\
		    %&=\frac{1}{\Gamma(1-\beta)}\left[\int_0^t (t-s)^{-\beta}s^{\alpha-1}E_{\alpha,\alpha}(-\lambda_i s^\alpha)\partial_s H_{k\tau,i}(0)\,ds
            %+\int_0^t (t-s)^{-\beta}\int_0^s \mu^{\alpha-1}E_{\alpha,\alpha}(-\lambda_i %\mu^\alpha)\partial_s^2 H_{k\tau,i}(s-\mu)\,d\mu\,ds\right]\\
            %=&\int_0^t \frac{(t-s)^{-\beta}}{\Gamma(1-\beta)}s^{\alpha-1}E_{\alpha,\alpha}(-\lambda_i s^\alpha)\partial_s H_{k\tau,i}(0)\,ds\\
            %&+\int_0^t \frac{(t-s)^{-\beta}}{\Gamma(1-\beta)}\int_0^s( s-\mu)^{\alpha-1}E_{\alpha,\alpha}(-\lambda_i ( s-\mu)^\alpha)\partial_\mu^2 H_{k\tau,i}(\mu)d\mu ds\\
            =&\int_0^t \frac{(t-s)^{-\beta}}{\Gamma(1-\beta)}B_i(s)\partial_s H_{k\tau,i}(0)\,ds+\int_0^t\frac{(t-s)^{-\beta}}{\Gamma(1-\beta)}\left[B_i(s-\mu)\partial_\mu H_{k\tau,i}(\mu)\bigg|_0^s\right.
            + \left.\int_0^s C_i(s-\mu)\partial_\mu H_{k\tau,i}(\mu)\,d\mu\right]
            ds\\
            %=&\int_0^t \frac{(t-s)^{-\beta}}{\Gamma(1-\beta)}s^{\alpha-1}E_{\alpha,\alpha}(-\lambda_i s^\alpha)\partial_s H_{k\tau,i}(0)\,ds\\
            %&-\frac{1}{\Gamma(1-\beta)}\int_0^t (t-s)^{-\beta}\left(-s^{\alpha-1}E_{\alpha,\alpha}(-\lambda_i s^\alpha)\partial_s H_{k\tau,i}(0)\right.\\
		%&\left.+\int_0^s (s-\mu)^{\alpha-2}E_{\alpha,\alpha-1}(-\lambda_i (s-\mu)^\alpha)\partial_s H_{k\tau,i}(\mu)\,d\mu\right)ds\\
            =&\int_0^t \frac{(t-s)^{-\beta}}{\Gamma(1-\beta)}\int_0^s C_i(\mu)\partial_s H_{k\tau,i}(s-\mu)\,d\mu ds
            \end{aligned}
	\end{equation}
and
\begin{equation}\label{eq2_64}
		\begin{aligned}
		\partial_t^2Q_{2,i}(t)
		%=&\partial_t\left[\frac{1}{\Gamma(1-\beta)}
		%\int_0^t (t-s)^{-\beta}\int_0^s \mu^{\alpha-2}E_{\alpha,\alpha-1}(-\lambda_i \mu^\alpha)\partial_s H_{k\tau,i}(s-\mu)\,d\mu\,ds
		%\right]\\
         =&\partial_t\left[\frac{-(t-s)^{1-\beta}}{\Gamma(2-\beta)}\int_0^s C_i(\mu)\partial_s H_{k\tau,i}(s-\mu)\,d\mu\bigg|_0^t\right.
         +\left.\int_0^t \frac{(t-s)^{1-\beta}}{\Gamma(2-\beta)}\left(\partial_s\int_0^s C_i(\mu)\partial_s H_{k\tau,i}(s-\mu)\,d\mu\right)ds\right]\\
        =&\int_0^t \frac{(t-s)^{-\beta}}{\Gamma(1-\beta)}
		\left(C_i(s)\partial_s H_{k\tau,i}(0) +\int_0^s C_i(\mu)\partial_s^2H_{k\tau,i}(s-\mu)\,d\mu\right)ds.
        % =&\frac{1}{\Gamma(1-\beta)}\int_0^t (t-s)^{-\beta}\left(s^{\alpha-2}E_{\alpha,\alpha-1}(-\lambda_i s^{\alpha-1})\partial_s H_{k\tau,i}(0)\,ds
        % +\int_0^s \mu^{\alpha-2}E_{\alpha,\alpha-1}(-\lambda_i \mu^\alpha)\partial_s^2 H_{k\tau,i}(s-\mu)\,d\mu\,\right)ds.
	\end{aligned}
	\end{equation}
	%\begin{equation}\label{eq2_64}
		%\begin{aligned}
		%\partial_t^2Q_{2,i}(t)
		%=&\partial_t\left[\frac{1}{\Gamma(1-\beta)}
		%\int_0^t (t-s)^{-\beta}\int_0^s \mu^{\alpha-2}E_{\alpha,\alpha-1}(-\lambda_i \mu^\alpha)\partial_s H_{k\tau,i}(s-\mu)\,d\mu\,ds
		%\right]\\
         %=&\partial_t\left[\frac{-(t-s)^{1-\beta}}{\Gamma(2-\beta)}\int_0^s \mu^{\alpha-2}E_{\alpha,\alpha-1}(-\lambda_i \mu^\alpha)\partial_s H_{k\tau,i}(s-\mu)\,d\mu\bigg|_0^t\right.\\
        % &+\left.\int_0^t \frac{(t-s)^{1-\beta}}{\Gamma(2-\beta)}\left(\partial_s\int_0^s \mu^{\alpha-2}E_{\alpha,\alpha-1}(-\lambda_i \mu^\alpha)\partial_s H_{k\tau,i}(s-\mu)\,d\mu\right)ds\right]\\
        %=&\int_0^t \frac{(t-s)^{-\beta}}{\Gamma(1-\beta)}
		%\left(s^{\alpha-2}E_{\alpha,\alpha-1}(-\lambda_i s^{\alpha-1})\partial_s H_{k\tau,i}(0) +\int_0^s \mu^{\alpha-2}E_{\alpha,\alpha-1}(-\lambda_i \mu^\alpha)\partial_s^2H_{k\tau,i}(s-\mu)\,d\mu\right)ds.
        % =&\frac{1}{\Gamma(1-\beta)}\int_0^t (t-s)^{-\beta}\left(s^{\alpha-2}E_{\alpha,\alpha-1}(-\lambda_i s^{\alpha-1})\partial_s H_{k\tau,i}(0)\,ds
        % +\int_0^s \mu^{\alpha-2}E_{\alpha,\alpha-1}(-\lambda_i \mu^\alpha)\partial_s^2 H_{k\tau,i}(s-\mu)\,d\mu\,\right)ds.
	%\end{aligned}
	%\end{equation}

For $\partial^l_tQ_{1,i}$, Lemma \ref{lemma2.1} and \eqref{eq2_54} indicate that
\begin{equation}\label{eqbu716A}
    |\partial^l_t Q_{1,i}|\leq Ct^{\beta-l}, l=0,1,2.
\end{equation}
Now we discuss the boundedness of $\sum_{i=1}^\infty|\partial^l_t Q_{2,i}|.$
When $l=0$ and $t\in ((k-1)\tau,k\tau]$, it follows from Lemma \ref{lemma2.1} and Lemma \ref{lemma2.3} that
\begin{equation}\label{7.16eq5}
    \begin{aligned}
\sum_{i=1}^\infty| Q_{2,i}|
=&\sum_{i=1}^\infty\left|\int_0^t \frac{(t-s)^{-\beta}}{\Gamma(1-\beta)}\int_0^s B_i(s-\mu)\partial_\mu H_{k\tau,i}(\mu)\,d\mu\,ds\right|\\
		\le &C\left[
            \int_0^{\min\{\tau,t\}}  (t-s)^{-\beta}\int_0^s(s-\mu)^{\alpha-1}\sum_{i=1}^\infty |\partial_\mu H_{\tau,i}(\mu)|\,d\mu\,ds\right.\\
            &+\left.\sum_{j=2}^{k}\int_{(j-1)\tau}^{\min\{j\tau,t\}}
            (t-s)^{-\beta}\int_0^s(s-\mu)^{\alpha-1}\sum_{i=1}^\infty |\partial_\mu H_{k\tau,i}(\mu)|\,d\mu\,ds
            \right] \\
            \leq& C\left[ \int_0^{\min\{\tau,t\}}  (t-s)^{-\beta}\int_0^s(s-\mu)^{\alpha-1} d\mu ds
            + \sum_{j=2}^{k}\int_{(j-1)\tau}^{\min\{j\tau,t\}} (t-s)^{-\beta}\left(\int_0^{\min\{\tau,s\}} (s-\mu)^{\alpha-1}d\mu\right.\right.\\
        &\left.\left.+\sum_{l_1=2}^{j}\int_{(l_1-1)\tau}^{\min\{l_1\tau,s\}}(s-\mu)^{\alpha-1}\left(1+\sum_{m=2}^{l_1}(\mu-(m-1)\tau)^{(m-1)\alpha-1}\right)d\mu\right) ds\right].\\
    \end{aligned}
\end{equation}
%\begin{equation}\label{7.16eq5}
  %  \begin{aligned}
%\sum_{i=1}^\infty| Q_{2,i}|
%=&\sum_{i=1}^\infty\left|\int_0^t \frac{(t-s)^{-\beta}}{\Gamma(1-\beta)}\int_0^s (s-\mu)^{\alpha-1}E_{\alpha,\alpha}(-\lambda_i (s-\mu)^\alpha)\partial_\mu H_{k\tau,i}(\mu)\,d\mu\,ds\right|\\
		%\le &C\left[
            %\int_0^{\min\{\tau,t\}}  (t-s)^{-\beta}\int_0^s(s-\mu)^{\alpha-1}\sum_{i=1}^\infty |\partial_\mu H_{\tau,i}(\mu)|\,d\mu\,ds\right.\\
            %&+\left.\sum_{j=2}^{k}\int_{(j-1)\tau}^{\min\{j\tau,t\}}
           % (t-s)^{-\beta}\int_0^s(s-\mu)^{\alpha-1}\sum_{i=1}^\infty |\partial_\mu H_{k\tau,i}(\mu)|\,d\mu\,ds
            %\right] \\
          %  \leq& C\left[ \int_0^{\min\{\tau,t\}}  (t-s)^{-\beta}\int_0^s(s-\mu)^{\alpha-1} d\mu ds+ \sum_{j=2}^{k}\int_{(j-1)\tau}^{\min\{j\tau,t\}} (t-s)^{-\beta}\left(\int_0^{\min\{\tau,s\}} (s-\mu)^{\alpha-1}d\mu\right.\right.\\
         %&\left.\left.+\sum_{l_1=2}^{j}\int_{(l_1-1)\tau}^{\min\{l_1\tau,s\}}(s-\mu)^{\alpha-1}\left(1+\sum_{m=2}^{l_1}(\mu-(m-1)\tau)^{(m-1)\alpha-1}\right)d\mu\right) ds\right].\\
    %\end{aligned}
%\end{equation}
What is more, by using Lemma \ref{lemma2.2}, it leads to
\begin{equation}\label{7.16eq6}
    \begin{aligned}
\sum_{i=1}^\infty| Q_{2,i}|
\leq& C\left[\int_0^{\min\{\tau,t\}}  (t-s)^{-\beta} s^{\alpha} ds
            + \sum_{j=2}^{k}\int_{(j-1)\tau}^{\min\{j\tau,t\}} (t-s)^{-\beta}\left(s^{\alpha}
            +\sum_{m=2}^{j}(s-(m-1)\tau)^{m\alpha-1}\right) ds\right]\\
			\leq& C\left[t^{\beta+1}+\sum_{m=2}^{k}
            (t-(m-1)\tau)^{(2m-1)\beta}\right].\\
			%\leq &C\left[1+(t-(k-1)\tau)^{(2k-1)\beta}\right].
    \end{aligned}
\end{equation}
Based on \eqref{eq2_56} and \eqref{eq2_64}, an analogous argument to $\sum_{i=1}^\infty| \partial_tQ_{2,i}|$ and $\sum_{i=1}^\infty| \partial^2_tQ_{2,i}|$ gives
\begin{equation}\label{7.16eq7}
    \begin{aligned}
\sum_{i=1}^\infty| \partial_tQ_{2,i}|
=&\sum_{i=1}^\infty\left|
\int_0^t \frac{(t-s)^{-\beta}}{\Gamma(1-\beta)}\int_0^s C_i(s-\mu)\partial_s H_{k\tau,i}(\mu)\,d\mu ds
        \right|
			\leq C\left[t^{\beta}+\sum_{m=2}^{k}(t-(m-1)\tau)^{(2m-1)\beta-1}\right],\\
			%&\leq C[1+(t-(k-1)\tau)^{(2k-1)\beta-1}]
    \end{aligned}
\end{equation}
%\begin{equation}\label{7.16eq7}
    %\begin{aligned}
%\sum_{i=1}^\infty| \partial_tQ_{2,i}|
%=&\sum_{i=1}^\infty\left|
%\int_0^t \frac{(t-s)^{-\beta}}{\Gamma(1-\beta)}\int_0^s (s-\mu)^{\alpha-2}E_{\alpha,\alpha-1}(-\lambda_i (s-\mu)^\alpha)\partial_s H_{k\tau,i}(\mu)\,d\mu ds
          %  \right|\\
			%\leq& C\left[t^{\beta}+\sum_{m=2}^{k}(t-(m-1)\tau)^{(2m-1)\beta-1}\right],\\
			%&\leq C[1+(t-(k-1)\tau)^{(2k-1)\beta-1}]
   % \end{aligned}
%\end{equation}
and
\begin{equation}\label{7.17eq1}
    \begin{aligned}
\sum_{i=1}^\infty| \partial_t^2Q_{2,i}|
=&\sum_{i=1}^\infty\left|
\int_0^t \frac{(t-s)^{-\beta}}{\Gamma(1-\beta)}
		C_i(s)\partial_s H_{k\tau,i}(0)ds\right.
        +\left.\int_0^t \frac{(t-s)^{-\beta}}{\Gamma(1-\beta)}\int_0^s C_i(s-\mu)\partial_s^2H_{k\tau,i}(\mu)\,d\mu ds\right|\\
        \leq &
        C \int_0^t \frac{(t-s)^{-\beta}}{\Gamma(1-\beta)}
		s^{\alpha-2}ds+\sum_{i=1}^\infty\left|\int_0^t \frac{(t-s)^{-\beta}}{\Gamma(1-\beta)}\int_0^s C_i(s-\mu)\partial_s^2H_{k\tau,i}(\mu)\,d\mu ds\right|\\
			\le &C\left[
            t^{\beta-1}+t^{\beta}
            +\sum_{m=2}^k(t-(m-1)\tau)^{(2m-1)\beta-2}\right].\\
			%\le& C(1+(t-(k-1)\tau)^{(2k-1)\beta-2}),
\end{aligned}
\end{equation}
%\begin{equation}\label{7.17eq1}
   % \begin{aligned}
%\sum_{i=1}^\infty| \partial_t^2Q_{2,i}|
%=&\sum_{i=1}^\infty\left|
%\int_0^t \frac{(t-s)^{-\beta}}{\Gamma(1-\beta)}
		%s^{\alpha-2}E_{\alpha,\alpha-1}(-\lambda_i s^{\alpha-1})\partial_s H_{k\tau,i}(0)ds\right.\\
        %&\left.+\int_0^t \frac{(t-s)^{-\beta}}{\Gamma(1-\beta)}\int_0^s (s-\mu)^{\alpha-2}E_{\alpha,\alpha-1}(-\lambda_i (s-\mu)^\alpha)\partial_s^2H_{k\tau,i}(\mu)\,d\mu ds\right|\\
        %\leq &
        %C \int_0^t \frac{(t-s)^{-\beta}}{\Gamma(1-\beta)}
		%s^{\alpha-2}ds+\sum_{i=1}^\infty\left|\int_0^t \frac{(t-s)^{-\beta}}{\Gamma(1-\beta)}\int_0^s (s-\mu)^{\alpha-2}E_{\alpha,\alpha-1}(-\lambda_i (s-\mu)^\alpha)\partial_s^2H_{k\tau,i}(\mu)\,d\mu ds\right|\\
			%\le &C\left[
            %t^{\beta-1}+t^{\beta}
            %+\sum_{m=2}^k(t-(m-1)\tau)^{(2m-1)\beta-2}\right].\\
			%\le& C(1+(t-(k-1)\tau)^{(2k-1)\beta-2}),
%\end{aligned}
%\end{equation}
Finally, the combination of \eqref{eqbu716A}--\eqref{7.17eq1} and $|X_i(x)|\le C$ \cite[p. 335]{Luchko2012Initial-boundary-value} immediately yields \eqref{eq2_42}.

\end{proof}
\section{Fully discrete finite element scheme and its stability}
Considering the weak singular behavior of the solution and the constant delay term involved in the considered problem, we employ the following mesh mentioned in \cite{Tan2022L1} to divide the interval $[-\tau,K\tau]$
\begin{equation}\label{eq3_1}
	t_{-2N}=-\tau,\quad t_{n}=
	\begin{cases}
		\displaystyle \frac{\tau}{2}\left( \frac{n-2(k-1)N}{N} \right)^{r}+(k-1)\tau, & 2(k-1)N+1\leq n<(2k-1)N, \\
		\displaystyle -\frac{\tau}{2}\left( \frac{2kN-n}{N} \right)^{r}+k\tau, & (2k-1)N\leq n\leq 2kN,
	\end{cases}
\end{equation}
where $0\leq k\leq K$, $r\geq 1$, and $N$ is a positive integer. Let $\rho_{n}=t_{n}-t_{n-1}$. The mesh \eqref{eq3_1} clearly satisfies
\begin{equation}\label{eq3_2}
	\rho_{1}=t_{1}\simeq N^{-r},\quad
	\rho_{n}\lesssim \left( \rho_{1} \right)^{1/r}\left( t_{n}-(k-1)\tau \right)^{1-1/r}.
\end{equation}
Denote $v^{n}=v(t_{n})$ and $\nabla_{t}v^{n}=v^{n}-v^{n-1}$. We use the nonuniform L1 formula
\cite{Liao2018Sharp}
%[Changpin Li book???]
to approximate the Caputo fractional derivative
${}_{0}^{C}D_{t}^{\beta}v(t_{n})$ having
\begin{equation}\label{eq3_3}	{}_{0}^{C}D_{t}^{\beta}v(t_{n})\approx\sum_{k=1}^{n}a_{n-k}^{(n)}\nabla_{t}v^{k}:=D^{\beta}v^{n},
\end{equation}
where
\begin{equation}\label{eq3_4}
	a_{n-k}^{(n)}=\frac{\omega_{2-\beta}(t_{n}-t_{k-1})-\omega_{2-\beta}(t_{n}-t_{k})}{\rho_k},\quad 1\leq n\leq 2KN.
\end{equation}
Furthermore,
let $X_h^q$ be  continuous piecewise polynomial finite element space of degree $q$ of $H_0^1([0,L])$ based on quasi-uniform partition of $[0,L]$ with maximum diameter $h.$
Now, by employing the above $L1$ approximation, we establish the following fully discrete finite element scheme for problem \eqref{eq2_36b1}--\eqref{eq2_36b6}: find $\left( \tilde{u}_{h}^{n},w_{h}^{n} \right)\in X^q_{h}\times X^q_{h}$, $1\leq n\leq 2KN$ such that for all $\chi\in X^q_{h}$ satisfying
\begin{align}
	& \left( D^{\beta}w_{h}^{n},\chi \right)+p\left( \nabla\tilde{u}_{h}^{n},\nabla \chi \right)
	= a\left( \tilde{u}_{h}^{n},\chi \right)+b\left( \tilde{u}_{h}^{n-2N},\chi \right)+\left( g^{n},\chi \right), \label{eq3_5}\\
	& \left( D^{\beta}\tilde{u}_{h}^{n},\chi \right)=\left( w_{h}^{n},\chi \right)\label{eq3_6},
\end{align}
where $\tilde{u}_{h}^{n}$ and $w_{h}^{n}$ are the approximations of $\tilde{u}(x,t_{n})$ and $w(x,t_{n})$,  $g^{n}=g(x,t_{n})$, and for $1\leq n\leq 2N$,
$\tilde{u}_{h}^{n-2N}\in X_h^q$ and $w_{h}^{n-2N}\in X_h^q$ are suitable approximations of $\tilde{u}(x,t_{n-2N})$ and $w(x,t_{n-2N})$.

In order to discuss the stability of the fully discrete finite element scheme \eqref{eq3_5}--\eqref{eq3_6}, we first introduce some lemmas.

\begin{lemma}\label{lemma3_1}\cite[Lemma 3.3]{Huang2023local}
	For $1\leq n\leq 2KN$, if $v^{n}\in L^{2}([0,L])$, then
	%\begin{equation}\label{eq3_8}
    $
		\left( D^{\beta}v^{n},v^{n} \right)\geq \left\| v^{n} \right\|\cdot D^{\beta}\left\| v^{n} \right\|.
        $
	%\end{equation}
\end{lemma}
\begin{lemma}\label{lemma3_2}\cite[Lemma 3.4]{Huang2023local} Let $v_{i}^{n}$, $1\leq n\leq 2KN,$ $1\leq i\leq l$ be nonnegative function. Then, for each $1\leq n\leq 2KN,$ one has
	%\begin{equation}\label{eq3_9}
    $
		\sum_{i=1}^{l}v_{i}^{n}\cdot D^{\beta}v _{i}^{n}
		\geq D^{\beta}\sqrt{\sum_{i=1}^{l}\left( v_{i}^{n} \right)^{2}}\sqrt{\sum_{i=1}^{l}\left( v_{i}^{n} \right)^{2}}.
        $
	%\end{equation}
	% where $n=1,2,\cdots,2KN$.
\end{lemma}

\begin{lemma}\label{lemma3_3} Assume that $\left( \tilde{u}_{h}^{n},w_{h}^{n} \right)\in X^q_{h}\times X^q_{h}$, $1\leq n\leq 2KN$ satisfies
	\begin{align}
		&\left( D^{\beta}w_{h}^{n},\chi \right)+p\left( \nabla\tilde{u}_{h}^{n},\nabla \chi \right)
		=a\left( \tilde{u}_{h}^{n},\chi \right)+b\left( \tilde{u}_{h}^{n-2N},\chi \right)+\left( h_{1}^{n},\chi \right),\quad \forall \chi\in X^q_{h},\label{eq3_10}\\
		&\left( D^{\beta}\tilde{u}_{h}^{n},\chi \right)-\left( w_{h}^{n},\chi \right)=\left( h_{2}^{n},\chi \right),\quad \forall \chi\in X^q_{h},
		\label{eq3_11}
	\end{align}
	 where $h_1^n\in L^2([0,L])$ and $h^n_2\in H^1([0,L]).$
     % \tilde{u}_{h}^{n-2N}=\phi(x,t_{n-2N})-t_{n-2N}\psi(x,0)$, $w_{h}^{n-2N}=\left. {}_{0}^{C}D_{t}^{\beta}\left(\phi(x,t)-t\psi(x,0)\right) \right|_{t=t_{n-2N}}$, $1\leq n\leq 2N$.
	Then we have
	\begin{equation}\label{eq3_12}
		\left( D^{\beta}+\frac{a}{\sqrt{p}} \right)\sqrt{\left\| w_{h}^{n} \right\|^{2}+p\left\| \nabla\tilde{u}_{h}^{n} \right\|^{2}}
		\lesssim \frac{b}{\sqrt{p}}\sqrt{\left\| w_{h}^{n-2N} \right\|^{2}+p\left\| \nabla\tilde{u}_{h}^{n-2N} \right\|^{2}}
		+\sqrt{\left\| h_{1}^{n} \right\|^{2}+p\left\| \nabla h_{2}^{n} \right\|^{2}}.
	\end{equation}
\end{lemma}
\begin{proof}
Define the discrete Laplacian operator $\Delta_h\colon X^q_h\to X^q_h$ by \cite[Eq. (1.33)]{Thomee2007Galerkin}
\begin{equation}\label{eq3_7}
	(\Delta_h v,\ w)=-(\nabla v,\ \nabla w), \forall\, v,\,w\in X_h^q.
\end{equation}
	Taking $\chi=w_{h}^{n}$ in \eqref{eq3_10} and $\chi=p\Delta_h\tilde{u}_{h}^{n}$ in \eqref{eq3_11} yields
	\begin{align}
		&\left( D^{\beta}w_{h}^{n},w_{h}^{n} \right)+p\left( \nabla\tilde{u}_{h}^{n},\nabla w_{h}^{n} \right)
		=a\left( \tilde{u}_{h}^{n},w_{h}^{n} \right)+b\left( \tilde{u}_{h}^{n-2N},w_{h}^{n} \right)+\left( h_{1}^{n},w_{h}^{n} \right)\label{eq3_13},\\
		&p\left( D^{\beta}\nabla\tilde{u}_{h}^{n},\nabla\tilde{u}_{h}^{n} \right)
		=p\left( \nabla w_{h}^{n},\nabla\tilde{u}_{h}^{n} \right)
		+p\left( \nabla h_{2}^{n},\nabla\tilde{u}_{h}^{n} \right).\label{eq3_14}
	\end{align}
	By using \eqref{eq3_14} to replace the second term on the left-hand side of \eqref{eq3_13}, we obtain
	\begin{equation}\label{eq3_15}
		\left( D^{\beta}w_{h}^{n},w_{h}^{n} \right)+p\left( D^{\beta}\nabla\tilde{u}_{h}^{n},\nabla\tilde{u}_{h}^{n} \right)
		=a\left( \tilde{u}_{h}^{n},w_{h}^{n} \right)+b\left( \tilde{u}_{h}^{n-2N},w_{h}^{n} \right)+\left( h_{1}^{n},w_{h}^{n} \right)+p\left( \nabla h_{2}^{n},\nabla\tilde{u}_{h}^{n} \right).
	\end{equation}
	Furthermore, according to Lemma \ref{lemma3_1}, the Cauchy--Schwarz inequality and the Poincar\'e inequality, \eqref{eq3_15} can be rewritten as
	% \begin{equation}\label{eq3_16}
	% 	D^{\beta}\left\| w_{h}^{n} \right\|\left\| w_{h}^{n} \right\|+pD^{\beta}\left\| \nabla\tilde{u}_{h}^{n} \right\|\left\| \nabla\tilde{u}_{h}^{n} \right\| \\
	% 	\leq (-a)\left\| \tilde{u}_{h}^{n} \right\|\left\| w_{h}^{n} \right\|+b\left\| \tilde{u}_{h}^{n-2N} \right\|\left\| w_{h}^{n} \right\|+\left\| h_{1}^{n} \right\|\left\| w_{h}^{n} \right\|+p\left\| \nabla h_{2}^{n} \right\|\left\| \nabla\tilde{u}_{h}^{n} \right\|
	% \end{equation}
	% By  $\left\| \tilde{u}_{h}^{n} \right\|\lesssim\left\| \nabla\tilde{u}_{h}^{n} \right\|$, \eqref{eq3_16} can be rewritten as
	\begin{equation}\label{eq3_17}
    \begin{aligned}
        \left\| w_{h}^{n} \right\|\cdot D^{\beta}\left\| w_{h}^{n} \right\|+p\left\| \nabla\tilde{u}_{h}^{n} \right\|\cdot D^{\beta}\left\| \nabla\tilde{u}_{h}^{n} \right\| \leq& (-a)\left\| \tilde{u}_{h}^{n} \right\|\left\| w_{h}^{n} \right\|+b\left\| \tilde{u}_{h}^{n-2N} \right\|\left\| w_{h}^{n} \right\|+\left\| h_{1}^{n} \right\|\left\| w_{h}^{n} \right\|+p\left\| \nabla h_{2}^{n} \right\|\left\| \nabla\tilde{u}_{h}^{n} \right\|\\
		\lesssim& (-a)\left\| \nabla\tilde{u}_{h}^{n} \right\|\left\| w_{h}^{n} \right\|+b\left\| \nabla\tilde{u}_{h}^{n-2N} \right\|\left\| w_{h}^{n} \right\|+\left\| h_{1}^{n} \right\|\left\| w_{h}^{n} \right\|+p\left\| \nabla h_{2}^{n} \right\|\left\| \nabla\tilde{u}_{h}^{n} \right\|.
    \end{aligned}
	\end{equation}
	From Lemma \ref{lemma3_2} and the H\"older inequality, \eqref{eq3_17} leads to
	\begin{equation}\label{eq3_18}
		\begin{aligned}
			&D^{\beta}\sqrt{\left\| w_{h}^{n} \right\|^{2}+p\left\| \nabla\tilde{u}_{h}^{n} \right\|^{2}}\cdot\sqrt{\left\| w_{h}^{n} \right\|^{2}+p\left\| \nabla\tilde{u}_{h}^{n} \right\|^{2}}
			+ \frac{a}{\sqrt{p}}\sqrt{\left\| w_{h}^{n} \right\|^{2}+p\left\| \nabla\tilde{u}_{h}^{n} \right\|^{2}}\cdot\sqrt{\left\| w_{h}^{n} \right\|^{2}+p\left\| \nabla\tilde{u}_{h}^{n} \right\|^{2}} \\
			&\lesssim\frac{b}{\sqrt{p}}\sqrt{\left\| w_{h}^{n-2N} \right\|^{2}+p\left\| \nabla\tilde{u}_{h}^{n-2N} \right\|^{2}}\cdot\sqrt{\left\| w_{h}^{n} \right\|^{2}+p\left\| \nabla\tilde{u}_{h}^{n} \right\|^{2}}
			+\sqrt{\left\| h_{1}^{n} \right\|^{2}+p\left\| \nabla h_{2}^{n} \right\|^{2}}\cdot\sqrt{\left\| w_{h}^{n} \right\|^{2}+p\left\| \nabla\tilde{u}_{h}^{n} \right\|^{2}}.
		\end{aligned}
	\end{equation}
    Therefore, by dividing both sides of \eqref{eq3_18} by $\sqrt{\left\| w_{h}^{n} \right\|^{2}+p\left\| \nabla\tilde{u}_{h}^{n} \right\|^{2}}$, we immediately obtain \eqref{eq3_12}.
	% Canceling  $$ on both sides yields \eqref{eq3_12}.
\end{proof}

\begin{lemma}\label{lemma3_4}\cite{Liao2018Sharp,Tan2022L1}
	For given nonnegative sequences $\left\{ y^{i},z^{i}\mid 1\leq i\leq 2KN \right\}$ and $\left\{ \lambda_{i}\mid 0\leq i\leq 2KN-1 \right\}$,
	suppose that there exists a constant $\Lambda$ independent of the time step-size such that $\Lambda\geq\sum_{i=0}^{2KN-1}\lambda_{i}$ and the maximum time step-size satisfies
	$\max_{1\leq n\leq 2KN}\rho_{n}\leq\frac{1}{\sqrt[\beta]{2\Gamma(2-\beta)\Lambda}}$.
	Then for any grid function $\left\{ V^{k}\mid k\geq 0 \right\}$ satisfying
	\begin{equation}\label{eq3_19}
		D^{\beta}V^{n}\leq\sum_{k=1}^{n}\lambda_{n-k}V^{k}+y^{n}+z^{n},\quad 1\leq n\leq 2KN,
	\end{equation}
	we have
	\begin{equation}\label{eq3_20}
		V^{n}\leq 2E_{\beta,1}\left( 2\max(1,\sigma)\Lambda(t_{n})^{\beta} \right)\left( V^{0}+\max_{1\leq k\leq n}\sum_{i=1}^{k}P_{k-i}^{(k)}y^{i}+\omega_{1+\beta}(t_{n})\max_{1\leq i\leq n}z^{i} \right),
	\end{equation}
	where $\sigma=\max \left\{ 3^{r-1},2^{r-1} \right\}$, $r\geq1$ and
	$
		P_{k-i}^{(k)}:=
		\begin{cases}
			\displaystyle \frac{1}{a_{0}^{(k)}}, & i=k, \\
			\displaystyle \frac{1}{a_{0}^{(i)}}\sum_{j=i+1}^{k}P_{k-j}^{(k)}\left( a_{j-(i+1)}^{(j)}-a_{j-i}^{(j)} \right), & 1\leq i\leq k-1.
		\end{cases}
	$
\end{lemma}

Now, we show the stability  result for the fully discrete finite element scheme \eqref{eq3_5}--\eqref{eq3_6}.
\begin{theorem}\label{theorem3_1}
	Let $\left( w_{h}^{n},\tilde{u}_{h}^{n} \right)$ be the numerical solution of the fully discrete finite element scheme \eqref{eq3_5}--\eqref{eq3_6}. Then
	\begin{equation}\label{eq3_21}
		\sqrt{\left\| w_{h}^{n} \right\|^{2}+p\left\| \nabla\tilde{u}_{h}^{n} \right\|^{2}}
		\leq 2E_{\beta,1}\left( 2\max(1,\sigma)|\Lambda|(t_{n})^{\beta} \right)\left( \sqrt{\left\| w_{h}^{0} \right\|^{2}+p\left\| \nabla\tilde{u}_{h}^{0} \right\|^{2}}+\omega_{1+\beta}(t_{n})\max_{1\leq i\leq n}z^{i} \right),
	\end{equation}
	where
	% \begin{equation*}
 %        % \quad g(x,t_{i}) =pt_{i}\Delta\psi(x,0)+a t_{i}\psi(x,0)+b(t_{i}-\tau)\psi(x,0)+f(x,t_{i}),
	% \end{equation*}
	$\Lambda=-\frac{a}{\sqrt{p}}$, $z^{i}=\frac{bc}{\sqrt{p}}+\left\| g^i \right\|$ for $1\leq i\leq 2N$ with $c=\sup_{0\leq k\leq 2N}\sqrt{\left\| w_{h}^{k-2N} \right\|^{2}+p\left\| \nabla\tilde{u}_{h}^{k-2N} \right\|^{2}}$,
	and
	$\Lambda=-\frac{a}{\sqrt{p}}+\frac{b}{\sqrt{p}}$, $z^{i}:=\left\| g^i \right\|$ for $2N+1\leq i\leq 2KN$.
\end{theorem}
\begin{proof}
	Taking $h_{1}^{n}=g^{n}$ and $h_{2}^{n}=0$ in Lemma \ref{lemma3_3}, we obtain
	\begin{equation}\label{eq3_22}
		\left(D^{\beta}+\frac{a}{\sqrt{p}}\right)\sqrt{\left\| w_{h}^{n} \right\|^{2}+p\left\| \nabla\tilde{u}_{h}^{n} \right\|^{2}}
		\lesssim \frac{b}{\sqrt{p}}\sqrt{\left\| w_{h}^{n-2N} \right\|^{2}+p\left\| \nabla\tilde{u}_{h}^{n-2N} \right\|^{2}}+\left\| g^{n} \right\|,\quad 1\leq n\leq 2KN.
	\end{equation}
	When $1\leq n\leq 2N$, \eqref{eq3_22} can be rewritten as
	\begin{equation}\label{eq3_23}
		\left(D^{\beta}+\frac{a}{\sqrt{p}}\right)\sqrt{\left\| w_{h}^{n} \right\|^{2}+p\left\| \nabla\tilde{u}_{h}^{n} \right\|^{2}}\lesssim\frac{bc}{\sqrt{p}}+\left\| g^{n} \right\|.
	\end{equation}
	Take $V^{n}=\sqrt{\left\| w_{h}^{n} \right\|^{2}+p\left\| \nabla\tilde{u}_{h}^{n} \right\|^{2}}$, $y^{n}=0$, $z^{n}=\frac{bc}{\sqrt{p}}+\left\| g^{n} \right\|$,
	$\lambda_{n-k}=\begin{cases} 0, & k\neq n, \\ \frac{-a}{\sqrt{p}}, & k=n, \end{cases}$ and $\Lambda=\frac{-a}{\sqrt{p}}$.
	According to Lemma \ref{lemma3_4} and \eqref{eq3_23}, it is clear that \eqref{eq3_21} holds.
	When $2N+1\leq n\leq 2KN$, we choose $V^{n}=\sqrt{\left\| w_{h}^{n} \right\|^{2}+p\left\| \nabla\tilde{u}_{h}^{n} \right\|^{2}}$, $y^{n}=0$,
	$z^{n}=\left\| g^{n} \right\|$,
	$\lambda_{n-k}=\begin{cases} \frac{b}{\sqrt{p}}, & k=n-2N, \\ \frac{-a}{\sqrt{p}}, & k=n, \\ 0, & \text{otherwise}, \end{cases}$
	and $\Lambda=\frac{-a}{\sqrt{p}}+\frac{b}{\sqrt{p}}$. It follows from Lemma
	\ref{lemma3_4} and \eqref{eq3_22} that \eqref{eq3_21} holds.
\end{proof}
\section{Convergence analysis}
% In this section, we discuss the convergence of the fully discrete finite element scheme \eqref{eq3_5}--\eqref{eq3_6}.
Define the Ritz projection operator $R_{h}:H_{0}^{1}([0,L])\rightarrow X^q_{h}$ satisfying
\begin{equation}\label{eq4_1}
	\left( \nabla R_{h}v,\nabla \chi \right)=\left( \nabla v,\nabla \chi \right),\quad \forall v\in H_{0}^{1}([0,L]),\forall \chi\in X^q_{h}.
\end{equation}
From \cite[Lemma 1.1]{Thomee2007Galerkin}, it has the following property
\begin{equation}\label{eq4_2}
	\left\| v-R_{h}v \right\|+h\left\|\nabla( v-R_{h}v) \right\|\leq C h^{q+1}\left\| v \right\|_{H^{q+1}([0,L])},\quad \forall v\in H^{q+1}([0,L])\cap H_{0}^{1}([0,L]).
\end{equation}
Let $e_{\tilde{u}}^{n}=\tilde{u}_{h}^{n}-R_{h}\tilde{u}^{n}$ and $e_{w}^{n}=w_{h}^{n}-R_{h}w^{n}$.
By taking $\chi=e_{w}^{n}$ in \eqref{eq3_5} and $\chi=p\Delta _he_{\tilde{u}}^{n}$ in \eqref{eq3_6}, it follows from \eqref{eq3_7} that
\begin{align}
	\left( D^{\beta}e_{w}^{n},e_{w}^{n} \right)+p\left( \nabla e_{\tilde{u}}^{n},\nabla e_{w}^{n} \right)
	=&a\left( e_{\tilde{u}}^{n},e_{w}^{n} \right)+b\left( e_{\tilde{u}}^{n-2N},e_{w}^{n} \right)+\left( g^{n},e_{w}^{n} \right)
	+\left( aR_{h}\tilde{u}^{n}+bR_{h}\tilde{u}^{n-2N},e_{w}^{n} \right)\nonumber\\
	&-\left( D^{\beta}R_{h}w^{n},e_{w}^{n} \right)-p\left( \nabla R_{h}\tilde{u}^{n},\nabla  e_{w}^{n} \right),\label{eq4_3}\\
	p\left( D^{\beta}\nabla e_{\tilde{u}}^{n},p\nabla e_{\tilde{u}}^{n} \right)
	=&p\left( \nabla e_{w}^{n},\nabla e_{\tilde{u}}^{n} \right)-p\left( D^{\beta}\nabla R_{h}\tilde{u}^{n},\nabla e_{\tilde{u}}^{n} \right)+p\left( \nabla R_{h}w^{n},\nabla e_{\tilde{u}}^{n} \right). \label{eq4_4}
\end{align}
In addition, according to \eqref{eq2_36b1}, \eqref{eq2_36b2} and \eqref{eq4_1}, one has
\begin{equation}
	p\left( \nabla R_{h}\tilde{u}^{n},\nabla e_{w}^{n} \right)
	=p\left( \nabla\tilde{u}^{n},\nabla e_{w}^{n} \right)
	=\left( a\tilde{u}^{n}+b\tilde{u}^{n-2N}+g^{n}-{}_{0}^{C}D_{t}^{\beta}w^{n},e_{w}^{n} \right),\label{eq4_5}
\end{equation}
\begin{equation}
	p\left( \nabla R_{h}w^{n},\nabla e_{\tilde{u}}^{n} \right)
	=p\left( \nabla w^{n},\nabla e_{\tilde{u}}^{n} \right)
	=p\left( {}_{0}^{C}D_{t}^{\beta}\nabla\tilde{u}^{n},\nabla e_{\tilde{u}}^{n} \right).\label{eq4_6}
\end{equation}
Therefore, substituting \eqref{eq4_5} and \eqref{eq4_6} into the last terms on the right-hand sides of \eqref{eq4_3} and \eqref{eq4_4}, respectively, and rearranging the resulting expressions, we obtain
\begin{align}
	\left( D^{\beta}e_{w}^{n},e_{w}^{n} \right)+p\left( \nabla e_{\tilde{u}}^{n},\nabla e_{w}^{n} \right)
	=&a\left( e_{\tilde{u}}^{n},e_{w}^{n} \right)+b\left( e_{\tilde{u}}^{n-2N},e_{w}^{n} \right)
	+a\left( R_{h}\tilde{u}^{n}-\tilde{u}^{n},e_{w}^{n} \right)+b\left( R_{h}\tilde{u}^{n-2N}-\tilde{u}^{n-2N},e_{w}^{n} \right) \nonumber\\
	&+\left( {}_{0}^{C}D_{t}^{\beta}w^{n}-D^{\beta}w^{n}+D^{\beta}w^{n}-D^{\beta}R_{h}w^{n},e_{w}^{n} \right),\label{eq4_7}\\
	p\left( D^{\beta}\nabla e_{\tilde{u}}^{n},\nabla e_{\tilde{u}}^{n} \right)
	=&p\left( \nabla e_{w}^{n},\nabla e_{\tilde{u}}^{n} \right)+p\left( {}_{0}^{C}D_{t}^{\beta}\nabla\tilde{u}^{n}-D^{\beta}\nabla\tilde{u}^{n}+D^{\beta}\nabla\tilde{u}^{n}-D^{\beta}\nabla R_{h}\tilde{u}^{n},\nabla e_{\tilde{u}}^{n} \right).\label{eq4_8}
\end{align}
For \eqref{eq4_7} and \eqref{eq4_8},
following the proof of Lemma \ref{lemma3_3}, if we replace $w_{h}^{n}$ in \eqref{eq3_13} and $\tilde{u}_{h}^{n}$ in \eqref{eq3_14} by $e_{w}^{n}$ and $e_{\tilde{u}}^{n}$, respectively, and take
$
	h_{1}^{n}=a( R_{h}\tilde{u}^{n}-\tilde{u}^{n} )+b( R_{h}\tilde{u}^{n-2N}-\tilde{u}^{n-2N} )+( {}_{0}^{C}D_{t}^{\beta}w^{n}-D^{\beta}w^{n} )+( D^{\beta}w^{n}-D^{\beta}R_{h}w^{n} )
$
and
$
	h_{2}^{n}=( {}_{0}^{C}D_{t}^{\beta}\tilde{u}^{n}-D^{\beta}\tilde{u}^{n} )+( D^{\beta}\tilde{u}^{n}-D^{\beta}R_{h}\tilde{u}^{n}),
$
then we can readily obtain the following inequality
\begin{equation}\label{eq4_9}
	\left( D^{\beta}+\frac{a}{\sqrt{p}} \right)\sqrt{\left\| e_{w}^{n} \right\|^{2}+p\left\| \nabla e_{\tilde{u}}^{n} \right\|^{2}}
	\lesssim\frac{b}{\sqrt{p}}\sqrt{\left\| e_{w}^{n-2N} \right\|^{2}+p\left\| \nabla e_{\tilde{u}}^{n-2N} \right\|^{2}}
	+ \sqrt{\left\| h_{1}^{n} \right\|^{2}+p\left\| \nabla h_{2}^{n} \right\|^{2}}.
\end{equation}

In order to investigate  the convergence of the fully discrete finite element scheme \eqref{eq3_5}--\eqref{eq3_6}, we first present some useful lemmas.
It is worth noting that some results are already available in \cite{Kopteva2025Pointwise} for
the problem without delay. Compared with the case considered there, however, the delay term in the present work will make the analysis of the problem fundamentally different.
For example, the comparison principle used in the above work involves only the current time level and the one immediately preceding it.
% , so that the induction never has to cross a delay interval;
In the present work, by contrast, the comparison principle  involves the data at the current time level and at the time level $2N$ steps earlier, and consequently both the comparison principle and the subsequent conclusions based on it have to be established afresh.
In addition, the presence of the delay term in the present work forces us to devise numerical scheme under the symmetric graded time mesh. In this situation, the family of barrier functions anchored at
$t=0$ cannot be shifted directly and has to be rebuilt on each time interval $((k-1)\tau,k\tau]$,
with $(k-1)\tau$ taken as the local origin, and the relevant properties have to be verified again.
In a word, we believe that it is clear that the following lemmas \ref{lemma4.1}--\ref{lemma4.7} cannot be directly observed from \cite{Kopteva2025Pointwise}. Therefore, for the sake of completeness and readability, We provide detailed proofs of these lemmas following an argument similar to that in \cite{Kopteva2025Pointwise}.
Since the proofs of these lemmas are
somewhat lengthy, we defer the proofs of Lemmas \ref{lemma4.1}--\ref{lemma4.6} to \ref{Proofs of Lemmas 4.1--4.6} for ease of reading and present only the corresponding results here.
\begin{lemma}\label{lemma4.1}
	Let $\bar{\rho}=\max\left\{\rho_n\right\}$. For a given $n$ with
    $2(k-1)N+1\leq n\leq 2kN,$
    $1\leq k\leq K$, if $\lambda_0,\lambda_1\geq 0$, $\lambda_0 \bar{\rho}^\beta<\left\{\Gamma(2-\beta)\right\}^{-1}$ and sequence
	$\{V^i\}$ satisfies
	\begin{empheq}[left=\empheqlbrace]{align}
		&(D^\beta-\lambda_0)V^n-\lambda_1 V^{n-2N}\leq 0,
		\label{eq4_10} \\
		&V^{n-2N}\leq 0,\quad 1\leq n\leq 2N,
		\label{eq4_11}
	\end{empheq}
	then
	\begin{equation}
		V^n\leq 0.
		\label{eq4_12}
	\end{equation}
\end{lemma}

Next, we give some definitions firstly. For fixed $k$ with $1\leq k\leq K,$ when $2(k-1)N+1\leq n\leq 2kN$, denote
\begin{equation}\label{eq4_13}
	\bar{t}_1=t_1,
	\bar{t}_j=\min\left\{t_n:\ t_n-(k-1)\tau\geq
	\left(\frac{C(j-2(k-1)N)}{N}\right)^r\right\},
\end{equation}
and construct a temporal submesh $\{t_j'\}\subset\{t_n\}$ satisfying
$
	t_j'=
	\begin{cases}
		\bar{t}_j, & \bar{t}_j\leq k\tau,\\
		k\tau, & \bar{t}_j>k\tau,
	\end{cases}
$
where $C$ is a sufficiently large constant. Let $\theta=1-\beta$. Define
\begin{equation}
	B(t,t_m')=\min\left\{
	\frac{t-(k-1)\tau}{t_m'-(k-1)\tau}\bigl(t_m'-(k-1)\tau\bigr)^{-\theta},
	\ (t-(k-1)\tau)^{-\theta}
	\right\},
	\label{eq4_14}
\end{equation}
where $(k-1)\tau\leq t\leq k\tau$ and $2(k-1)N+1\leq m\leq 2kN.$
Choose a sufficiently large integer $p$ with $2(k-1)N+1\leq p\leq 2kN$, $1\leq k\leq K$,
and for a given $\gamma \in \mathbb{R}$ and any fixed integer $l$ with $0\leq l\leq L$, $L\in  \mathbb{N}$, define
$
	p_l=2^l(p-2(k-1)N)$,
	$c_l=(\frac{t_p'-(k-1)\tau}{t_{p_l}'})^\gamma$ and $C_L(t)=\sum_{l=0}^{L}c_l\,B(t,t_{p_l}'+(k-1)\tau).$
Now, for fixed $k$ with $1\leq k \leq K,$ we state the following lemmas.
\begin{lemma}\label{lemma4.2}
	Under the symmetric graded temporal mesh \eqref{eq3_1}, the following conclusions hold:
	\begin{enumerate}
		\item[(i)]
		\begin{equation}\label{eq4_15}
			t_j'-(k-1)\tau\approx
			\left(\frac{C(j-2(k-1)N)}{N}\right)^r.
		\end{equation}
		\item[(ii)]
		\begin{empheq}[left={\Gamma(\theta)D_t^\beta B(t,t_m')\geq\empheqlbrace}]{align}
			&\theta\left(\frac{t_m'-(k-1)\tau}{t-(k-1)\tau}\right)^{1-\theta}
			\bigl(t-(k-1)\tau\bigr)^{-1},
			\quad t>t_m', \label{eq4_16} \\
			&\theta^{-1}\left(\frac{t-(k-1)\tau}{t_m'-(k-1)\tau}\right)^\theta
			\bigl(t_m'-(k-1)\tau\bigr)^{-1},
			\quad t\leq t_m'. \label{eq4_17}
		\end{empheq}
		\item[(iii)]
		From $p_l+2(k-1)N<j\leq p_{l+1}+2(k-1)N$, it follows that
		\begin{equation}\label{eq4_18}
			\frac{t_{p_l}'}{t_j'-(k-1)\tau}\geq 2^{-r},\quad
			c_l\simeq\left(\frac{t_p'-(k-1)\tau}{t_j'-(k-1)\tau}\right)^\gamma,
		\end{equation}
		and
		\begin{equation}\label{eq4_19}
			\sum_{l=0}^{L}c_l\lesssim S_L:=
			\begin{cases}
				1, & \gamma>0,\\
				1+\ln\!\left(\dfrac{t_{p_L}'}{t_p'-(k-1)\tau}\right), & \gamma=0,\\
				\left(\dfrac{t_p'-(k-1)\tau}{t_{p_L}'}\right)^\gamma, & \gamma<0.
			\end{cases}
		\end{equation}
	\end{enumerate}
\end{lemma}

\begin{lemma}\label{lemma4.3}
	Let $\gamma^*=\min\left\{\beta,\frac{2-\beta}{r}\right\}$.
	If $\gamma\in \mathbb{R}$ satisfies $\gamma\leq\gamma^*$, then we have
	\begin{empheq}[left={\left|\left.(D_t^\beta-D^\beta) C_L (t)\right|_{t=t_j'}\right|\leq\empheqlbrace}]{align}
		&0, && j\leq p, \label{eq4_20} \\
		&S_J(t_j'-(k-1)\tau)^{-1}
		(p-2(k-1)N)^{\beta-2}
		\left(\dfrac{t_p'-(k-1)\tau}{t_j'-(k-1)\tau}\right)^{\gamma^*},
		&& j>p, \label{eq4_21}
	\end{empheq}
	where $2(k-1)N+1\leq j\leq 2kN$,
	$J:=J(j)$ denotes the maximum value of  set $\{l:p_l+2(k-1)N<j\}$.
\end{lemma}

\begin{lemma}\label{lemma4.4}
	 If $\gamma\in\mathbb{R}$ satisfying $\gamma\le \gamma^*$, then
	\begin{align}
		&C_L(t_n) \lesssim S_L \bigl(t_n-(k-1)\tau\bigr)^{-\theta}\label{eq4_22}
	\end{align}
    and
    \begin{align}
		D^\beta C_L(t_n) \gtrsim \bigl(t_n-(k-1)\tau\bigr)^{-1}
		\left( \frac{t_p'-(k-1)\tau}{t_n-(k-1)\tau} \right)^{\gamma}\label{eq4_23},
	\end{align}
	where
    % $n\leq2(k-1)N + p_{L+1}$,
    $2(k-1)N+1\leq n\leq min\{n\leq2(k-1)N + p_{L+1},2kN\}$, and the constant $C$ hidden in \eqref{eq4_23} depends on $p-2(k-1)N$ but is independent of $L$.
\end{lemma}
\begin{lemma}\label{lemma4.5}
	Let $\left\{V^i \right\}$ be a sequence. If $\gamma\in\mathbb{R}$ and
	\begin{empheq}[left=\empheqlbrace]{align}
		&D^\beta V^n \lesssim \left(t_n-(k-1)\tau\right)^{-1}
		\left( \frac{\rho_1}{t_n-(k-1)\tau} \right)^\gamma,
		\label{eq4_24} \\
		&V^{n-2N} = 0,\qquad 1\le n\le 2N,
		\label{eq4_25}
	\end{empheq}	
	then
	\begin{equation}
		\label{eq4_26}
		V^n \lesssim \left(t_n-(k-1)\tau\right)^{\beta-1}
		l_\gamma^n \left( \frac{\rho_1}{t_n-(k-1)\tau} \right)^{\min\left\{0,\,\gamma\right\}},
	\end{equation}
	where
	\begin{equation}
		\label{eq4_27}
		l_\gamma^n :=
		\begin{cases}
			1 + \ln\left( \dfrac{t_n-(k-1)\tau}{\rho_1} \right), & \gamma=0,\\
			1, & \gamma\ne 0.
		\end{cases}
	\end{equation}
\end{lemma}

\begin{lemma}\label{lemma4.6}
	For  fixed constant $a_0$ and mesh point $t_n$ with
	$0<a_0<\dfrac{2}{5}\left\{\lambda\Gamma(2-\beta)\right\}^{-1/\beta}$ and  $0\le t_n\le K\tau$,
	if $\bar{\rho}\le \dfrac{a_0}{2}$,
    %$2(k-1)N+1\le n\le 2kN$,
	then there exists a nondecreasing discrete sequence $\left\{E^i\right\}$ satisfying $0\le E^i\lesssim1$, $E^i=0$ for  $i\le n$
	and
	\begin{equation*}
		\left(D^\beta - \lambda\right)E^i\gtrsim Q^i :=
		\left\{
		\begin{aligned}
			&0, \quad t_i < t_n + a_0, \\
			&1, \quad t_i \ge t_n + a_0,
		\end{aligned}
		\right.
		\text{\ for\ } 1 \le i \le 2KN.
	\end{equation*}
\end{lemma}

\begin{lemma}\label{lemma4.7}
	If $\gamma\in\mathbb{R}$, $\lambda_0,\lambda_1\ge 0$,
    $\lambda_0\bar{\rho}^{\beta}<\left\{\Gamma(2-\beta)\right\}^{-1}$, and for fixed $n$ with $2\left(k-1\right)N+1\le n\le 2kN$ such that sequence $\{V^{i}\}$ satisfies
	\begin{empheq}[left=\empheqlbrace]{align}
		&\left(D^{\beta}-\lambda_0\right)V^{n}-\lambda_1V^{n-2N}\lesssim\left(\frac{\rho_1}{t_n-\left(k-1\right)\tau}\right)^{1+\gamma},
		\label{eq4_28} \\
		&V^{n-2N} = 0,\quad 1\le n\le 2N.
		\label{eq4_29}
	\end{empheq}
	then
	\begin{equation}\label{eq4_30}
		V^{n}\lesssim\rho_1\,\left(t_n-\left(k-1\right)\tau\right)^{\beta-1}l_\gamma^n
		\left(\frac{\rho_1}{t_n-\left(k-1\right)\tau}\right)^{\min\{0,\gamma\}}.
	\end{equation}
	% where
	% \begin{equation*}
	% 	l_\gamma^n:=\begin{cases}
	% 		1+\ln\dfrac{t_n-\left(k-1\right)\tau}{\rho_1},&\gamma=0,\\
	% 		1,&\gamma\ne 0.
	% 	\end{cases}
	% \end{equation*}
\end{lemma}

\begin{proof}
For convenience, we denote $t_n^{[k]}=t_n-(k-1)\tau$.
	Let sequence $\left\{E_\gamma^{i}\right\}$ satisfy
	\begin{empheq}[left=\empheqlbrace]{align}
		&D^{\beta}E_\gamma^{n}
		=\left(\frac{\rho_1}{t_n^{[k]}}\right)^{\gamma+1},
		\label{eq4_31} \\
		&E_\gamma^{n-2N}=0,\qquad 1\le n\le 2N.
		\label{eq4_32}
	\end{empheq}
	According to Lemma \ref{lemma4.5}, it is clear that
	\begin{equation}\label{eq4_33}
		0\leq E_\gamma^{n}\lesssim\left(t_n^{[k]}\right)^{\beta-1}l_\gamma^n
		\left(\frac{\rho_1}{t_n^{[k]}}\right)^{\min\{0,\gamma\}}\rho_1\\
		=\rho_1^{\beta}l_\gamma^n
		\left(\frac{\rho_1}{t_n^{[k]}}\right)^{\min\{0,\gamma\}+1-\beta}.
	\end{equation}
    In addition, \eqref{eq4_31} also implies that
	\begin{equation}\label{eqs_0914A}
		D^{\beta}E_\gamma^{n-2N}=\left(\frac{\rho_1}{t_{n-2N}^{[k-1]}}\right)^{\gamma+1}=\left(\frac{\rho_1}{t_n^{[k]}}\right)^{\gamma+1}=D^{\beta}E_\gamma^{n}.
	\end{equation}
    When $1\le n\le 2N$, noting that $0=E_\gamma^{n-2N}\le E_\gamma^n$, it follows from Lemma \ref{lemma4.1} and \eqref{eqs_0914A} that
    $E_\gamma^{n-2N}\le E_\gamma^n$.
    Let
	\begin{equation}\label{eq4_34}
		\lambda = \lambda_0 + \lambda_1,\quad\bar{\gamma} = \min \{0, \gamma \} -\beta,\quad \bar{l_\gamma} = \max_{2(k-1)N+1 \leq n \leq 2kN} l_\gamma^n=l_\gamma^{2kN}.
	\end{equation}
	and choose a sufficiently large constant $C_\gamma$ such that $\lambda_0E_\gamma^{n}+\lambda_1E_\gamma^{n-2N}\le\lambda C_\gamma E_\gamma^{n}$.
	By combining with \eqref{eq4_31}, \eqref{eq4_33} and \eqref{eq4_34}, one has
	\begin{equation}\label{eq4_35}
		%\begin{aligned}
			D^{\beta}E_\gamma^{n}-\left(\lambda_0E_\gamma^{n}+\lambda_1E_\gamma^{n-2N}\right)
			\ge D^{\beta}E_\gamma^{n}-\lambda C_\gamma E_\gamma^{n}
			\ge\left(\frac{\rho_1}{t_n^{[k]}}\right)^{\gamma+1}-\lambda C_\gamma\rho_1^{\beta}\bar l_\gamma
			\left(\frac{\rho_1}{t_n^{[k]}}\right)^{\min\{0,\gamma\}+1-\beta}.
		%\end{aligned}
	\end{equation}
    If $\gamma$ is set equal to $\bar{\gamma}$ in \eqref{eq4_35}, then, noting that $\bar l_{\bar\gamma}=1$ (as $\bar\gamma\neq0)$ and $\min\{0,\bar{\gamma}\}-\beta=\bar{\gamma}-\beta$, we immediately obtain
	\begin{equation}\label{eq4_36}
		\begin{aligned}
			\left(D^{\beta}-\lambda_0\right)E_{\bar\gamma}^{n}-\lambda_1E_{\bar\gamma}^{n-2N}
			&\ge\left(\frac{\rho_1}{t_n^{[k]}}\right)^{\bar\gamma+1}
			-\lambda C_{\bar\gamma} \rho_1^{\beta}
			\left(\frac{\rho_1}{t_n^{[k]}}\right)^{\bar\gamma+1-\beta}
			=\frac12\left(\frac{\rho_1}{t_n^{[k]}}\right)^{\bar\gamma+1}
			+\frac12\left(\frac{\rho_1}{t_n^{[k]}}\right)^{\bar\gamma+1}
			\Big[1-2\lambda C_{\bar\gamma}\left(t_n^{[k]}\right)^{\beta}\Big].
		\end{aligned}
	\end{equation}

    It is worth noting that, for given time mesh \eqref{eq3_1}, if there exists $n'$ such that time step-size $\rho_{n'}>a_0/2 $, then it follows from \eqref{eq3_2} that $a_0/2<\rho_{n'}\lesssim \rho_1^{1/r}(t_{n'}-(k-1)\tau)^{1-1/r}$, which implies $\rho_1^{1/r}\gtrsim (t_{n'}-(k-1)\tau)^{1/r-1}$ and hence $N\simeq 1$.
    Thus, it is clear that \eqref{eq4_30} holds directly in this case. Now we consider the case for $\bar{\rho}\le a_0/2.$
    First, in this case, Lemma \ref{lemma4.6} indicates
    we can construct a discrete nondecreasing sequence $\{E^i\}$ satisfying $E^i=0$ for $-2N\leq i<0$,
   $0\le E^{i}\lesssim 1,$
	\begin{equation}\label{eq4_38}
		 E^{i}=0,\ 0\leq i\leq m,\ \text{and} \ 	
		\left(D^{\beta}-\lambda\right) E^{i}\gtrsim Q^{i}:=
		\begin{cases}
			0,&t_i<t_m+a_0,\\
			1,&t_i\ge t_m+a_0,
		\end{cases}
		\quad 1\le i\le 2KN,
	\end{equation}
	where $m$ is chosen such that $t_m$ as the largest mesh point satisfying $t_m\le (k-1)\tau+a_0$. Since $\lambda_0E^{n}+\lambda_1E^{n-2N}\le\lambda E^{n}$, \eqref{eq4_38} gives
	\begin{equation}\label{eq4_39}
		\left(D^{\beta}-\lambda_0\right)E^{n}-\lambda_1E^{n-2N}
		\ge D^{\beta}E^{n}-\lambda E^{n}\gtrsim Q^{n}.
	\end{equation}
	Let
	\begin{equation}\label{eq4_40}
		F^{n}=E_\gamma^{n}+2\lambda C_\gamma\bar l_\gamma\rho_1^{\beta}
		\big[E_{\bar\gamma}^{n}+C_{*}\rho_1^{1+\bar\gamma}E^{n}\big],
	\end{equation}
	where $C_{*}$ is a sufficiently large constant. The combination of \eqref{eq4_35}, \eqref{eq4_36} and \eqref{eq4_39} leads to
	\begin{equation}\label{eq4_42}
		\begin{aligned}
			\left(D^{\beta}-\lambda_0\right)F^{n}-\lambda_1F^{n-2N}
			% =&\left(D^{\beta}-\lambda_0\right)E_\gamma^{n}-\lambda_1E_\gamma^{n-2N}
			% +2\lambda C_\gamma\bar l_\gamma\rho_1^{\beta}
			% \left[\left(D^{\beta}-\lambda_0\right)E_{\bar\gamma}^{n}-\lambda_1E_{\bar\gamma}^{n-2N}\right]\\
			% & +2\lambda C_\gamma\bar l_\gamma\rho_1^{\beta}C_{*}\rho_1^{1+\bar\gamma}
			% \left[\left(D^{\beta}-\lambda_0\right)E^{n}-\lambda_1E^{n-2N}\right]\\
			\gtrsim&\left(\frac{\rho_1}{t_n^{[k]}}\right)^{1+\gamma}
			-\lambda C_{\gamma}\bar l_\gamma \rho_1^{\beta}\left(\frac{\rho_1}{t_n^{[k]}}\right)^{1+\bar\gamma}
			+2\lambda C_{\gamma}\bar l_\gamma \rho_1^{\beta}\\
			& \cdot\left\{
			\frac12\left(\frac{\rho_1}{t_n^{[k]}}\right)^{1+\bar\gamma}
			+\frac12\left(\frac{\rho_1}{t_n^{[k]}}\right)^{1+\bar\gamma}
			\left[1-2\lambda C_{\bar\gamma}\left(t_n^{[k]}\right)^\beta\right]
			\right\}
			 +2\lambda C_{\gamma}\bar l_\gamma \rho_1^{\beta}C_*\rho_1^{1+\bar\gamma}Q^n\\
            \gtrsim&\left(\frac{\rho_1}{t_n^{[k]}}\right)^{1+\gamma}
			+\lambda C_{\gamma}\bar l_\gamma\rho_1^{\beta}
			\left\lbrace
			\left(\frac{\rho_1}{t_n^{[k]}}\right)^{1+\bar\gamma}
			\left[1-2\lambda C_{\bar\gamma}\left(t_n^{[k]}\right)^\beta\right]
			+2C_*\rho_1^{1+\bar\gamma}Q^n
			\right\rbrace.
		\end{aligned}
	\end{equation}		
	% After rearrangement, \eqref{eq4_42} can be rewritten as
	% \begin{equation}\label{eq4_43}
	% 	\begin{aligned}
	% 		\left(D^{\beta}-\lambda_0\right)F^n-\lambda_1F^{n-2N}
	% 		\gtrsim&\left(\frac{\rho_1}{t_n-\left(k-1\right)\tau}\right)^{1+\gamma}
	% 		+\lambda C_{\gamma}\bar l_\gamma\rho_1^{\beta}\\
	% 		&\times\left\lbrace
	% 		\left(\frac{\rho_1}{t_n-\left(k-1\right)\tau}\right)^{1+\bar\gamma}
	% 		\Big[1-2\lambda C_{\bar\gamma}\left(t_n-\left(k-1\right)\tau\right)^\beta\Big]
	% 		+2C_*\rho_1^{1+\bar\gamma}Q^n
	% 		\right\rbrace.
	% 	\end{aligned}
	% \end{equation}
    From \eqref{eq4_42}, the fact that
	\begin{equation}\label{eq4_41}
		\left(D^{\beta}-\lambda_0\right)F^{n}-\lambda_1F^{n-2N}
		\gtrsim\left(\frac{\rho_1}{t_n^{[k]}}\right)^{1+\gamma}
	\end{equation}
   can be demonstrated by the following two cases:
	If $1-2\lambda C_{\bar\gamma}(t_n^{[k]})^\beta\ge 0$, then \eqref{eq4_41} obviously holds;
	If $1-2\lambda C_{\bar\gamma}(t_n^{[k]})^\beta<0$, then $t_n^{[k]}>(2\lambda C_{\bar\gamma})^{-1/\beta}$. In addition, since $a_0$ is a sufficiently small constant, we have $2\lambda C_{\bar\gamma}(2a_0)^\beta\le 1$. Therefore, it is clear that $t_n^{[k]}> 2a_0\geq t_m-(k-1)\tau+a_0$, i.e., $t_n>t_m+a_0$. According to  \eqref{eq4_38}, the above result indicates that $Q^n=1$. Moreover, since $C_*$ is a sufficiently large constant, we have $C_*\geq\lambda C_{\bar\gamma}(t_n^{[k]})^{\beta-1-\bar\gamma}$, i.e., 	
	$
		2C_*\rho_1^{1+\bar\gamma}Q^n
		+\left(\rho_1/(t_n^{[k]}\right)^{1+\bar\gamma}
		[1-2\lambda C_{\bar\gamma}(t_n^{[k]})^\beta]\ge 0.
	$
	Thus, it follows from the foregoing discussion that \eqref{eq4_41} holds under the case $1-2\lambda C_{\bar\gamma}(t_n^{[k]})^\beta< 0$.
	
	Now, by combining \eqref{eq4_41} and \eqref{eq4_28}, we can obtain
	$
		(D^{\beta}-\lambda_0)V^n-\lambda_1V^{n-2N}
		\lesssim (D^{\beta}-\lambda_0)F^n-\lambda_1F^{n-2N}.
	$
	Since $E_{\gamma}^{n-2N}=E_{\bar\gamma}^{n-2N}=E^{n-2N}=0$ for $ 1\le n\le 2N$ implies  $F^{n-2N}=V^{n-2N}=0$ for $ 1\le n\le 2N$, Lemma \ref{lemma4.1} gives
	\begin{equation}\label{eq4_44}
		V^n\le F^n
		\lesssim E_{\gamma}^n+\bar{l_{\gamma}}\rho_1^{\beta}
\left[E_{\bar\gamma}^n+\rho_1^{1+\bar\gamma}E^n\right].
	\end{equation}
    Based on \eqref{eq4_33} and $0\le E^n\lesssim 1$, the fact that
    \begin{equation}\label{eq4_45}
		V^n\lesssim\rho_1^{\beta}l_\gamma^n
		\left(\frac{\rho_1}{t_n^{[k]}}\right)^{\min\{0,\gamma\}+1-\beta}+\bar l_\gamma \rho_1^{\beta}
		\left[\rho_1^{\beta}\left(\frac{\rho_1}{t_n^{[k]}}\right)^{1+\bar\gamma-\beta}\right]
	\end{equation}
    can be concluded by the following two cases: If $E^n= 0$, then \eqref{eq4_44} can be rewritten into \eqref{eq4_45} obviously;
	If $E^n\neq 0$, then \eqref{eq4_44} can also be easily rewritten as \eqref{eq4_45} by
	$
\rho_1^\beta\left(\rho_1/t_n^{[k]}\right)^{1+\bar\gamma-\beta}
=\rho_1^{1+\bar\gamma}\left(t_n^{[k]}\right)^{2\beta-1-\min\{0,\gamma\}}
		\simeq \rho_1^{1+\bar\gamma}E^n.
	$
	Since
	\begin{equation}\label{eq4_46}
		\begin{aligned}
			\bar{l}_\gamma\left(t_n^{[k]}\right)^\beta
			% &=l_\gamma^{2kN}\left(t_n-(k-1)\tau\right)^\beta
			&=\left[1+\ln\left(\frac{\tau}{\rho_1}\right)\right]\left(t_n^{[k]}\right)^\beta
			=\left(1+\ln\frac{t_n^{[k]}}{\rho_1}+\ln\frac{\tau}{t_n^{[k]}}\right)\left(t_n^{[k]}\right)^\beta
			\lesssim l_\gamma^n,
		\end{aligned}
	\end{equation}	
    it follows from
	\eqref{eq4_45} and \eqref{eq4_46} that
	\begin{equation}\label{eq4_47}
		\begin{aligned}
			V^n
			\lesssim& \rho_1^\beta l_\gamma^n\left(\frac{\rho_1}{t_n^{[k]}}\right)^{\min\{0,\gamma\}+1-\beta}+l_\gamma^n\left(t_n^{[k]}\right)^{-\beta}
			\rho_1^\beta\left[\left(\frac{\rho_1}{t_n^{[k]}}\right)^{1+\bar\gamma-\beta}\rho_1^\beta\right]
			= \rho_1^\beta l_\gamma^n\left(\frac{\rho_1}{t_n^{[k]}}\right)^{\min\{0,\gamma\}+1-\beta}.
  %           \lesssim\rho_1\,\left(t_n-\left(k-1\right)\tau\right)^{\beta-1}l_\gamma^n
		% \left(\frac{\rho_1}{t_n-\left(k-1\right)\tau}\right)^{\min\{0,\gamma\}}.
		\end{aligned}
	\end{equation}	
    The above result clearly show that
	  \eqref{eq4_30} holds.
\end{proof}
\begin{lemma}\label{lemma4.8}
	Suppose that $(\tilde{u}(x,t),w(x,t))$ is the solution of the problem \eqref{eq2_36b1}--\eqref{eq2_36b6} satisfying the regularity conditions \eqref{eq1_6} and \eqref{eq2_42}. Then, for fixed $x\in[0,L]$ under the time mesh \eqref{eq3_1}, we have
	\begin{equation}\label{eq4_48}
		\left|\left(D^\beta-_0^CD_t^\beta\right)w^n\right|
		\lesssim \sum_{l=1}^{k}\left(\frac{\rho_1}{t_n-\left(l-1\right)\tau}\right)^{\frac{2-\beta}{r}}
		\left(t_n-\left(l-1\right)\tau\right)^{\left(2l-2\right)\beta},
	\end{equation}
	\begin{equation}\label{eq4_49}
		\left|\left(D^\beta-_0^CD_t^\beta\right)\tilde{u}^n\right|
		\lesssim \sum_{l=1}^{k}\left(\frac{\rho_1}{t_n-\left(l-1\right)\tau}\right)^{\frac{2-\beta}{r}}
		\left(t_n-\left(l-1\right)\tau\right)^{\left(2l-1\right)\beta},
	\end{equation}
	where $2\left(k-1\right)N+1\le n\le 2kN$, $1\le k\le K$.
\end{lemma}

\begin{proof}
	Similar to \cite{Bu2024Finite}, we define
	$
		\sigma_{2(l-1)N+1}=\sup_{v\in(t_{2(l-1)N}, t_{2(l-1)N+1})}
		(v-(l-1)\tau)^{1-\beta}
		|\frac{\nabla_tw_{2(l-1)N+1}}{\rho_{2(l-1)N+1}}-\partial_vw(v)|, 1\le l\le k,
	$
	and for $i\in\Omega_1:=\{2(l-1)N+1\le i\le 2lN,\ 1\le l\le k-1\}\cup\{2(k-1)N+1\le i\le n,\ l=k\},$ denote
	$
		\sigma_i'=(t_i^{[l]})^{2-\beta}\sup_{v\in(t_{i-1},t_i)}|\partial_v^2w(v)|.
	$
	% \begin{equation*}
	% 	i\in\Omega_1:=\left\{2\left(l-1\right)N+1\le i\le 2lN,\ 1\le l\le k-1\right\}\cup\left\{2\left(k-1\right)N+1\le i\le n,\ l=k\right\}.
	% \end{equation*}
    Let $w_L$ be a piecewise linear Lagrange interpolant of $w$ on time mesh \eqref{eq3_1}  and $\xi=w-w_L$.
	Then, according to \cite[Lemma 5.1, Lemma 5.2]{Bu2024Finite}, it holds
    \begin{equation}\label{eqs4_0916A}
    \sigma_{2(l-1)N+1}\lesssim \rho_1^{\left(2l-2\right)\beta}, 1\le l\le k \text{ and } \sigma_i'\lesssim \left(t_i^{[l]}\right)^{\left(2l-2\right)\beta},  i\in\Omega_1,
    \end{equation}
	%Let $w_L$ be a piecewise linear Lagrange interpolant of $w$ on time mesh \eqref{eq3_1}  and $\xi=w-w_L$.
	\begin{equation}\label{eq4_50}
			|\xi(v)|
			\lesssim \left(t_{2(l-1)N+1}-v\right)\left(v-(l-1)\tau\right)^{\beta-1}\sigma_{2(l-1)N+1},\quad v\in(t_{2(l-1)N},t_{2(l-1)N+1}),\ 1\le l\le k
	\end{equation}
    and \begin{equation}\label{eq4_51}
		|\xi(v)|
		\lesssim \rho_i^2\left(t_i^{[l]}\right)^{\beta-2}\sigma_i',\quad v\in(t_{i-1},t_i),\ i\in\Omega_1.
	\end{equation}
    %When $v\in(t_{2(l-1)N},t_{2(l-1)N+1})$,$\ 1\le l\le k$, we have
	%\begin{equation}\label{eq4_50}
		%\begin{aligned}
			%|\xi(v)|&=\left|\int_{v}^{t_{2(l-1)N+1}}\xi'(\eta)\,\mathrm{d}\eta\right|
			%\leq \int_{v}^{t_{2(l-1)N+1}}
			%\left|\partial_\eta w(\eta)-\frac{\nabla_tw_{2(l-1)N+1}}{\rho_{2(l-1)N+1}}\right|\mathrm{d}\eta
			%\lesssim \bigl(t_{2(l-1)N+1}-v\bigr)\bigl(v-(l-1)\tau\bigr)^{\beta-1}\sigma_{2(l-1)N+1}.
		%\end{aligned}
	%\end{equation}
	%When $v\in(t_{i-1},t_i),\ i\in\Omega_1$, the error estimate theory of linear Lagrange interpolation demonstrates
	%\begin{equation}\label{eq4_51}
		%|\xi(v)|\lesssim \rho_i^2\sup_{v\in(t_{i-1},t_i)}\left|\partial_v^2w(v)\right|
		%\lesssim \rho_i^2\bigl(t_i-(l-1)\tau\bigr)^{\beta-2}\sigma_i'.
	%\end{equation}
   Since $\left({}^C_0D_t^\beta-D^\beta\right)w^n
			=\frac{\beta}{\Gamma\left(1-\beta\right)}\int_{0}^{t_n}\left(t_n-v\right)^{-\beta-1}\xi(v)\,\mathrm{d}v$,
	from \eqref{eq4_50} and \eqref{eq4_51}, we have
	\begin{equation}\label{eq4_52}
		\begin{aligned}
			\left|\left({}^C_0D_t^\beta-D^\beta\right)w^n\right|
			\lesssim& \sum_{l=1}^{k}\int_{t_{2(l-1)N}}^{t_{2(l-1)N+1}}
			\left(t_n-v\right)^{-\beta-1}\left(v-\left(l-1\right)\tau\right)^{\beta-1}
			\left(t_{2(l-1)N+1}-v\right)\sigma_{2(l-1)N+1}\mathrm{d}v\\
			&+\sum_{l=1}^{k}\sum_{{i=2(l-1)N+2}}^{2lN}\int_{t_{i-1}}^{t_i}
			\left(t_n-v\right)^{-\beta-1}\left(\rho_i\right)^2\left(t_i^{[l]}\right)^{\beta-2}\sigma_i'\mathrm{d}v\\
			&+\sum_{{i=2(k-1)N+2}}^{n}\int_{t_{i-1}}^{t_i}
			\left(t_n-v\right)^{-\beta-1}\left(\rho_i\right)^2\left(t_i^{[k]}\right)^{\beta-2}
			\min\left\{1,\frac{t_n-v}{\rho_n}\right\}\sigma_i'\mathrm{d}v
			:=I+II+III.
		\end{aligned}
	\end{equation}

	Next, we estimate $I$, $II$ and $III$ separately. For $I$, noting that if $v\in(t_{2(l-1)N},t_{2(l-1)N+1})$, then
	$
		(t_{2(l-1)N+1}-v)/(t_n-v)\le (t_{2(l-1)N+1}-(l-1)\tau)/(t_n-\left(l-1\right)\tau).
	$
	Therefore, it holds
	$
		\left(t_n-v\right)^{-\beta-1}
		\le[(t_n-(l-1)\tau)(t_{2(l-1)N+1}-v)/\rho_1]^{-\beta-1}.
	$
    By using this inequality and \eqref{eqs4_0916A}, we have
	% where $\rho_1=t_{2(l-1)N+1}-\left(l-1\right)\tau$.
    % Combining with $\sigma_{2(l-1)N+1}\lesssim \rho_1^{\left(2l-2\right)\beta},\ 1\le l\le k$, the term $I$ can be rewritten as
	\begin{equation}\label{eq4_53}
		\begin{aligned}
			I&\lesssim \sum_{l=1}^{k}\left[\int_{t_{2(l-1)N}}^{t_{2(l-1)N+1}}
			\left(t_{2(l-1)N+1}-v\right)^{-\beta}\left(v-\left(l-1\right)\tau\right)^{\beta-1}\mathrm{d}v\right]
			\times\frac{\rho_1^{\left(2l-1\right)\beta+1}}{\left(t_n^{[l]}\right)^{\beta+1}}\\
	&=\sum_{l=1}^{k}\int_{0}^{\rho_1}\left(\rho_1-x\right)^{-\beta}x^{\beta-1}\mathrm{d}x
			\cdot\frac{\rho_1^{\left(2l-1\right)\beta+1}}{\left(t_n^{[l]}\right)^{\beta+1}}
			% &=\sum_{l=1}^{k}\int_{0}^{1}\rho_1^{1-\beta+\beta-1}\left(1-\widetilde{x}\right)^{-\beta}\widetilde{x}^{\beta-1}\mathrm{d}\widetilde{x}
			% \cdot\frac{\rho_1^{\left(2l-1\right)\beta+1}}{\left(t_n-\left(l-1\right)\tau\right)^{\beta+1}}\\
			\lesssim \sum_{l=1}^{k}\left(\frac{\rho_1}{t_n^{[l]}}\right)^{\left(2l-1\right)\beta+1}
			\left(t_n^{[l]}\right)^{\left(2l-2\right)\beta}.
		\end{aligned}
	\end{equation}
	% where \eqref{eq4_53} uses the substitutions $x=v-(l-1)\tau$ and $\widetilde{x}=\dfrac{x}{\rho_1}$.
	For $II$, it follows from $t_i-\left(l-1\right)\tau\simeq v-\left(l-1\right)\tau,\ v\in(t_{i-1},t_i)$ and \eqref{eqs4_0916A},
    %$\sigma_i'\lesssim \left(t_i-\left(l-1\right)\tau\right)^{\left(2l-2\right)\beta},\ i\in\Omega_1$
    that
	\begin{equation}\label{eq4_54}
		\begin{aligned}
			II\lesssim& \sum_{l=1}^{k-1}\int_{t_{2(l-1)N+1}}^{t_{2lN}}
			\left(t_n-v\right)^{-\beta-1}\left(v-\left(l-1\right)\tau\right)^{-\frac{\beta}{r}-1}\mathrm{d}v\times\max_{i\in\Omega_1}\left\{\left(\rho_i\right)^2\left(t_i^{[l]}\right)^{\beta-2+\frac{\beta}{r}+1+\left(2l-2\right)\beta}\right\}\\
			\lesssim& \sum_{l=1}^{k-1}\left(\frac{\rho_1}{t_n^{[l]}}\right)^{-\frac{\beta}{r}}
			\left(t_n^{[l]}\right)^{-\beta-\frac{\beta}{r}-1}
			\times\max_{i\in\Omega_1}\left\{\rho_1^\frac{2}{r}\left(t_i^{[l]}\right)^{\beta-2+\frac{\beta}{r}+1+\left(2l-2\right)\beta+2-\frac{2}{r}}\right\}\\
			\lesssim& \sum_{l=1}^{k-1}\left(\frac{\rho_1}{t_n^{[l]}}\right)^{\frac{2-\beta}{r}}
			\left(t_n^{[l]}\right)^{\left(2l-2\right)\beta},
		\end{aligned}
	\end{equation}
	where the facts
	$
		\int_{t_{2(l-1)N+1}}^{t_{2lN}}
		(t_n-v)^{-\beta-1}(v-\left(l-1\right)\tau)^{-\beta/r-1}\mathrm{d}v
		\lesssim(\frac{\rho_1}{t_n-(l-1)\tau})^{-\beta/r}
		(t_n-\left(l-1\right)\tau)^{-\beta-\beta/r-1}
	$
	and $(\rho_i)^2\lesssim \rho_1^{2/r}\left(t_i-(l-1\right)\tau)^{2-2/r}$ are used (see \cite[Theorem 7]{Tan2022L1}). Analogous to \eqref{eq4_54},
	$III$ can be bounded by
	\begin{equation}\label{eq4_55}
		\begin{aligned}
			III&\lesssim \int_{t_{2(k-1)N+1}}^{t_n}
			\left(t_n-v\right)^{-\beta-1}\left(v-\left(k-1\right)\tau\right)^{-\frac{\beta}{r}-1}
			\min\left\{1,\frac{t_n-v}{\rho_n}\right\}\mathrm{d}v\times\max_{\substack{2(k-1)N+1\le i\le n}}
			\left\{\left(\rho_i\right)^2\left(t_i^{[k]}\right)^{\beta-2+\frac{\beta}{r}+1+\left(2k-2\right)\beta}\right\}\\
			&\lesssim \left(\frac{\rho_1}{t_n^{[k]}}\right)^{\frac{2-\beta}{r}}
			\left(t_n^{[k]}\right)^{\left(2k-2\right)\beta},
		\end{aligned}
	\end{equation}
	where the fact
	$
		\int_{t_{2(k-1)N+1}}^{t_n}
		(t_n-v)^{-\beta-1}(v-\left(k-1\right)\tau)^{-\beta/r-1}
		\min\{1,\frac{t_n-v}{\rho_n}\}\mathrm{d}v
		\lesssim(\rho_1/t_n^{[k]})^{-\beta/r}
		(t_n^{[k]})^{-\beta-\beta/r-1}
	$
	is used (see  \cite[Theorem 7]{Tan2022L1}).
	Since $\min\left\{\left(2l-1\right)\beta+1,(2-\beta)/r\right\}=(2-\beta)/r$, The combination of \eqref{eq4_52}, \eqref{eq4_54} and \eqref{eq4_55} immediately yields \eqref{eq4_48}. Following a similar argument as in the proof of \eqref{eq4_48}, we can also obtain \eqref{eq4_49}.	
\end{proof}

Now we state the convergence result of the proposed numerical scheme.
\begin{theorem}\label{theorem4.9}
	Assume that $\left(w(x,t),\widetilde{u}(x,t)\right)$ is the solution of the problem \eqref{eq2_36b1}--\eqref{eq2_36b6} satisfying the regularity conditions \eqref{eq1_6} and \eqref{eq2_42}. Let $\frac{-a}{\sqrt{p}} \bar{\rho}^\beta<\left\{\Gamma(2-\beta)\right\}^{-1}$. Then, under the time mesh \eqref{eq3_1}, the numerical solution $\left(w^n_h,\widetilde{u}^n_h\right)$ of \eqref{eq3_5}--\eqref{eq3_6} satisfies
	\begin{equation}\label{eq4_56}
		\left\|w^n_h-w^n\right\|+p\left\|\nabla\widetilde{u}^n_h-\nabla\widetilde{u}^n\right\|
		\lesssim
		\begin{cases}
			&\mathcal{E}_1^n+\rho_1^{\frac{2-\beta}{r}} t_n^{2\beta-\frac{2-\beta}{r}} +h^q, \quad k=1,\\
		&\rho_1^{\frac{2-\beta}{r}}+h^q, \quad k\ge 2,
       % \displaystyle\sum_{l=1}^{k}\rho_1^{\frac{2-\beta}{r}}\left(t_n-\left(l-1\right)\tau\right)^{\left(2l-1\right)\beta-\frac{2-\beta}{r}}
			%+\sum_{l=1}^{k}\rho_1^{\frac{2-\beta}{r}}\left(t_n-\left(l-1\right)\tau\right)^{2l\beta-\frac{2-\beta}{r}}\\
			%&+h^q\left(t_n-\left(k-1\right)\tau\right)^{\beta}+h^q, \quad k\ge 2,
		\end{cases}
	\end{equation}
	where $2(k-1)N+1\le n\le 2kN,\ 1\le k\le K$, and
	$
		\mathcal{E}_1^n:=
		\begin{cases}
			\rho_1^{\frac{2-\beta}{r}}t_n^{\beta-\frac{2-\beta}{r}}, & \dfrac{2-\beta}{r}<1,\\
			\rho_1 t_n^{\beta-1}\left(1+\ln\frac{t_n}{\rho_1}\right), & \dfrac{2-\beta}{r}=1,\\
			\rho_1t_n^{\beta-1}, & \dfrac{2-\beta}{r}>1.
		\end{cases}
	$
\end{theorem}

\begin{proof}
	Take $R_1^n=\left({}^C_0D_t^\beta-D^\beta\right)w^n$, $R_2^n=\left({}^C_0D_t^\beta-D^\beta\right)\nabla\widetilde{u}^n$, $R_3^n=a\left(R_h\widetilde{u}^n-\widetilde{u}^n\right)$,
	$R_4^n=b\left(R_h\widetilde{u}^{n-2N}-\widetilde{u}^{n-2N}\right)$, $R_5^n=D^\beta w^n-D^\beta R_h w^n$ and $R_6^n=D^\beta \nabla\widetilde{u}^n-D^\beta\nabla R_h\widetilde{u}^n$.
	Then \eqref{eq4_9} can be rewritten as
	\begin{equation}\label{eq4_57}
		\left(D^\beta+\frac{a}{\sqrt{p}}\right)\sqrt{\left\|e_w^n\right\|^2+p\left\|\nabla e_{\widetilde{u}}^n\right\|^2}-\frac{b}{\sqrt{p}}\sqrt{\left\|e_w^{n-2N}\right\|^2+p\left\|\nabla e_{\widetilde{u}}^{n-2N}\right\|^2}
		\lesssim \sum_{i=1}^{6}\left\|R_i^n\right\|.
	\end{equation}
	Let the sequences $\left\{E_1^n\right\}, \left\{E_2^n\right\},\cdots,\left\{E_6^n\right\}$ satisfy
	\begin{empheq}[left=\empheqlbrace]{align}
		&\left(D^\beta+\frac{a}{\sqrt{p}}\right)E_i^n-\frac{b}{\sqrt{p}}E_i^{n-2N}=\left\|R_i^n\right\|,
		\label{eq4_58} \\
		&E_i^{n-2N}=0,\quad 1\le n\le 2N,
		\label{eq4_59}
	\end{empheq}
	where $i=1,2,\cdots,6.$ Now we consider the case $i=1$ firstly.
	By \eqref{eq4_48}, one has
	\begin{equation}\label{eq4_60}
		\left(D^\beta+\frac{a}{\sqrt{p}}\right)E_1^n-\frac{b}{\sqrt{p}}E_1^{n-2N}
		\lesssim \sum_{l=1}^{k}\left(\frac{\rho_1}{t_n^{[l]}}\right)^{\frac{2-\beta}{r}-\left(2l-2\right)\beta}\rho_1^{\left(2l-2\right)\beta}
		:=\sum_{l=1}^{k}r_l^n.
	\end{equation}
	When $k=1$, i.e., $1\leq n\leq 2N$, \eqref{eq4_60} implies
	\begin{equation}\label{eq4_61}
		\left(D^\beta+\frac{a}{\sqrt{p}}\right)E_1^n-\frac{b}{\sqrt{p}}E_1^{n-2N}
		\lesssim \left(\frac{\rho_1}{t_n}\right)^{\frac{2-\beta}{r}}.
	\end{equation}
	Therefore, applying Lemma \ref{lemma4.7} with $\gamma=\frac{2-\beta}{r}-1$, we obtain
	\begin{equation}\label{eq4_62}
		E_1^n\lesssim \mathcal{E}_1^n
		=
		\begin{cases}
			\rho_1^{\frac{2-\beta}{r}} t_n^{\beta-\frac{2-\beta}{r}}, & \dfrac{2-\beta}{r}<1,\\
			\rho_1 t_n^{\beta-1}\left(1+\ln\dfrac{t_n}{\rho_1}\right), & \dfrac{2-\beta}{r}=1,\\
			\rho_1 t_n^{\beta-1}, & \dfrac{2-\beta}{r}>1.
		\end{cases}
	\end{equation}
	When $k\ge2$, define sequences $\left\{e_1^n\right\}, \left\{e_2^n\right\}, \cdots, \left\{e_k^n\right\}$ satisfying
	\begin{empheq}[left=\empheqlbrace]{align}
		&\left(D^\beta+\dfrac{a}{\sqrt{p}}\right)e_l^n-\dfrac{b}{\sqrt{p}}e_l^{n-2N}=r_l^n,
		\label{eq4_63} \\
		&e_l^{n-2N}=0,\quad 1\le n\le 2N,
		\label{eq4_64}
	\end{empheq}
    where $l=1,2,\cdots,k.$
	If $1\le l\le k-1$, then $r_l^n=\left(\rho_1/t_n^{[l]}\right)^{(2-\beta)/r-(2l-2)\beta}\rho_1^{(2l-2)\beta}$.
	Applying Lemma \ref{lemma4.7} with $\gamma=-1$ to \eqref{eq4_63} and \eqref{eq4_64} yields
	\begin{equation}\label{eq4_65}
		\begin{aligned}
			e_l^n&\lesssim \rho_1\left(t_n^{[k]}\right)^{\beta-1}
			\left(\frac{\rho_1}{t_n^{[k]}}\right)^{-1} r_l^n
			\lesssim \rho_1^{\frac{2-\beta}{r}} \left(t_n^{[l]}\right)^{(2l-2)\beta-\frac{2-\beta}{r}}\left(t_n^{[k]}\right)^\beta
            \lesssim\rho_1^{\frac{2-\beta}{r}}\left(t_n^{[l]}\right)^{(2l-1)\beta-\frac{2-\beta}{r}}.
		\end{aligned}
	\end{equation}
	If $l=k$, then $r_k^n=\left(\rho_1/t_n^{[k]}\right)^{(2-\beta)/r-(2k-2)\beta}\rho_1^{(2k-2)\beta}$. Applying Lemma \ref{lemma4.7} with $\gamma=(2-\beta)/r-(2k-2)\beta-1$ to \eqref{eq4_63} and \eqref{eq4_64}, and noting that
	  $\gamma<0$, we have
	\begin{equation}\label{eq4_66}
		\begin{aligned}
			e_k^n \lesssim \rho_1\left(t_n^{[k]}\right)^{\beta-1}
			\left(\frac{\rho_1}{t_n^{[k]}}\right)^{\frac{2-\beta}{r}-(2k-2)\beta-1}\rho_1^{(2k-2)\beta}
			=\rho_1^{\frac{2-\beta}{r}} \left(t_n^{[k]}\right)^{(2k-1)\beta-\frac{2-\beta}{r}}.
		\end{aligned}
	\end{equation}
	Now, it follows from \eqref{eq4_60} and \eqref{eq4_63} that
   $
		\left(D^\beta+\frac{a}{\sqrt{p}}\right)E_1^n-\frac{b}{\sqrt{p}}E_1^{n-2N}
		\lesssim \left(D^\beta+\frac{a}{\sqrt{p}}\right)\sum_{l=1}^k e_l^n-\frac{b}{\sqrt{p}}\sum_{l=1}^k e_l^{n-2N}.
	$
	Moreover, according to \eqref{eq4_59}, \eqref{eq4_64} and  Lemma \ref{lemma4.1}, the above inequality indicates $E_1^n\lesssim\sum_{l=1}^k e_l^n$.
    Thus, from $t_n\geq (k-1)\tau$, $(2k-1)\beta-(2-\beta)/r\geq0$, \eqref{eq4_65} and \eqref{eq4_66}, we get
	\begin{equation}\label{eq4_67}
		\begin{aligned}
			E_1^n &\lesssim \sum_{l=1}^{k-1} \rho_1^{\frac{2-\beta}{r}} \left(t_n^{[l]}\right)^{(2l-1)\beta-\frac{2-\beta}{r}}
			+\rho_1^{\frac{2-\beta}{r}}\left(t_n^{[k]}\right)^{(2k-1)\beta-\frac{2-\beta}{r}}
			\lesssim\rho_1^{\frac{2-\beta}{r}}.
		\end{aligned}
	\end{equation}
	For $E^n_2$, an analogous argument gives
	\begin{equation}\label{eq4_68}
		E_2^n \lesssim  \rho_1^{\frac{2-\beta}{r}}.
	\end{equation}
	Since the approximation property \eqref{eq4_2} of $R_h$ implies that
	\begin{equation*}
		\begin{aligned}
			\|R_3^n\| &\lesssim h^{q+1} \|\widetilde{u}^n\|_{H^{q+1}([0,L])} \lesssim h^{q+1},
			\|R_4^n\| \lesssim h^{q+1} \|\widetilde{u}^{n-2N}\|_{H^{q+1}([0,L])} \lesssim h^{q+1},\\
			\|R_5^n\| &\lesssim \int_0^{t_n}\left\|\partial_t\big(w^n-R_h w^n\big)\right\|dt
			\lesssim h^{q+1}\int_0^{t_n}\left\|\partial_t w^n\right\|_{H^{q+1}([0,L])}dt
			\lesssim h^{q+1},\\
			\|R_6^n\| &\lesssim \int_0^{t_n}\left\|\partial_t\big(\nabla\widetilde{u}^n-\nabla R_h \widetilde{u}^n\big)\right\|dt
			\lesssim h^{q}\int_0^{t_n}\left\|\partial_t \widetilde{u}^n\right\|_{H^{q+1}([0,L])}dt
			\lesssim h^{q},
		\end{aligned}
	\end{equation*}
	the combination of \eqref{eq4_58} and \eqref{eq4_59} and  Lemma \ref{lemma4.7} with $\gamma=-1$ can easily obtain
	\begin{equation}\label{eq4_69}
		E_i^n \lesssim h^{q+1} \rho_1 \left(t_n^{[k]}\right)^{\beta-1}
		\left(\frac{\rho_1}{t_n^{[k]}}\right)^{-1}
	=h^{q+1}\left(t_n^{[k]}\right)^{\beta},\quad i=3,4,5.
	\end{equation}
	and
	\begin{equation}\label{eq4_70}
		E_6^n \lesssim h^{q} \rho_1 \left(t_n^{[k]}\right)^{\beta-1}
		\left(\frac{\rho_1}{t_n^{[k]}}\right)^{-1}
		=h^{q}\left(t_n^{[k]}\right)^{\beta}.
	\end{equation}
	By applying Lemma \ref{lemma4.1} with \eqref{eq4_59} and $\|e_w^{n-2N}\|+p\|\nabla e_{\widetilde{u}}^{n-2N}\|=0,\ 1\leq n\leq 2N$, we obtain $\sqrt{\left\|e_w^n\right\|^2+p\left\|\nabla e_{\widetilde{u}}^n\right\|^2}
	\lesssim \sum_{i=1}^{6} E_i^n$.
	Since $w_h^n-w^n=e_w^n+R_h w^n-w^n$ and $\nabla\widetilde{u}_h^n-\nabla\widetilde{u}^n=\nabla e_{\widetilde{u}}^n+\nabla R_h \widetilde{u}^n-\nabla\widetilde{u}^n$, it leads to
	\begin{equation}\label{eq4_71}
		\begin{aligned}
			\left\|w_h^n-w^n\right\|+p\left\|\nabla\widetilde{u}_h^n-\nabla\widetilde{u}^n\right\|
			&\lesssim \sqrt{\left\|e_w^n\right\|^2+p\left\|\nabla e_{\widetilde{u}}^n\right\|^2}+h^{q+1}+h^{q}
			\lesssim \sum_{i=1}^6 E_i^n +h^{q}.
		\end{aligned}
	\end{equation}
	Thus, the combination of  \eqref{eq4_62}, \eqref{eq4_67}, \eqref{eq4_68}, \eqref{eq4_69} and \eqref{eq4_70} immediately yields \eqref{eq4_56}.
\end{proof}
% \begin{remark}
% Denote $E_h^n=\max_{2(k-1)N+1\le n\le 2kN}
% \left(\left\|w_h^n-w^n\right\|+p\left\|\nabla \widetilde{u}_h^n-\nabla \widetilde{u}^n\right\|\right).$
% When $k=1$, we have the following two cases：
% if $\dfrac{2-\beta}{r}<1$, then
% $
% E_h^n
% \lesssim\rho_1^{\min\left\{\beta,\frac{2-\beta}{r}\right\}}+\rho_1^{\min\left\{2\beta,\frac{2-\beta}{r}\right\}}+h^{q};
% $
% if $\dfrac{2-\beta}{r}=1$, then
% $
% E_h^n
% 	\lesssim\rho_1^{\beta}\left(1+\ln\frac{t_n}{\rho_1}\right)+\rho_1^{\min\left\{2\beta,\frac{2-\beta}{r}\right\}}+h^{q};
% $
% and if $\dfrac{2-\beta}{r}>1$, then
% $
% E_h^n
% 	\lesssim\rho_1^{\beta}+\rho_1^{\min\left\{2\beta,\frac{2-\beta}{r}\right\}}+h^{q}.
% $

% Thus, Theorem \ref{theorem4.9} implies
% \begin{equation*}
% 	\max_{2(k-1)N+1\le n\le 2kN}
% 	\left(\left\|w_h^n-w^n\right\|+p\left\|\nabla \widetilde{u}_h^n-\nabla \widetilde{u}^n\right\|\right)
% 	\lesssim
% 	\begin{cases}
% 		\rho_1^{\min\left\{\beta,\frac{2-\beta}{r}\right\}}+h^{q}, &k=1,\\
% 		\rho_1^{\frac{2-\beta}{r}}+h^{q}, &k\ge 2.
% 	\end{cases}
% \end{equation*}

% \end{remark}

\section{Numerical experiments}
In this section, some numerical tests are carried out to verify the theoretical results established in the above sections. For the problem \eqref{eq1_1}--\eqref{eq1_4}, we take $p=1$, $a=0$, $b=1$, $\tau=1$, $L=1$ and $K=3$. Choose the initial values $\phi(x,t)=\psi(x,t)=0$ and the source term
\begin{equation*}
	f(x,t)=
	\begin{cases}
		\pi^2 t^\alpha+\Gamma(\alpha+1), & 0<t\le 1,\\[4pt]
		\pi^2 t^\alpha+\Gamma(\alpha+1)+\pi^2\left(t-1\right)^{2\alpha}-\left(t-1\right)^\alpha+\dfrac{\Gamma(2\alpha+1)}{\Gamma(\alpha+1)}\left(t-1\right)^\alpha, & 1<t\le 2,\\
		\pi^2 t^\alpha+\Gamma(\alpha+1)+\pi^2\left(t-1\right)^{2\alpha}-\left(t-1\right)^\alpha+\dfrac{\Gamma(2\alpha+1)}{\Gamma(\alpha+1)}\left(t-1\right)^\alpha\\
		\quad+\dfrac{\Gamma(3\alpha+1)}{\Gamma(2\alpha+1)}\left(t-2\right)^{2\alpha}+\pi^2\left(t-2\right)^{3\alpha}-\left(t-2\right)^{2\alpha}, & 2<t\le 3.
	\end{cases}
\end{equation*}
It is easy to know that the exact solution of problem \eqref{eq1_1}--\eqref{eq1_4} is
\begin{equation*}
	u(x,t)=
	\begin{cases}
		t^\alpha \sin\left(\pi x\right), & 0<t\le 1,\\[4pt]
		\big(t^\alpha+\left(t-1\right)^{2\alpha}\big)\sin\left(\pi x\right), & 1<t\le 2,\\[4pt]
		\big(\left(t-2\right)^{3\alpha}+\left(t-1\right)^{2\alpha}+t^\alpha\big)\sin\left(\pi x\right), & 2<t\le 3.
	\end{cases}
\end{equation*}
According to the symmetric fractional-order reduction technique, it is clear that $\widetilde{u}(x,t)=u(x,t)$ and $g(x,t)=f(x,t)$ for the above example. Define
$
	E(h,N,i)=\max_{2(i-1)N+1\le n\le 2iN} \|\nabla \widetilde{u}^n-\nabla \widetilde{u}_h^n\|_0$,
	$\mathrm{rate}_t=\log_2(\dfrac{E(h,N,i)}{E(h,2N,i)})$ and
   $ \mathrm{rate}_s=\log_2(\dfrac{E(h,N,i)}{E(\frac{h}{2},N,i)}).$
We take $q=3$ in the finite element space $X_h^q$ to test the temporal error, choosing a sufficiently small $h=1/100$ so that the spatial error can be neglected. In Tables \ref{Table1}--\ref{Table4}, the temporal convergence orders at $t=1,2,3$ are presented.
By taking different $\alpha$ and $r$, the results in Tables \ref{Table1} and \ref{Table2} show that the temporal convergence order is near $r$ at $t=1$ and near $2-\alpha/2$ at $t=2,3$. The results in Tables \ref{Table3} and \ref{Table4}, which are obtained by choosing $\alpha=1.5$, $r=5/3$ and $\alpha=1.8$, $r=11/9$, show that the temporal convergence order is near $2-\alpha/2$ at $t=1,2,3$.
In Tables \ref{Table6}--\ref{Table10}, we test the local maximum temporal convergence orders on the time intervals $((i-1)\tau,i\tau],\ i=1,2,3$.
By taking different $\alpha$ and $r$,
these obtained results demonstrate that the temporal convergence orders are near $\min\{r\alpha/2,2-\alpha/2\}$ on $\left(0,\tau\right]$, and near $2-\alpha/2$ on the $\left(\tau,2\tau\right]$ and $\left(2\tau,3\tau\right]$. Moreover, when $r\ge (4-\alpha)/\alpha$, the optimal convergence is achieved.
In order to examine the spatial convergence rate, we choose $N=100$ and linear element (i.e. $q=1$) in Table \ref{Table11}, and $N=1000$ and $q=2$ (i.e. quadratic elements) in Table \ref{Table12}.
By choosing different values of $h$ to test the spatial convergence orders on $\big((i-1)\tau,i\tau\big],\ i=1,2,3$ respectively, the results indicate that the spatial convergence orders are near 1 and 2, respectively. These results are consistent with what is predicted by Theorem \ref{theorem4.9}.
\begin{table}[htbp]
	\centering
	\caption{Temporal accuracy and errors at $t_{iN}$ with $\alpha=1.5$ and $r=1$
		\label{Table1}}
	\begin{tabular}{c|cc|cc|cc}
		\toprule
		\multirow{2}{*}{$N$} & \multicolumn{2}{c}{$i=1$} & \multicolumn{2}{c}{$i=2$} & \multicolumn{2}{c}{$i=3$} \\
		\cmidrule(lr){2-3}\cmidrule(lr){4-5}\cmidrule(lr){6-7}
		& $\|\nabla\widetilde u^{2N}-\nabla\widetilde u_h^{2N}\|_0$ & $rate_t$ & $\|\nabla\widetilde u^{4N}-\nabla\widetilde u_h^{4N}\|_0$ & $rate_t$ & $\|\nabla\widetilde u^{6N}-\nabla\widetilde u_h^{6N}\|_0$ & $rate_t$ \\
		\midrule
		400  & 6.5353e-05 &        & 2.6529e-04 &        & 1.2932e-03 &        \\
		800  & 3.7011e-05 & 0.8203 & 1.1223e-04 & 1.2411 & 5.4367e-04 & 1.2501 \\
		1600 & 2.0313e-05 & 0.8655 & 4.7529e-05 & 1.2396 & 2.2841e-04 & 1.2511 \\
		3200 & 1.0914e-05 & 0.8962 & 2.0164e-05 & 1.2370 & 9.5940e-05 & 1.2514 \\
		%\midrule
		%\multicolumn{1}{c}{--} & -- & 1 & -- & 1.25& -- &1.25 \\
		\bottomrule
	\end{tabular}
\end{table}
\begin{table}[htbp]
	\centering
	\caption{Temporal accuracy and errors at $t_{iN}$ with $\alpha=1.8$ and $r=2-\dfrac{\alpha}{2}=1.1$\label{Table2}}
	\begin{tabular}{c|cc|cc|cc}
		\toprule
		\multirow{2}{*}{$N$} & \multicolumn{2}{c}{$i=1$} & \multicolumn{2}{c}{$i=2$} & \multicolumn{2}{c}{$i=3$} \\
		\cmidrule(lr){2-3}\cmidrule(lr){4-5}\cmidrule(lr){6-7}
		& $\|\nabla\widetilde u^{2N}-\nabla\widetilde u_h^{2N}\|_0$ & $rate_t$ & $\|\nabla\widetilde u^{4N}-\nabla\widetilde u_h^{4N}\|_0$ & $rate_t$ & $\|\nabla\widetilde u^{6N}-\nabla\widetilde u_h^{6N}\|_0$ & $rate_t$ \\
		\midrule
		400  & 1.3900e-04 &        & 2.1190e-03 &        & 7.7940e-03 &        \\
		800  & 6.0392e-05 & 1.2027 & 9.8433e-04 & 1.1062 & 3.6396e-03 & 1.0986 \\
		1600 & 2.6030e-05 & 1.2142 & 4.5704e-04 & 1.1068 & 1.6994e-03 & 1.0987 \\
		3200 & 1.1131e-05 & 1.2257 & 2.1217e-04 & 1.1071 & 7.9348e-04 & 1.0988 \\
		%\midrule
		%\multicolumn{1}{c}{\text{Theoretical order}} & -- & 1.1 & -- & 1.1 & -- & 1.1 \\
		\bottomrule
	\end{tabular}
\end{table}
\begin{table}[htbp]
	\centering
	\caption{Temporal accuracy and errors at $t_{iN}$ with $\alpha=1.5$ and $r=\dfrac{2-\alpha/2}{\alpha/2}=\dfrac{5}{3}$\label{Table3}}
	\begin{tabular}{c|cc|cc|cc}
		\toprule
		\multirow{2}{*}{$N$} & \multicolumn{2}{c}{$i=1$} & \multicolumn{2}{c}{$i=2$} & \multicolumn{2}{c}{$i=3$} \\
		\cmidrule(lr){2-3}\cmidrule(lr){4-5}\cmidrule(lr){6-7}
		& $\|\nabla\widetilde u^{2N}-\nabla\widetilde u_h^{2N}\|_0$ & $rate_t$ & $\|\nabla\widetilde u^{4N}-\nabla\widetilde u_h^{4N}\|_0$ & $rate_t$ & $\|\nabla\widetilde u^{6N}-\nabla\widetilde u_h^{6N}\|_0$ & $rate_t$ \\
		\midrule
		400  & 5.3471e-05 &        & 1.1261e-04 &        & 7.1479e-04 &        \\
		800  & 2.2377e-05 & 1.2568 & 4.6754e-05 & 1.2682 & 2.9875e-04 & 1.2586 \\
		1600 & 9.3737e-06 & 1.2553 & 1.9526e-05 & 1.2597 & 1.2519e-04 & 1.2549 \\
		3200 & 3.9338e-06 & 1.2527 & 8.2085e-06 & 1.2502 & 5.2595e-05 & 1.2511 \\
		%\midrule
		%\multicolumn{1}{c}{\text{Theoretical order} } & -- & 1.25 & -- & 1.25 & -- & 1.25 \\
		\bottomrule
	\end{tabular}
\end{table}
\begin{table}[htbp]
	\centering
	\caption{Temporal accuracy and errors at $t_{iN}$ with $\alpha=1.8$ and $r=\dfrac{2-\alpha/2}{\alpha/2}=\dfrac{11}{9}$\label{Table4}}
	\begin{tabular}{c|cc|cc|cc}
		\toprule
		\multirow{2}{*}{$N$} & \multicolumn{2}{c}{$i=1$} & \multicolumn{2}{c}{$i=2$} & \multicolumn{2}{c}{$i=3$} \\
		\cmidrule(lr){2-3}\cmidrule(lr){4-5}\cmidrule(lr){6-7}
		& $\|\nabla\widetilde u^{2N}-\nabla\widetilde u_h^{2N}\|_0$ & $rate_t$ & $\|\nabla\widetilde u^{4N}-\nabla\widetilde u_h^{4N}\|_0$ & $rate_t$ & $\|\nabla\widetilde u^{6N}-\nabla\widetilde u_h^{6N}\|_0$ & $rate_t$ \\
		\midrule
		400  & 1.8010e-04 &        & 2.1241e-03 &        & 7.5427e-03 &        \\
		800  & 8.2072e-05 & 1.1338 & 9.8924e-04 & 1.1024 & 3.5201e-03 & 1.0995 \\
		1600 & 3.7390e-05 & 1.1342 & 4.6061e-04 & 1.1028 & 1.6427e-03 & 1.0995 \\
		3200 & 1.7045e-05 & 1.1333 & 2.1447e-04 & 1.1028 & 7.6659e-04 & 1.0996 \\
		%\midrule
		%\multicolumn{1}{c}{\text{Theoretical order}} & -- & 1.1 & -- & 1.1 & -- & 1.1 \\
		\bottomrule
	\end{tabular}
\end{table}
%
%\begin{table}[htbp]
	%\centering
	%\caption{Temporal accuracy and errors with $\alpha=1.5$ and $r=1$\label{Table5}}
	%\begin{tabular}{c|cc|cc|cc}
		%\toprule
		%\multirow{2}{*}{$N$} & \multicolumn{2}{c}{$i=1$} & \multicolumn{2}{c}{$i=2$} & \multicolumn{2}{c}{$i=3$} \\
		%\cmidrule(lr){2-3}\cmidrule(lr){4-5}\cmidrule(lr){6-7}
		%& $E(h,N,i)$ & $rate_t$ & $E(h,N,i)$ & $rate_t$ & $E(h,N,i)$ & $rate_t$ \\
		%\midrule
		%400  & 4.0935e-04 &        & 2.6529e-04 &        & 1.2932e-03 &        \\
		%800  & 2.2072e-04 & 0.8911 & 1.1223e-04 & 1.2411 & 5.4367e-04 & 1.2501 \\
		%1600 & 1.1729e-04 & 0.9121 & 4.7529e-05 & 1.2396 & 2.2841e-04 & 1.2511 \\
		%3200 & 6.1620e-05 & 0.9286 & 2.0164e-05 & 1.2370 & 9.5940e-05 & 1.2514 \\
		%\midrule
		%\multicolumn{1}{c}{\text{Theoretical order}} & -- & 0.75 & -- & 1.25 & -- & 1.25 \\
		%\bottomrule
	%\end{tabular}
%\end{table}
%
\begin{table}[htbp]
	\centering
	\caption{Temporal accuracy and errors with $\alpha=1.5$ and $r=\dfrac{2-\alpha/2}{\alpha/2}=\dfrac{5}{3}$\label{Table6}}
	\begin{tabular}{c|cc|cc|cc}
		\toprule
		\multirow{2}{*}{$N$} & \multicolumn{2}{c}{$i=1$} & \multicolumn{2}{c}{$i=2$} & \multicolumn{2}{c}{$i=3$} \\
		\cmidrule(lr){2-3}\cmidrule(lr){4-5}\cmidrule(lr){6-7}
		& $E(h,N,i)$ & $rate_t$ & $E(h,N,i)$ & $rate_t$ & $E(h,N,i)$ & $rate_t$ \\
		\midrule
		400  & 5.3471e-05 &        & 4.5157e-04 &        & 1.4127e-03 &        \\
		800  & 2.2377e-05 & 1.2568 & 1.9002e-04 & 1.2488 & 5.9437e-04 & 1.2490 \\
		1600 & 9.3737e-06 & 1.2553 & 7.9924e-05 & 1.2495 & 2.4999e-04 & 1.2495 \\
		3200 & 3.9338e-06 & 1.2527 & 3.3612e-05 & 1.2497 & 1.0514e-04 & 1.2495 \\
		%midrule
		%\multicolumn{1}{c}{\text{Theoretical order}} & -- & 1.25 & -- & 1.25 & -- & 1.25 \\
		\bottomrule
	\end{tabular}
\end{table}
\begin{table}[htbp]
	\centering
	\caption{Temporal accuracy and errors with $\alpha=1.5$ and $r=3$\label{Table7}}
	\begin{tabular}{c|cc|cc|cc}
		\toprule
		\multirow{2}{*}{$N$} & \multicolumn{2}{c}{$i=1$} & \multicolumn{2}{c}{$i=2$} & \multicolumn{2}{c}{$i=3$} \\
		\cmidrule(lr){2-3}\cmidrule(lr){4-5}\cmidrule(lr){6-7}
		& $E(h,N,i)$ & $rate_t$ & $E(h,N,i)$ & $rate_t$ & $E(h,N,i)$ & $rate_t$ \\
		\midrule
		400  & 1.3508e-04 &        & 9.0919e-04 &        & 2.7652e-03 &        \\
		800  & 5.6984e-05 & 1.2452 & 3.8322e-04 & 1.2464 & 1.1647e-03 & 1.2474 \\
		1600 & 2.4003e-05 & 1.2474 & 1.6133e-04 & 1.2482 & 4.9016e-04 & 1.2487 \\
		3200 & 1.0104e-05 & 1.2483 & 6.7877e-05 & 1.2490 & 2.0619e-04 & 1.2493 \\
		%\midrule
		%\multicolumn{1}{c}{\text{Theoretical order}} & -- & 1.25 & -- & 1.25 & -- & 1.25 \\
		\bottomrule
	\end{tabular}
\end{table}
\begin{table}[htbp]
	\centering
	\caption{Temporal accuracy and errors  with $\alpha=1.7$ and $r=1$\label{Table8}}
	\begin{tabular}{c|cc|cc|cc}
		\toprule
		\multirow{2}{*}{$N$} & \multicolumn{2}{c}{$i=1$} & \multicolumn{2}{c}{$i=2$} & \multicolumn{2}{c}{$i=3$} \\
		\cmidrule(lr){2-3}\cmidrule(lr){4-5}\cmidrule(lr){6-7}
		& $E(h,N,i)$ & $rate_t$ & $E(h,N,i)$ & $rate_t$ & $E(h,N,i)$ & $rate_t$ \\
		\midrule
		400  & 2.3641e-04 &        & 1.0655e-03 &        & 4.5924e-03 &        \\
		800  & 1.3651e-04 & 0.7923 & 4.7341e-04 & 1.1704 & 2.0721e-03 & 1.1481 \\
		1600 & 7.7113e-05 & 0.8240 & 2.0991e-04 & 1.1733 & 9.3467e-04 & 1.1486 \\
		3200 & 4.2722e-05 & 0.8520 & 9.2870e-05 & 1.1765 & 4.2154e-04 & 1.1488 \\
		%\midrule
		%\multicolumn{1}{c}{\text{Theoretical order}} & -- & 0.85 & -- &1.15& -- & 1.15 \\
		\bottomrule
	\end{tabular}
\end{table}
\begin{table}[htbp]
	\centering
	\caption{Temporal accuracy and errors  with $\alpha=1.7$ and $r=\dfrac{2-\alpha/2}{\alpha/2}=\dfrac{23}{17}$\label{Table9}}
	\begin{tabular}{c|cc|cc|cc}
		\toprule
		\multirow{2}{*}{$N$} & \multicolumn{2}{c}{$i=1$} & \multicolumn{2}{c}{$i=2$} & \multicolumn{2}{c}{$i=3$} \\
		\cmidrule(lr){2-3}\cmidrule(lr){4-5}\cmidrule(lr){6-7}
		& $E(h,N,i)$ & $rate_t$ & $E(h,N,i)$ & $rate_t$ & $E(h,N,i)$ & $rate_t$ \\
		\midrule
		400  & 1.2924e-04 &        & 1.2468e-03 &        & 4.2394e-03 &        \\
		800  & 5.6947e-05 & 1.1824 & 5.6139e-04 & 1.1511 & 1.9112e-03 & 1.1494 \\
		1600 & 2.5129e-05 & 1.1803 & 2.5272e-04 & 1.1515 & 8.6142e-04 & 1.1497 \\
		3200 & 1.1110e-05 & 1.1774 & 1.1376e-04 & 1.1515 & 3.8824e-04 & 1.1498 \\
		%\midrule
		%\multicolumn{1}{c}{\text{Theoretical order}} & -- & 1.15 & -- & 1.15 & -- & 1.15 \\
		\bottomrule
	\end{tabular}
\end{table}
\begin{table}[htbp]
	\centering
	\caption{Temporal accuracy and errors  with $\alpha=1.7$ and $r=3$\label{Table10}}
	\begin{tabular}{c|cc|cc|cc}
		\toprule
		\multirow{2}{*}{$N$} & \multicolumn{2}{c}{$i=1$} & \multicolumn{2}{c}{$i=2$} & \multicolumn{2}{c}{$i=3$} \\
		\cmidrule(lr){2-3}\cmidrule(lr){4-5}\cmidrule(lr){6-7}
		& $E(h,N,i)$ & $rate_t$ & $E(h,N,i)$ & $rate_t$ & $E(h,N,i)$ & $rate_t$ \\
		\midrule
		400  & 3.9265e-04 &        & 2.6887e-03 &        & 8.6474e-03 &        \\
		800  & 1.7788e-04 & 1.1423 & 1.2144e-03 & 1.1467 & 3.9015e-03 & 1.1482 \\
		1600 & 8.0377e-05 & 1.1461 & 5.4779e-04 & 1.1485 & 1.7591e-03 & 1.1492 \\
		3200 & 3.6271e-05 & 1.1479 & 2.4697e-04 & 1.1493 & 7.9288e-04 & 1.1496 \\
		%\midrule
		%\multicolumn{1}{c}{\text{Theoretical order}} & -- & 1.15 & -- & 1.15 & -- & 1.15 \\
		\bottomrule
	\end{tabular}
\end{table}
\begin{table}[htbp]
	\centering
	\caption{Spatial accuracy and errors with $\alpha=1.5$,$r=\dfrac{2-\alpha/2} {\alpha/2}=\dfrac{5}{3}$ and $q=1$\label{Table11}}
	\begin{tabular}{c|cc|cc|cc}
		\toprule
		\multirow{2}{*}{$h$} & \multicolumn{2}{c}{$i=1$} & \multicolumn{2}{c}{$i=2$} & \multicolumn{2}{c}{$i=3$} \\
		\cmidrule(lr){2-3}\cmidrule(lr){4-5}\cmidrule(lr){6-7}
		& $E(h,N,i)$ & $rate_s$ & $E(h,N,i)$ & $rate_s$ & $E(h,N_h,i)$ & $rate_s$ \\
		\midrule
		1/8   & 2.5121e-01 &        & 9.6173e-01 &        & 3.5662e+00 &        \\
		1/16  & 1.2584e-01 & 0.9974 & 4.8176e-01 & 0.9973 & 1.7864e+00 & 0.9973 \\
		1/32  & 6.2947e-02 & 0.9993 & 2.4099e-01 & 0.9993 & 8.9363e-01 & 0.9993 \\
		1/64  & 3.1478e-02 & 0.9998 & 1.2051e-01 & 0.9998 & 4.4688e-01 & 0.9998 \\
		%\midrule
		%\multicolumn{1}{c}{\text{Theoretical order}} & -- & 1 & -- & 1 & -- & 1 \\
		\bottomrule
	\end{tabular}
\end{table}
\begin{table}[htbp]
	\centering
	\caption{Spatial accuracy and errors with $\alpha=1.5$, $r=\dfrac{2-\alpha/2}{\alpha/2}=\dfrac{5}{3}$ and $q=2$\label{Table12}}
	\begin{tabular}{c|cc|cc|cc}
		\toprule
		\multirow{2}{*}{$h$} & \multicolumn{2}{c}{$i=1$} & \multicolumn{2}{c}{$i=2$} & \multicolumn{2}{c}{$i=3$} \\
		\cmidrule(lr){2-3}\cmidrule(lr){4-5}\cmidrule(lr){6-7}
		& $E(h,N,i)$ & $rate_s$ & $E(h,N,i)$ & $rate_s$ & $E(h,N,i)$ & $rate_s$ \\
		\midrule
		1/8   & 1.2739e-02 &        & 4.8770e-02 &        & 1.8084e-01 &        \\
		1/16  & 3.1900e-03 & 1.9976 & 1.2213e-02 & 1.9976 & 4.5286e-02 & 1.9976 \\
		1/32  & 7.9800e-04 & 1.9991 & 3.0546e-03 & 1.9993 & 1.1328e-02 & 1.9991 \\
		1/64  & 2.0019e-04 & 1.9950 & 7.6450e-04 & 1.9984 & 2.8408e-03 & 1.9956 \\
		%\midrule
		%\multicolumn{1}{c}{\text{Theoretical order}} & -- & 2 & -- & 2 & -- & 2 \\
		\bottomrule
	\end{tabular}
\end{table}
\section{Conclusions}
A local error estimate has been established for the fully discrete finite element scheme of
time fractional diffusion-wave equations with a constant delay and weakly singular solutions.
By introducing the auxiliary variable $w={}^{C}_{0}D_{t}^{\alpha/2}(u-t\,\psi(x,0))$, we recast the
problem as an equivalent coupled system and the regularity of $w$ is discussed in detail.
On a symmetric graded time mesh, we constructed a fully discrete finite element scheme for the coupled
problem and established the stability of the resulting scheme.
Since the barrier function framework of \cite{Kopteva2025Pointwise} cannot be applied directly once a delay term is present, we construct new barrier functions and develop the associated theory, thereby obtaining a stability property for the discrete fractional derivative operator.
Combining the analysis of the local truncation error in time yields a local error estimate for the proposed scheme. Numerical experiments verify the validity of the theoretical results. The method developed can be further extended to other fractional
evolution equations with weakly singular solutions and time delays, while nonlinear delay terms and
the variable-order case are left for future work.

\appendix
\section{Proofs of Lemmas 4.1--4.6}
\label{Proofs of Lemmas 4.1--4.6}
%In this section,
We now provide detailed proofs for Lemmas 4.1--4.6.
\subsection{Proof of Lemma \ref{lemma4.1}}
\begin{proof}
	From \eqref{eq3_3} and \eqref{eq4_10}, it is clear that
	% \begin{equation*}
	% 	a_0^{(n)}V^n-\sum_{i=1}^{n-1}\left(a_{n-i-1}^{(n)}-a_{n-i}^{(n)}\right)V^i
	% 	-a_{n-1}^{(n)}V^0-\lambda_0 V^n-\lambda_1 V^{n-2N}\leq 0,
	% \end{equation*}
	% i.e.,
	\begin{equation}
		\left(a_0^{(n)}-\lambda_0\right)V^n\leq \lambda_1 V^{n-2N}
		+\sum_{i=1}^{n-1}\left(a_{n-i-1}^{(n)}-a_{n-i}^{(n)}\right)V^i
		+a_{n-1}^{(n)}V^0,
		\label{eq5_1}
	\end{equation}
    where \cite{Tan2022L1} and $\lambda_0 \bar{\rho}^\beta<\left\{\Gamma(2-\beta)\right\}^{-1}$ imply that the coefficients $a^{(n)}_i$ satisfies
	\begin{equation}	\label{eq5_2}
		a_{n-i-1}^{(n)}-a_{n-i}^{(n)}>0,\quad a_{n-1}^{(n)}>0,\quad	a_0^{(n)}-\lambda_0=\frac{\rho_n^{-\beta}}{\Gamma\left(2-\beta\right)}-\lambda_0>0.
	\end{equation}
	For the case $1\leq n\leq 2N$, i.e., $k=1$, since $\lambda_1\geq0,$ $V^0\leq0$ and $V^{n-2N}\leq0$, by mathematical induction, it is easy to know that \eqref{eq5_1} indicates $V^n\leq 0$. For the case $n>2N,$ based on the above result, we assume that $V^n\leq 0$ holds for $2(m-1)N+1\leq n\leq 2mN$. Then, when $2mN+1\leq n\leq 2(m+1))N,$ we have $V^{i},i\leq 2mN.$ Based on these conditions, applying mathematical induction to \eqref{eq5_1} yields
    $V^n\leq0, 2mN+1\leq n\leq 2(m+1))N.$
    These above discussions demonstrate, for each $k\geq0,$
	$V^n\leq 0, 2(k-1)N+1\leq n\leq 2kN.$
\end{proof}

\subsection{Proof of Lemma \ref{lemma4.2}}
\begin{proof}
Let $t_j'-(k-1)\tau=t_j'^{[k]}$.
From \eqref{eq4_13}, for each $j$, there exists a $n$ such that $t_j'=t_n$. Then we have $\left(\rho_1/t_j'^{[k]}\right)^{1/r}=\left(\rho_1/t_n^{[k]}\right)^{1/r}$,
and using this and the properties of the symmetric graded mesh \eqref{eq3_2} yields
		\begin{equation}\label{eq5_3}
			%\begin{aligned}
				\rho_n/t_n^{[k]}
				\lesssim \left(\rho_1/t_n^{[k]}\right)^{1/r}
				\simeq\left[\frac{N^{-r}}{(C\left(j-2(k-1)N\right))^r/N^r}\right]^{1/r}
				=C^{-1}\left[j-2(k-1)N\right]^{-1}.
			%\end{aligned}
		\end{equation}
		Since $t_{n-1}^{[k]}<(\frac{C(j-2(k-1)N)}{N})^r$ and $C$ is a sufficiently large constant,
		\eqref{eq5_3} gives
		$t_j'^{[k]}=t_n^{[k]}\ \approx t_{n-1}^{[k]}
		<(\frac{C(j-2(k-1)N)}{N})^r$.
		According to \eqref{eq4_13}, this implies that $t_j'^{[k]}$ is sufficiently close to
		$(\frac{C(j-2(k-1)N)}{N})^r$. Hence, \eqref{eq4_15} holds.

For $t\leq t_m'$, a direct comparison gives $(t-(k-1)\tau)/(t_m'^{[k]})\leq 1$ and
		$(t_m'^{[k]})^{-\theta}\leq(t-(k-1)\tau)^{-\theta}$,
		so
        $
			B(t,t_m')
			=(t-(k-1)\tau)(t_m'^{[k]})^{-\theta-1}.
		$
		This leads to
		\begin{equation*}
			D_t^\beta B(t,t_m')
			=
			\frac{1}{\Gamma\left(1-\beta\right)}\left(t_m'^{[k]}\right)^{-\theta-1}\int_{\left(k-1\right)\tau}^{t}\left(t-s\right)^{-\beta}\,ds
			=
			\frac{1}{\Gamma\left(\theta\right)\,\theta}\left(t_m'^{[k]}\right)^{-\theta-1}\left(t-(k-1)\tau\right)^\theta,
		\end{equation*}
		which confirms \eqref{eq4_17}.
		For $t>t_m'$, the inequalities $(t-(k-1)\tau)/(t_m'^{[k]})\geq 1$,
		$(t_m'^{[k]})^{-\theta}>(t-(k-1)\tau)^{-\theta}$
		imply $B(t,t_m')=(t-(k-1)\tau)^{-\theta}$. Consequently,
		\begin{equation}\label{eq5_4}
			\begin{aligned}
				D_t^\beta B\left(t,t_m'\right)
				&=\frac{1}{\Gamma\left(1-\beta\right)}\int_{\left(k-1\right)\tau}^{t}\left(t-s\right)^{-\beta}
				\partial_s B\left(s,t_m'\right)\,ds\\
				&=\frac{1}{\Gamma\left(1-\beta\right)}\int_{\left(k-1\right)\tau}^{t_m'}\left(t-s\right)^{-\beta}
				\left(t_m'^{[k]}\right)^{-\theta-1}\,ds
				-\frac{1-\beta}{\Gamma\left(1-\beta\right)}
				\int_{t_m'}^{t}\left(t-s\right)^{-\beta}\left(s-\left(k-1\right)\tau\right)^{-\theta-1}\,ds
				:=I_1+I_2.
			\end{aligned}
		\end{equation}
		Next, we estimate the bounds of $I_1$ and $I_2$, respectively.
        Noting that $(t-s)^{-\beta}\geq(t-(k-1)\tau)^{-\beta}$, then
		\begin{equation}\label{eq5_5}
				I_1\geq \frac{1}{\Gamma\left(1-\beta\right)}\int_{\left(k-1\right)\tau}^{t_m'}
				\left(t-(k-1)\tau\right)^{-\beta}\left(t_m'^{[k]}\right)^{-\theta-1}\,ds
				=\frac{\left(t-(k-1)\tau\right)^{-\beta}}{\Gamma\left(1-\beta\right)}
				\left(t_m'^{[k]}\right)^{-\theta}.
		\end{equation}
		By applying the substitutions $u=s-(k-1)\tau$, $\widehat{u}=u/(t-(k-1)\tau)$ and $\beta+\theta=1$ to $I_2$, we have
		\begin{equation}\label{eq5_6}
			\begin{aligned}
				I_2=\frac{-\left(1-\beta\right)}{\Gamma\left(1-\beta\right)}
				\int_{t_m'^{[k]}}^{t-(k-1)\tau}
				\left(t-(k-1)\tau-u\right)^{-\beta}u^{-\theta-1}\,du
				=\frac{-\left(1-\beta\right)}{\Gamma\left(1-\beta\right)}
				\int_{\frac{t_m'^{[k]}}{t-(k-1)\tau}}^{1}
				\left(t-(k-1)\tau\right)^{-1}(1-\widehat{u})^{-\beta}
				\widehat{u}^{-\theta-1}\,d\widehat{u}.
			\end{aligned}
		\end{equation}
		Setting $x=\dfrac{1}{\widehat{u}}$, \eqref{eq5_6} gives
		\begin{equation}\label{eq5_7}
			\begin{aligned}
				I_2&=\frac{-\left(1-\beta\right)}{\Gamma\left(1-\beta\right)}\left(t-\left(k-1\right)\tau\right)^{-1}
				\int_{1}^{\frac{t-\left(k-1\right)\tau}{t_m'^{[k]}}}\left(x-1\right)^{-\beta}\,dx=-\frac{\left(t-\left(k-1\right)\tau\right)^{-1}}{\Gamma\left(1-\beta\right)}
				\left(\frac{t-\left(k-1\right)\tau}{t_m'^{[k]}}-1\right)^{1-\beta}\\
				&=\frac{-\left(t-\left(k-1\right)\tau\right)^{-1}}{\Gamma\left(1-\beta\right)}
				\left(\frac{t_m'^{[k]}}{t-\left(k-1\right)\tau}\right)^{-\left(1-\beta\right)}
				\left(1-\frac{t_m'^{[k]}}{t-\left(k-1\right)\tau}\right)^{1-\beta}\\
				&\geq -\frac{\left(t-(k-1)\tau\right)^{-1}}{\Gamma\left(1-\beta\right)}
				\left(\frac{t_m'^{[k]}}{t-(k-1)\tau}\right)^{\beta-1}
				\left[1-(1-\beta)\frac{t_m'^{[k]}}{t-(k-1)\tau}\right],
			\end{aligned}
		\end{equation}
		where the Bernoulli's inequality $(1-y)^{1-\beta}\leq 1-(1-\beta) y$ is used.
		In view of \eqref{eq5_5} and \eqref{eq5_7}, \eqref{eq5_4} indicates
		\begin{equation}\label{eq5_8}
			\begin{aligned}
				D_t^\beta B(t,t_m')\geq &
				\frac{\left(t-(k-1)\tau\right)^{-\beta}
				\left(t_m'^{[k]}\right)^{-\theta}}{\Gamma\left(1-\beta\right)}-\frac{\left(t-(k-1)\tau\right)^{-\beta}
				\left(t_m'^{[k]}\right)^{\beta-1}}{\Gamma\left(1-\beta\right)}
				+\frac{(1-\beta)}{\Gamma\left(1-\beta\right)}\left(t-(k-1)\tau\right)^{-\beta-1}
				\left(t_m'^{[k]}\right)^\beta\\
				&=\frac{\theta}{\Gamma(\theta)}
				\left(\frac{t_m'^{[k]}}{t-(k-1)\tau}\right)^\beta
				\left(t-(k-1)\tau\right)^{-1}.
			\end{aligned}
		\end{equation}
		So, \eqref{eq4_16} is true.
		
Since $p_l<j-2(k-1)N\leq p_{l+1}$,
		by \eqref{eq4_15} and the expression of $p_l$, we have
		\begin{equation}\label{eq5_9}
			\frac{t_{p_l}'}{t_j'^{[k]}}\geq
			\frac{t_{p_l}'}{t_{p_{l+1}}'}\simeq
			\frac{\left(C2^l(p-2(k-1)N)\right)^r/N^r}
			{\left(C2^{l+1}(p-2(k-1)N)\right)^r/N^r}
			=2^{-r}.
		\end{equation}
		A direct consequence of \eqref{eq5_9} is that $t_{p_{l+1}}'\simeq 2^r t_{p_l}'$.
		It then follows that $t_{p_l}'<t_j'^{[k]}\leq t_{p_{l+1}}'\simeq 2^r t_{p_l}'$.
		Therefore $c_l=\left(t_p'^{[k]}/t_{p_l}'\right)^\gamma
		\simeq\left(t_p'^{[k]}/t_j'^{[k]}\right)^\gamma,$
		which yields \eqref{eq4_18}. Besides,
		from \eqref{eq4_15}, we can deduce that $c_l=2^{-rl\gamma}.$
		For $L>0$, if $\gamma>0$, then
		$\sum_{l=0}^{L}2^{-r\gamma l}\lesssim 1$;
		if $\gamma=0$, then $\displaystyle\sum_{l=0}^{L}c_l=L+1$
		and $\ln\left(t_{p_L}'/t_p'^{[k]}\right)
		=\ln 2^{rL}\simeq L$;
		and if $\gamma<0$, then it follows that $\displaystyle\sum_{l=0}^{L}c_l\lesssim c_L$. Thus, the above discussion confirms
		\eqref{eq4_19}.
\end{proof}

\subsection{Proof of Lemma \ref{lemma4.3}}
\begin{proof}
	When $t_j'\leq t_{p}'\leq t_{p_l}'+(k-1)\tau$, since
	$B(t_j',t_{p_l}'+(k-1)\tau)
	= (t_j'^{[k]})(t_{p_l}')^{-\theta-1}$
	depends linearly on $t_j'$, it suggests
	$
		(D_t^\beta-D^\beta)B(t,t_{p_l}'+(k-1)\tau)|_{t=t_j'} = 0.
	$
	Therefore, it follows from  $C_L(t_j') = \sum_{l=0}^{L}c_l\,B(t_j',t_{p_l}'+(k-1)\tau)$ that  \eqref{eq4_20} is true.
	
	Next, we discuss $(D_t^\beta-D^\beta)C_L(t)|_{t=t_j'}$ under the case $t_j'>t_p'$.
	Note that, for $t_j'>t_{p_l}'+(k-1)\tau$, it holds
	\begin{equation}\label{eq5_10}
		\begin{aligned}
			\left.(D_t^\beta-D^\beta)B\left(t,t_{p_l}'+(k-1)\tau\right)\right|_{t=t_j'}
			=& \frac{1}{\Gamma\left(1-\beta\right)}\int_{(k-1)\tau}^{t_j'}
			\left(t_j'-s\right)^{-\beta}\partial_s\left(B-B^I\right)\left(s,t_{p_l}'+(k-1)\tau\right)\,ds\\
			=&\frac{\beta}{\Gamma\left(1-\beta\right)}\int_{(k-1)\tau}^{t_{p_l}'+(k-1)\tau}
			\left(t_j'-s\right)^{-\beta-1}\left(B-B^I\right)\left(s,t_{p_l}'+(k-1)\tau\right)\,ds\\
			&+\frac{\beta}{\Gamma\left(1-\beta\right)}\int_{t_{p_l}'+(k-1)\tau}^{t_j'}
			\left(t_j'-s\right)^{-\beta-1}\left(B-B^I\right)\left(s,t_{p_l}'+(k-1)\tau\right)\,ds,
		\end{aligned}
	\end{equation}
	where $B^I$ denotes the standard piecewise linear interpolant of $B$ on temporal submesh $\{t_j'\}$.
	When $s\leq t_{p_l}'+(k-1)\tau$, one has
	$\left(B^I-B\right)\left(s,t_{p_l}'+(k-1)\tau\right)=0$,
	so the second line of \eqref{eq5_10} vanishes.
	When $s>t_{p_l}'+(k-1)\tau$, since
	$B\left(s,t_{p_l}'+(k-1)\tau\right) = \left(s-(k-1)\tau\right)^{-\theta}=B\left(s,t_{p}'\right)$,
	we have
	\begin{equation}
		\label{eq5_11}
		\left(B^I-B\right)\left(s,t_{p_l}'+(k-1)\tau\right) = \left(B^I-B\right)\left(s,t_{p}'\right).
	\end{equation}
	By using \eqref{eq5_10}, \eqref{eq5_11} and the expression of the $C_L\left(t_j'\right)$ for $t_j'>t_p'$, it leads to
	\begin{equation}
		\label{eq5_12}
		\begin{aligned}
			\left.\left(D_t^\beta - D^\beta\right) C_L\left(t\right)\right|_{t=t_j'}
			&= \frac{\beta}{\Gamma\left(1-\beta\right)}
			\sum_{l:\,p_l+2(k-1)N<j}c_l
			\int_{t_{p_l}'+(k-1)\tau}^{t_j'}\left(t_j'-s\right)^{-\beta-1}
			\left(B^I-B\right)\left(s,t_{p}'\right)\,ds\\
			&\lesssim\sum_{l:\,p_l+2(k-1)N<j}c_l
			\int_{t_p'}^{t_j'}\left(t_j'-s\right)^{-\beta-1}\left(B^I-B\right)\left(s,t_{p}'\right)\,ds\\
			&= S_J\sum_{i=p+1}^{j}\int_{t_{i-1}'}^{t_i'}
			\left(t_j'-s\right)^{-\beta-1}\left(B^I-B\right)\left(s,t_{p}'\right)\,ds,
		\end{aligned}
	\end{equation}
    where $\sum_{l:p_l+2(k-1)N<j}c_l\lesssim S_J$ in view of \eqref{eq4_19} is used.
    Let $\rho_i'=t_i'-t_{i-1}'$.
	If $p+1\leq i<j-1$, then
	$|(B-B^I)(s,t_p')|\leq \rho_i'^2|B''(t,t_p')|_{t=t_{i-1}'}|$ for $s\in(t_{i-1}',t_i')$;
	and if $i=j$, then
	$|(B-B^I)(s,t_p')| \leq \rho_i'(t_j'-s)B''(t,t_p')|_{t=t_{i-1}'}|$ for $s\in(t_{i-1}',t_i')$.
	Consequently, we have
	$
		|(B-B^I)(s,t_p')| \leq \rho_i'^2\min\{1,\frac{t_j'-s}{\rho_j'}\}|B''(t,t_p')|_{t=t_{i-1}'}|
	$
     for $t_j'>t_p'$.
	Since $t_i'\simeq t_{i-1}'\simeq s$, we obtain
	$|B''(t,t_p')|_{t=t_{i-1}'}|\simeq |B''(s,t_p')|_{t=t_{i-1}'}|\simeq (s-(k-1)\tau)^{-\theta-2}$ and $\rho_i' \lesssim \rho_1^{1/r}(t_i'-(k-1)\tau)^{1-1/r}
	\simeq \rho_1^{1/r}\left(s-(k-1)\tau\right)^{1-1/r}$.
	based on the above estimates, \eqref{eq5_12} leads to
	\begin{equation}
		\label{eq5_13}
		\begin{aligned}
			\left|\left.\left(D_t^\beta - D^\beta\right) C_L\left(t\right)\right|_{t=t_j'} \right|
			&\lesssim S_J\sum_{i=p+1}^{j}\int_{t_{i-1}'}^{t_i'}
			\left(t_j'-s\right)^{-\beta-1}\rho_i'^2
			\min\left\{1,\frac{t_j'-s}{\rho_j'}\right\}
			\left(s-(k-1)\tau\right)^{-\theta-2}\,ds\\
			&= S_J\rho_1^{2/r}\int_{t_p'}^{t_j'}
			\left(t_j'-s\right)^{-\beta-1}\left(s-(k-1)\tau\right)^{-\frac{2}{r}-\theta}
			\min\left\{1,\frac{t_j'-s}{\rho_j'}\right\}ds\\
			&= S_J\rho_1^{2/r}\left(t_j'^{[k]}\right)^{-1-\frac{2}{r}}
			\int_{t_p'^{[k]}/t_j'^{[k]}}^{1}
			\left(1-\widehat{v}\right)^{-\beta-1}\widehat{v}^{-\frac{2}{r}-\theta}
			\min\left\{1,\frac{\left(1-\widehat{v}\right)\left(t_j'^{[k]}\right)}{\rho_j'}\right\}
			d\widehat{v}\\
			&:= S_J\rho_1^{2/r}\left(t_j'^{[k]}\right)^{-1-\frac{2}{r}}I_3,
		\end{aligned}
	\end{equation}
	where the substitution
	$\widehat{v}=(s-(k-1)\tau)/(t_j'^{[k]})$ and $\beta+\theta=1$ are used.

	In order to estimate $I_3$, we divide its integration interval into
			$(t_p'^{[k]}/t_j'^{[k]},
			\max\{\frac{1}{2},t_p'^{[k]}/t_j'^{[k]}\})$,
            $
			(\max\{\frac{1}{2},t_p'^{[k]}/t_j'^{[k]}\},
            1-\rho_j'/t_j'^{[k]})$ and
			$(1-\rho_j'/t_j'^{[k]},1)$
	and denote $I_3=I_{31}+I_{32}+I_{33}$.
	Since
	$
		\rho_j'/t_j'^{[k]}
		\lesssim (\rho_1/t_j'^{[k]})^{1/r}
		\simeq (C(j-2(k-1)N))^{-1}
		\leq (C(p-2(k-1)N))^{-1}
	$
	for sufficiently large $p$, it is clear that
	%\begin{equation}\label{eq5_14}
    $
		\rho_j'/(t_j'^{[k]}) \leq \frac{1}{2}.
        $
	%\end{equation}

	For $I_{31}$, if $t_p'^{[k]}/t_j'^{[k]}\geq1/2$,
	then $I_{31}=0$. Otherwise,
	we have $\widehat{v}\leq1/2$, i.e.,
	$(1-\widehat{v})^{-\beta-1}<(1/2)^{-\beta-1}$.
	From $\rho_j'/t_j'^{[k]} \leq \frac{1}{2}$
    %\eqref{eq5_14}
    and $\beta=-\theta+1$, it readily follows that $(1-\widehat{v})t_j'^{[k]}/\rho_j'\geq 1$.
	Thus we have
	%\begin{equation}\label{eq5_15}
    $
			|I_{31}| \leq |\int_{t_p'^{[k]}/t_j'^{[k]}}^{1/2}
			\widehat{v}^{-\theta-\frac{2}{r}}\,d\widehat{v}|
			\lesssim (t_p'^{[k]}/t_j'^{[k]})^{\beta-\frac{2}{r}}.
            $
	For $I_{32}$, since $\widehat{v}\geq1/2$ and
	$\widehat{v}\leq 1-\rho_j'/t_j'^{[k]}$,
	we can deduce
 %    $\widehat{v}^{-\theta-2/r}\leq\left(1/2\right)^{-\theta-2/r}$ and
	% $1-\widehat{v}\geq\dfrac{\rho_j'}{t_j'-(k-1)\tau}$,
	% which implies
    $(1-\widehat{v})t_j'^{[k]}/\rho_j'\geq 1$.
	This gives
	%\begin{equation}
		%\label{eq5_16}
        $
			|I_{32}| \leq |\int_{1/2}^{1-\rho_j'/t_j'^{[k]}}
			(1-\widehat{v})^{-\beta-1}\,d\widehat{v}|
			\lesssim (\rho_j'/t_j'^{[k]})^{-\beta}.
        $
	%\end{equation}
	For $I_{33}$, in view of $\widehat{v}\geq 1-\rho_j'/t_j'^{[k]}$ and $\rho_j'/t_j'^{[k]} \leq \frac{1}{2}$,
    %\eqref{eq5_14},
    it is easy to know that $1-\widehat{v}\leq\rho_j'/t_j'^{[k]}\leq1/2$, i.e.,
	$(1-\widehat{v})t_j'^{[k]}/\rho_j'\leq 1$ and
	$\widehat{v}^{-\theta-2/r}\leq\left(1/2\right)^{-\theta-2/r}$.
	Thus, we have
	%\begin{equation}
		%\label{eq5_17}
        $
			|I_{33}| \leq |\int_{1-\rho_j'/t_j'^{[k]}}^{1}
			(1-\widehat{v})^{-\beta-1}
			(1-\widehat{v})t_j'^{[k]}/\rho_j'\,d\widehat{v}|
			\lesssim (\rho_j'/{t_j'^{[k]}})^{-\beta}.
            $
	%\end{equation}
	Now, by combining above estimates of $I_{31}$, $I_{32}$ and $I_{33}$, \eqref{eq5_13} can be rewritten as
	\begin{equation}
		\label{eq5_18}
		\begin{aligned}
			\left|\left.\left(D_t^\beta - D^\beta\right) C_L\left(t\right)\right|_{t=t_j'} \right|
			\lesssim &S_J\rho_1^{2/r}\left(t_j'^{[k]}\right)^{-1-\frac{2}{r}}
			\left[
			\left(\frac{t_p'^{[k]}}{t_j'^{[k]}}\right)^{\beta-\frac{2}{r}}
			+\left(\frac{\rho_1^{\frac{1}{r}}(t_j'^{[k]})^{1-\frac{1}{r}}}{t_j'^{[k]}}\right)^{-\beta}
			\right]\\
			=& S_J\left(t_j'^{[k]}\right)^{-1}
			\left(\frac{\rho_1}{t_p'^{[k]}}\right)^{\frac{2-\beta}{r}}
			\left[
			\left(\frac{t_p'^{[k]}}{t_j'^{[k]}}\right)^{\beta}
			\left(\frac{\rho_1}{t_p'^{[k]}}\right)^{\frac{\beta}{r}}
			+
			\left(\frac{t_p'^{[k]}}{t_j'^{[k]}}\right)^{\frac{2-\beta}{r}}
			\right]\\
			\lesssim& S_J\left(t_j'^{[k]}\right)^{-1}
			\left(p-2(k-1)N\right)^{\beta-2}
			\left(\frac{t_p'^{[k]}}{t_j'^{[k]}}\right)^{\gamma^*},
		\end{aligned}
	\end{equation}
	where
$\rho_j' \lesssim \rho_1^{1/r}\left(t_j'^{[k]}\right)^{1-1/r}$,
		$(\rho_1/(t_p'^{[k]}))^{(2-\beta)/r}
		\simeq \left[\left(p-2(k-1)N\right)\right]^{-(2-\beta)}$
        and
        $
		(t_p'^{[k]})/(t_j'^{[k]}) \leq 1$
    are used.
	Therefore, \eqref{eq4_21} holds.
\end{proof}

\subsection{Proof of Lemma \ref{lemma4.4}}
\begin{proof}
First, based on the submesh $\left\{t_j'\right\}$, we derive the following inequalities
\begin{align}
	&C_L\left(t_j'\right) \lesssim S_L \left(t_j'^{[k]}\right)^{-\theta},\label{eq5_19}\\
	&D^\beta C_L\left(t_j'\right) \gtrsim \left(t_j'^{[k]}\right)^{-1}
	\left( \frac{t_p'^{[k]}}{t_j'^{[k]}} \right)^{\gamma}.\label{eq5_20}
\end{align}
From the expression of $B(t,t_m')$, it is easy to see that $B(t_j',\, t_{pl}' + (k-1)\tau) \le (t_j'^{[k]})^{-\theta}$.
Therefore, the expression of $C_L\left(t_j'\right)$ with $\sum_{l=0}^{L}c_l\lesssim S_L$ clearly leads to \eqref{eq5_19}.
In order to establish \eqref{eq5_20}, we firstly consider the following estimate
\begin{equation}\label{eq5_21}
	\left| \left.\left(D_t^\beta - D^\beta\right) C_L\left.\left(t\right)\right|_{t=t_j'} \right|
	\le \frac{1}{2} D_t^\beta C_L\left(t\right)\right|_{t=t_j'}.
\end{equation}	
Recalling conclusion $(ii)$ of Lemma \ref{lemma4.2} yields
$D_t^\beta B\left.\left(t,\, t_{p_l}' + (k-1)\tau\right)\right|_{t=t_j'} \ge 0$.
Therefore, for each $0\le l\le L$: if
 $p_l + 2(k-1)N < j < p_{l+1} + 2(k-1)N$, then it follows from \eqref{eq4_16} and \eqref{eq4_18} that
\begin{equation}
	\label{eq5_22}
	\begin{aligned}
		D_t^\beta C_L\left.\left(t\right)\right|_{t=t_j'}
		&= D_t^\beta\left[ \sum_{l=0}^{L} c_l \left.B\left(t, t_{p_l}' + (k-1)\tau\right)\right|_{t=t_j'} \right]
		\ge c_l D_t^\beta \left.B\left(t, t_{p_l}' + (k-1)\tau\right)\right|_{t=t_j'} \\
		&\ge \left( \frac{t_p'^{[k]}}{t_j'^{[k]}} \right)^{\gamma}
		\frac{\theta}{\Gamma(\theta)}
		\left( \frac{t_{p_l}'}{t_j'^{[k]}} \right)^{1-\theta}
		\left(t_j'^{[k]}\right)^{-1}
		\gtrsim \left( \frac{t_p'^{[k]}}{t_j'^{[k]}} \right)^{\gamma}
		\left(t_j'^{[k]}\right)^{-1};
	\end{aligned}
\end{equation}
and if $1+2(k-1)N \le j \le p$,
then \eqref{eq4_17} implies
\begin{equation}
	\label{eq5_23}
	\begin{aligned}
		D_t^\beta \left.C_L\left(t\right)\right|_{t=t_j'}
		&\ge c_0 D_t^\beta \left.B\left(t, t_p'\right)\right|_{t=t_j'}
		\ge \frac{1}{\Gamma(\theta)\theta}
		\left(t_p'^{[k]}\right)^{-1}
		\left( \frac{t_j'^{[k]}}{t_p'^{[k]}} \right)^{\theta}
		\gtrsim \left(t_j'^{[k]}\right)^{-1}
		\left( \frac{t_p'^{[k]}}{t_j'^{[k]}} \right)^{\gamma}
		\left( \frac{t_j'^{[k]}}{t_p'^{[k]}} \right)^{\gamma+\theta+1} \\
		&\gtrsim \left(t_j'^{[k]}\right)^{-1}
		\left( \frac{t_p'^{[k]}}{t_j'^{[k]}} \right)^{\gamma}
		\left(p-2(k-1)N\right)^{-r\max\left\{0,\gamma+\theta+1\right\}},
	\end{aligned}
\end{equation}
where
$
	t_j'^{[k]}/t_p'^{[k]}
	\ge \rho_1/t_p'^{[k]}
	\simeq
	 \left[C\left(p-2(k-1)N\right)\right]^{-r}$ and $c_0=1$
 are used in \eqref{eq5_23}.
%if $\gamma+\theta+1 \ge 0$, then
%$\left( \frac{t_j'-(k-1)\tau}{t_p'-(k-1)\tau} \right)^{\gamma+\theta+1}
%\gtrsim\left(p-2(k-1)N\right)^{-r(\gamma+\theta+1)}$;
%and if $\gamma+\theta+1 < 0$, then
%$\left( \frac{t_j'-(k-1)\tau}{t_p'-(k-1)\tau} \right)^{\gamma+\theta+1}\ge
%\left( \frac{t_j'-(k-1)\tau}{t_p'-(k-1)\tau} \right)^0 = 1$.

When $t_j' \le t_p'$, we can clearly obtain
$|(D_t^\beta-D^\beta)C_L(t)|_{t=t_j'}| = 0$ by Lemma \ref{lemma4.3} and $D_t^\beta C_L(t)|_{t=t_j'}>0$ from \eqref{eq5_23}. Thus,  \eqref{eq5_21} holds.
When $t_j' > t_p'$, in order to prove \eqref{eq5_21}, comparing \eqref{eq4_21} with \eqref{eq5_22} implies that we suffice to show that
$
	S_J \left(t_j'^{[k]}\right)^{-1}
	\left(p-2(k-1)N\right)^{\beta-2}
	\left( t_p'^{[k]}/t_j'^{[k]}\right)^{\gamma^*}
	\lesssim\left( t_p'^{[k]}/t_j'^{[k]} \right)^{\gamma}
	\left(t_j'^{[k]}\right)^{-1}.
$
Since $p$ is a sufficiently large integer, it means that we suffice to prove
\begin{equation}\label{eq5_24}
	S_J \left( \frac{t_p'^{[k]}}{t_j'^{[k]}} \right)^{\gamma^*}
	\lesssim\left( \frac{t_p'^{[k]}}{t_j'^{[k]}} \right)^{\gamma},\  \gamma \le \gamma^*.
\end{equation}
If $0 < \gamma \le\gamma ^*$, then $t_p'^{[k]}/t_j'^{[k]} \le 1$ and $S_J=1$, which follows from  \eqref{eq4_19}, indicate that \eqref{eq5_24} holds;
and if $\gamma \le 0$, then, by recalling the definition of $J = J(j)$ which implies $p_J+2(k-1)N < j \le p_{J+1}+2(k-1)N$, we have
$
t_{p_J}'/t_p'^{[k]}\le t_j'^{[k]}/t_p'^{[k]}.
$
Furthermore, for \(\gamma=0\), we deduce
$S_J = 1 + \ln (t_{p_J}'/t_p'^{[k]} )
\le 1 + \ln(t_j'^{[k]}/t_p'^{[k]}),$
which suggests $
S_J (t_j'^{[k]}/t_p'^{[k]} )^{-\gamma^*} \lesssim 1=( t_p'^{[k]}/t_j'^{[k]})^{0};
$
for \(\gamma < 0\), since
$
t_p'^{[k]}/t_{p_J}'
\ge t_p'^{[k]}/t_j'^{[k]}
$
and $t_p'^{[k]}/t_j'^{[k]} \le 1$,
one has
$
S_J =  (t_p'^{[k]}/t_{p_J}' )^{\gamma}
\le  ( t_p'^{[k]}/t_j'^{[k]} )^{\gamma}.
$
Since $( t_p'^{[k]}/t_j'^{[k]} )^{\gamma^*} \le 1$, it is clear that \eqref{eq5_24} holds. It means that \eqref{eq5_21} is true.

According to \eqref{eq5_21}, we clearly have $D^\beta C_L(t_j') \ge \frac{1}{2} D_t^\beta C_L(t)|_{t=t_j'}$.
It follows from \eqref{eq5_22} and \eqref{eq5_23} that \eqref{eq5_20}  can be obtained obviously. By an analogous derivation to that of \eqref{eq5_19} and \eqref{eq5_20}, we can readily arrive at \eqref{eq4_22} and \eqref{eq4_23}.
\end{proof}

\subsection{Proof of Lemma \ref{lemma4.5}}
\begin{proof}
We divide this proof into two different cases: $\gamma \le \gamma^*$ and $\gamma > \gamma^*$. First, we consider the case $\gamma \le \gamma^*$.
From \eqref{eq4_23} and \eqref{eq4_24}, it immediately follows that
	\begin{equation}
		\label{eq5_25}
			D^\beta V^n
			\le \left(t_n^{[k]}\right)^{-1}
			\left( \frac{\rho_1}{t_n^{[k]}} \right)^\gamma
			\lesssim \left(t_n^{[k]}\right)^{-1}
			\left( \frac{\left(C\left(p-2(k-1)N\right)\right)^r\rho_1}{t_n^{[k]}} \right)^\gamma
			\simeq \left(t_n^{[k]}\right)^{-1}
			\left( \frac{t_p'^{[k]}}{t_n^{[k]}} \right)^\gamma
			\lesssim D^\beta C_L(t_n),
	\end{equation}
	where
	$t_p'^{[k]} \simeq ( \frac{C\left(p-2(k-1)N\right)}{N} )^r
	\simeq [C(p-2(k-1)N)]^r \rho_1$ is used.
	Since $V^0 = C_L(0) = \sum\limits_{l=0}^{L} C_l B\left(0,\, t_{p_l}'\right) = 0$,
	applying Lemma \ref{lemma4.1} and utilizing \eqref{eq4_22} gives
	\begin{equation}
		\label{eq5_26}
		V^n \le C_L(t_n) \lesssim \left(t_n^{[k]}\right)^{-\theta} S_L,
		\ n\leq2(k-1)N + p_{L+1}.
	\end{equation}	
	When $t_n \ge t_p'$, if $\gamma>0$, then $S_L=1$; if $\gamma= 0$, then $S_L=1+ \ln( t_{p_L}'/t_p'^{[k]})\le 1+ \ln(t_n^{[k]}/\rho_1)$
    %for $p_L + 2(k-1)N < n \le p_{L+1} + 2(k-1)N$
    ; and if $\gamma<0$, then
	$S_L = ( t_p'^{[k]}/t_{p_L}' )^\gamma
	\le ( \rho_1/t_n^{[k]})^\gamma$,
	where the fact
    %$\frac{t_p'-(k-1)\tau}{t_{p_L}'} \ge \frac{\rho_1}{t_n-(k-1)\tau}$, and when $p_L + 2(k-1)N < n < p_{L+1} + 2(k-1)N$, we have
    $t_{p_L}'/t_p'^{[k]}
	\le t_n^{[k]}/\rho_1$ for $p_L + 2(k-1)N < n < p_{L+1} + 2(k-1)N$ is used.
    Since the error constants in \eqref{eq4_22} and \eqref{eq4_23} are independent of $L$, it follows from \eqref{eq4_27} and \eqref{eq5_26} that
	  \eqref{eq4_26} holds for fixed $L>0$.
	When $t_n \le t_p'$, taking $L=0$ in \eqref{eq5_26} and noting that $S_0=1$, we obtain
	\begin{equation}
		\label{eq5_27}
		V^n \le C_0(t_n) \lesssim \left(t_n^{[k]}\right)^{-\theta} S_0
		= \left(t_n^{[k]}\right)^{-\theta}.
	\end{equation}
 If $\gamma>0$, then \eqref{eq4_26} follows directly from \eqref{eq5_27}.
	If $\gamma=0$, then we have
	$V^n \lesssim \left(t_n^{[k]}\right)^{-\theta}
	\le (t_n^{[k]})^{-\theta}
	( 1 + \ln(t_n^{[k]}/\rho_1 )) $.
	If $\gamma<0$, then
	$V^n \lesssim (t_n^{[k]})^{-\theta}
	\le (t_n^{[k]})^{-\theta}
	( \rho_1/t_n^{[k]})^\gamma$.
	According to \eqref{eq4_27}, it is easy to see that \eqref{eq4_26} also holds when $t_n \le t_p'$.
	
	Next, we consider the case $\gamma > \gamma^*$.
	Let $\left\{\overline{V}^n\right\}$ satisfy
	\begin{empheq}[left=\empheqlbrace]{align}
	&D^\beta \overline{V}^n = \left(t_n^{[k]}\right)^{-1}
	\left( \frac{\rho_1}{t_n^{[k]}} \right)^{\gamma^*},
		\label{eq5_28} \\
		&\overline{V}^{n-2N} = 0,\qquad 1\le n\le 2N.
		\label{eq5_29}
	\end{empheq}
	By applying the established \eqref{eq4_26} for $\gamma=\gamma^*$, so one has $\overline{V}^n \lesssim \left(t_n^{[k]}\right)^{\beta-1}$.
    In addition,
	combining \eqref{eq4_24} and \eqref{eq5_28} yields
	\begin{equation}\label{eq5_30}
	D^\beta V^n
	\le \left(t_n^{[k]}\right)^{-1}
	\left( \frac{\rho_1}{t_n^{[k]}} \right)^\gamma
	\le \left(t_n^{[k]}\right)^{-1}
	\left( \frac{\rho_1}{t_n^{[k]}} \right)^{\gamma^*}
	= D^\beta \overline{V}^n.
	\end{equation}
So, it follows from \eqref{eq4_25}, \eqref{eq5_29}, \eqref{eq5_30} and Lemma \ref{lemma4.1} that
	$V^n \le \overline{V}^n$. It indicates $V^n\lesssim \left(t_n-(k-1)\tau\right)^{\beta-1}$.
	Thus, \eqref{eq4_26} holds for $\gamma > \gamma^*$.
\end{proof}

\subsection{Proof of Lemma \ref{lemma4.6}}
\begin{proof}
	We first consider the case $t_n=0$. Let
	\begin{equation}\label{eq5_31}
	E(t) := \sum_{q=0}^Q \left(\bar{a}\right)^q E_q(t),
	E_q(t) = \max\left\{0,\, t-b_q\right\},
	b_0=0,b_q \in \left[a_0 q - \bar{\rho},\, a_0 q\right],
    \chi_I\left(t\right) =
	\begin{cases}
		1, & t\in I, \\
		0, & \text{otherwise}.
	\end{cases}
\end{equation}
%and
%\begin{equation}\label{eq5_32}
%	\chi_I\left(t\right) =
	%\begin{cases}
	%	1, & t\in I, \\
		%0, & \text{otherwise},
	%\end{cases}
%\end{equation}
	where $\bar{a}$ is a sufficiently large number. Moreover, we choose $0\le Q\lesssim1$ such that $K\tau \in (b_{Q+1},\, b_{Q+2}]$,
	i.e., $Q+2 = \lceil K\tau/a_0 \rceil$, unless $a_0\ge K\tau$, in which case we take $Q=0$.
	For $q=0$ in \eqref{eq5_31}, applying the continuous operator $(D_t^\beta-\lambda)$ to $E_0(t)=t$ yields
	$(D_t^\beta - \lambda)E_0(t)
	= t^{1-\beta}( 1/(\Gamma(2-\beta)) - \lambda t^\beta ).$
	If $t \le t^* := (\Gamma(2-\beta)\lambda)^{-1/\beta}$, then we have	
	$
	(D_t^\beta - \lambda)E_0(t) \ge 0;
	$
    and if $t > t^*$, then one has $(D_t^\beta - \lambda)E_0(t) < 0$.
    Therefore,
	the fact that  $t^* >5a_0/2 > b_2$ indicates
	$
	(D_t^\beta - \lambda)E_0(t) \ge 0
	$
    for $t\in(0,\, b_2]$,
	 and in view of $(D_t^\beta - \lambda)E_0(t) \ge 0$ for $t\in(b_2,\, t^*]$,
	we have
	$
	(D_t^\beta - \lambda)E_0(t)\gtrsim -1
	$ for $t\in(b_2,\, K\tau]$.
By \eqref{eq5_31}, it is clear that
	\begin{equation}\label{eq5_33}
	\left(D_t^\beta - \lambda\right)E_0(t) \gtrsim\chi_{\left[b_1,b_2\right]}(t)-\chi_{\left(b_2,K\tau\right]}(t).
	\end{equation}
	For each $q$: if $t\le b_q$, then $E_q(t)=0$; and if $t>b_q$, then $E_q(t)=t-b_q$. Hence, by a similar argument, applying $(D_t^\beta-\lambda)$ to $t-b_q$ for $t>b_q$ yields
	$
		(D_t^\beta - \lambda)\left(t-b_q\right)=(t-b_q)^{1-\beta}
		( 1/(\Gamma(2-\beta)) - \lambda\left(t-b_q\right)^\beta ) \gtrsim\chi_{[b_{q+1},b_{q+2}]}(t)-\chi_{(b_{q+2},K\tau]}(t),
	$
	which leads to
	\begin{equation}\label{eq5_34}
	(D_t^\beta - \lambda)E_q(t) \gtrsim\chi_{[b_{q+1},b_{q+2}]}(t)-\chi_{(b_{q+2},K\tau]}(t).
	\end{equation}
	Using \eqref{eq5_33} and \eqref{eq5_34}, and noting that $Q\lesssim 1$, $\bar{a}$ is sufficiently large, and $K\tau \in \left(b_{Q+1}, b_{Q+2}\right]$ implies $\left(b_{Q+2}, K\tau\right]=\emptyset$  for $q=Q$, we obtain
	\begin{equation}\label{eq5_35}
	\left(D_t^\beta - \lambda\right)E(t)=\sum_{q=0}^Q \left(\bar{a}\right)^q\left(D_t^\beta - \lambda\right) E_q(t)\gtrsim \chi_{\left(b_{1},K\tau\right]}(t).
	\end{equation}
	Choose $t_{b_q}$ to be the largest grid point in $\left[a_0 q - \bar{\rho},\, a_0 q\right]$, i.e., $\{b_q\}_{q=0}^{Q}\subset\{t_i\}_{i=0}^{2KN}$. Then we have
	\begin{equation}\label{eq5_36}
	\left(D_t^\beta - \lambda\right)E(t)\biggr|_{t=t_i} = \left(D^\beta - \lambda\right)E(t_i).
	\end{equation}
	Let $E^i= E(t_i)$.  From $E(t_i) = \sum_{q=0}^Q (\bar{a})^q \max\{0,\, t_i - b_q\}$ and $\max\{0,\, t_{i+1}-b_q\} \ge \max\{0,\, t_i-b_q\}$,
	it immediately follows that $E(t_{i+1}) \ge E(t_i)$, i.e., $\left\{E^i\right\}$ is a nondecreasing sequence.
	Therefore, for $t_n=0$, it follows from \eqref{eq5_35} and \eqref{eq5_36} that we have found a nondecreasing sequence $E^i=E(t_i)$ satisfying $0\le E^i\lesssim 1$ and
	$E^i=0,\ i=0,
	\left(D^\beta - \lambda\right)E^i \gtrsim
	\begin{cases}
		0, & t_i < a_0, \\
		1, & t_i \ge a_0,
	\end{cases},\quad1\le i\le 2KN.
	$
	
For the case $t_n>0$, we introduce the submesh $\{t_i\}$ and define
	$
	\tilde{t}_l = t_{n+l} - t_n,\ 0\leq l\leq 2KN-n,
	$
	$
	\tilde{\rho}_l = \tilde{t}_{l} - \tilde{t}_{l-1}
	= \rho_{n+l} \le \bar{\rho} \le a_0/2.
	$
	 Based on the above discussion on $t_n=0$, it is obvious that, for $\tilde{t}_l=0$, there exists a nondecreasing sequence
	$\left\{\tilde{E}^l\right\}$ satisfying $0\le \tilde{E}^l\lesssim 1$ and
\begin{equation}\label{eq5_37}
\tilde{E}^l=0,\quad l=0,\quad
	(D^\beta-\lambda)\tilde{E}^l\gtrsim
	\begin{cases}
		0,& \tilde{t_l}<a_0,\\
		1,& \tilde{t_l}\ge a_0,
	\end{cases}
	,\quad 1\le l\le 2KN-n.
\end{equation}
	Define
	$
		E^i =
		\begin{cases}
			0, & 0\le i< n, \\
			\tilde{E}^{i-n}, & n \leq i \le 2KN.
		\end{cases}
	$
	Therefore, for $t_n>0$,  \eqref{eq5_37} suggests that there exists a nondecreasing sequence
	$\left\{E^i\right\}$ such that Lemma 4.6 holds.
 %    $0\le E^i\lesssim 1$ and
	% $$
	% E^i=0,\quad \forall i\le n,\quad
	% \left(D^\beta - \lambda\right)E^i=
	% \begin{cases}
	% 	0,& i<n,\\
	% 	(D^\beta-\lambda)\tilde{E}^{\,i-n},& n\le i\le 2KN,
	% \end{cases}\gtrsim
	% \begin{cases}
	% 	0, & t_i < t_n + a_0, \\
	% 	1, & t_i \ge t_n + a_0,
	% \end{cases}
	% ,\quad 1\leq i\leq 2KN.
	% $$
\end{proof}
%\section{References}

  \end{document}